\documentclass[10pt,leqno,english]{amsart}
\usepackage{hyperref} 
\usepackage{amssymb} 
\usepackage[square,numbers]{natbib} 
\usepackage[all]{xy} 

\renewcommand{\thesection}{\arabic{section}}

\numberwithin{subsubsection}{section}

\renewcommand{\theenumi}{\arabic{enumi}}

\renewcommand{\theequation}{\arabic{equation}}
\numberwithin{equation}{subsubsection}

\theoremstyle{plain}

\newtheorem{mainthm}{Theorem}

\newtheorem{mainlemm}{Lemma}

\newtheorem*{externalthm}{Theorem}

\swapnumbers
\newtheorem{thm}[subsubsection]{Theorem}
\newtheorem{lemm}[subsubsection]{Lemma}
\newtheorem{prop}[subsubsection]{Proposition}

\newtheorem{fact}[subsubsection]{Fact}
\newtheorem{obsv}[subsubsection]{Observation}

\DeclareMathOperator{\Mor}{Mor} 
\DeclareMathOperator{\End}{\mathit{End}} 
\DeclareMathOperator{\Hom}{\mathit{Hom}} 
\DeclareMathOperator{\Ex}{Ex} 

\DeclareMathOperator{\unit}{\mathsf{1}} 

\DeclareMathOperator{\NN}{\mathbb{N}} 
\DeclareMathOperator{\kk}{\Bbbk} 
\DeclareMathOperator{\QQ}{\mathbb{Q}} 

\DeclareMathOperator{\AOp}{\mathsf{A}} 

\DeclareMathOperator{\POp}{\mathsf{P}}
\DeclareMathOperator{\QOp}{\mathsf{Q}}

\DeclareMathOperator{\Sym}{Sym} 
\DeclareMathOperator{\Free}{\mathit{F}} 
\DeclareMathOperator{\Env}{\mathit{E}} 

\DeclareMathOperator{\X}{\mathcal{X}} 
\DeclareMathOperator{\A}{\mathcal{A}} 
\DeclareMathOperator{\C}{\mathcal{C}} 

\DeclareMathOperator{\I}{\mathcal{I}} 
\DeclareMathOperator{\J}{\mathcal{J}} 
\DeclareMathOperator{\K}{\mathcal{K}} 

\DeclareMathOperator{\Top}{\mathcal{T}} 
\DeclareMathOperator{\Bij}{\mathcal{B}ij} 
\DeclareMathOperator{\Simp}{\mathcal{S}} 
\DeclareMathOperator{\Mod}{Mod} 

\DeclareMathOperator{\Prop}{\mathcal{P}} 
\DeclareMathOperator{\Graph}{\mathcal{G}} 

\DeclareMathOperator{\ICat}{\mathbb{I}} 
\DeclareMathOperator{\JCat}{\mathbb{J}}
\DeclareMathOperator{\VCat}{\mathbb{V}}
\DeclareMathOperator{\YCat}{\mathbb{Y}}

\DeclareMathOperator*{\bigtimes}{\times} 

\DeclareMathOperator{\BCat}{\mathbb{B}}
\DeclareMathOperator{\ACat}{\mathbb{A}}

\DeclareMathOperator{\id}{id}
\DeclareMathOperator{\coker}{coker} 
\DeclareMathOperator*{\colim}{colim} 

\DeclareMathOperator{\dg}{\mathit{dg}} 
\DeclareMathOperator{\s}{\mathit{s}} 
\DeclareMathOperator{\Char}{char} 

\DeclareMathOperator*{\fib}{\twoheadrightarrow} 
\DeclareMathOperator*{\cofib}{\rightarrowtail} 
\DeclareMathOperator{\trivialcofib}{{\displaystyle\cofib^{\sim}}} 

\DeclareMathOperator{\In}{\mathit{In}} 
\DeclareMathOperator{\Out}{\mathit{Out}} 

\DeclareMathOperator*{\SCopy}{\mathit{S}} 
\DeclareMathOperator*{\TzeroCopy}{\mathit{T}_0}
\DeclareMathOperator*{\ToneCopy}{\mathit{T}_1}

\title[Props in model categories and homotopy invariance of structures]{Props in model categories\\and\\homotopy invariance of structures}
\author{Benoit Fresse}
\date{15 December 2008 (corrections and updates on 7 January, 22 December 2009)}
\address{Laboratoire Painlev\'e\\
UMR 8524 de l'Universit\'e Lille 1 - Sciences et Technologies - et du CNRS\\
Cit\'e Scientifique -- B\^atiment M2\\
F-59655 Villeneuve d'Ascq C\'edex (France)}
\email{Benoit.Fresse@math.univ-lille1.fr}
\urladdr{http://math.univ-lille1.fr/\~{ }fresse}
\thanks{Research supported in part by ANR grant JCJC-06-0042 OBTH}

\begin{document}

\begin{abstract}
We prove that any category of props in a symmetric monoidal model category
inherits a model structure.
We devote an appendix, about half the size of the paper, to the proof of the model category axioms in a general setting.
We need the general argument to address the case of props in topological spaces and dg-modules over an arbitrary ring,
but we give a less technical proof which applies to the category of props in simplicial sets, simplicial modules,
and dg-modules over a ring of characteristic $0$.

We apply the model structure of props to the homotopical study of algebras over a prop.
Our goal is to prove that an object $X$
homotopy equivalent to an algebra $A$ over a cofibrant prop $\POp$
inherits a $\POp$-algebra structure so that $X$ defines a model of $A$
in the homotopy category of $\POp$-algebras.
In the differential graded context,
this result leads to a generalization of Kadeishvili's minimal model of $A_\infty$-algebras.
\end{abstract}

\maketitle

\begin{center}\emph{In honor of Tornike Kadeishvili}\end{center}

\tableofcontents

\part*{Foreword}

The notion of prop has been introduced by F. Adams and S. Mac Lane as a conceptual device
to handle the operations of homotopical structures.
The idea, explained in~\cite{MacLaneProps}, is to consider all natural operations $p: X^{\otimes m}\rightarrow X^{\otimes n}$
with $m$ inputs and $n$ outputs (where both $m$ and $n$ run over $\NN$)
that can be formed in a category of algebras
within an ambient category $\C$
equipped with a symmetric tensor product $\otimes: \C\times\C\rightarrow\C$.
A prop $\POp$ is the abstract structure formed by the double sequence of objects $\POp(m,n)$
whose elements represent such operations $p: X^{\otimes m}\rightarrow X^{\otimes n}$.
The category of algebras associated to a prop $\POp$
consists of pairs $(X,\phi)$,
where $X$ is an object of the base category $\C$
and~$\phi$ is a map which associates an actual morphim $p: X^{\otimes m}\rightarrow X^{\otimes n}$
to each abstract operation $p\in\POp(m,n)$.

The structure of a prop is necessary to model operations of bialgebras,
but many algebra structures
are fully determined by operations $p: X^{\otimes m}\rightarrow X$
with only $1$ output.
In this situation,
any operation $p: X^{\otimes m}\rightarrow X^{\otimes n}$
with $n$ outputs has a decomposition
\begin{equation*}
X^{\otimes m}\xrightarrow{\sigma^*} X^{\otimes m}
= X^{\otimes m_1}\otimes\dots\otimes X^{\otimes m_n}\xrightarrow{q_1\otimes\dots\otimes q_n} X\otimes\dots\otimes X
= X^{\otimes n}
\end{equation*}
so that each $q_i$, $i = 1,\dots,n$, is an operation with $1$ output
and where $\sigma^*$ is a tensor permutation.
The work of J. Boardman and R. Vogt~\cite{BoardmanVogtAnnouncement,BoardmanVogt},
has brought out applications of props of this form
in homotopy theory.
The work of P. May~\cite{May}
has highlighted the role of operations with a single output, which define the core of these prop structures,
and the structured object, named an operad by him,
formed by the collection of these operations $p: X^{\otimes m}\rightarrow X$.

\medskip
The general background for the definition of props and operads is given by the notion of symmetric monoidal category.
Briefly, a symmetric monoidal category is a category $\C$
equipped with a unit object $\unit\in\C$
and a tensor product $\otimes: \C\times\C\rightarrow\C$
satisfying natural symmetry, associativity, and unit constraints.
The examples considered in this paper include the category of differential graded modules over a ground ring $\kk$
together with the tensor structure derived from the tensor product of $\kk$-modules,
the categories of simplicial modules,
as well as the categories of simplicial sets and topological spaces whose tensor product
is defined by the cartesian product of these categories.
In all cases, the category of props is richer than the category of operads.
Note however that any category of algebras over a prop is identified with a category of algebras over an operad
when the tensor product is given by a categorical cartesian product.
Indeed, in that situation, any operation with multiple outputs $p: X^{\times m}\rightarrow X^{\times n}$
is fully determined by its components $p_i: X^{\times m}\rightarrow X$, $i = 1,\dots,n$,
which define operations with a single output.
Thus any algebra over a prop $\POp$ in a cartesian category $\C$
is equivalent to an algebra over an operad $\POp^{\times}$
associated to $\POp$.
But when a structure is naturally encoded by a prop
we will find more convenient to do constructions within the category of props, and not operads,
even in the context of cartesian categories.
For our purpose,
one might observe that the correspondence between props and operads does not fit the standard settings
for the application of homotopical algebra methods.

\medskip
Initially,
the objective of Boardman-Vogt~\cite{BoardmanVogtAnnouncement,BoardmanVogt},
and May~\cite{May}
was to extend results of M. Sugawara~\cite{Sugawara}
and J. Stasheff~\cite{Stasheff}
on single loop spaces to iterated loop spaces.
The structure of single loop spaces is modelled, according to Sugawara and Stasheff,
by the notion of $A_\infty$-space,
a strong homotopical version of the notion of associative monoid.
The structure of an $A_\infty$-space consists roughly of a product
and a full set of homotopies that make this product associative
in the strongest homotopical sense.
This structure is associated to an operad.
There is also a notion of $A_\infty$-equivalence used to represent chains of equivalences of~$A_\infty$-spaces.

The present paper is mainly concerned with Boardman-Vogt' approach
and with homotopy invariance results rather than applications to iterated loop spaces.
Let $\POp$ be an operad in topological spaces.
The main device of Boardman-Vogt' work is an operad $W(\POp)$
naturally associated to $\POp$.
They also define a notion of $W(\POp)$-equivalence
in order to model chains of homotopy equivalences of $W(\POp)$-algebras.
For the operad of associative monoids $\AOp$,
the operad~$W(\AOp)$ is identified with a cubical subdivision of the operad of $A_\infty$-spaces
and the notion of $W(\AOp)$-equivalence
is the same as the notion of $A_\infty$-equivalence.

Boardman-Vogt prove that certain homotopy invariance properties of $A_\infty$-structures
can be generalized to $W(\POp)$-algebras.
Among other things:

\begin{externalthm}[J. Boardman, R. Vogt~\cite{BoardmanVogt}]
Let $(Y,\psi)$ be a $W(\POp)$-algebra.
If we have a homotopy equivalence of topological spaces $f: X\xrightarrow{\sim} Y$,
then the action of~$W(\POp)$ on~$Y$ can be transported to~$X$
and the $W(\POp)$-algebra $(X,\phi)$ obtained by this process
is $W(\POp)$-equivalent to the original one $(Y,\psi)$.
\end{externalthm}

Boardman and Vogt have settled their constructions in the topological context,
but we have already mentioned that the notion of algebra over an operad, and the notion of algebra over a prop similarly,
make sense in any symmetric monoidal category.
In particular,
we have an analogue of the notion of $A_\infty$-object
in the context of differential graded algebra
(in the sequel, we will use the prefix dg as an abbreviation for differential graded).

Parallel to the study of homotopy invariant structures in topology,
techniques were being developed in homological algebra to transport perturbations of differentials
through chain-equivalences.
The main construction, called the basic perturbation lemma,
has been brought out in~\cite{RBrown}.
This homological perturbation theory has been applied to chain models of fibrations.
The most basic example is the reduced bar complex~$B(C^*(X))$ of a cochain algebra~$C\*(X)$
which models the cochain complex of the loop space $C^*(\Omega X)$
in the path-space fibration $\Omega X\rightarrow P X\rightarrow X$.

Recall that the reduced bar complex of an associative dg-algebra $A$
is the twisted chain complex $B(A) = (T^c(\Sigma A),\partial)$
formed by the tensor coalgebra on the suspension of $A$
together with a differential $\partial: T^c(\Sigma A)\rightarrow T^c(\Sigma A)$
determined by the product of $A$.
Suppose we have a dg-module $M$ connected to $A$
by a homotopy equivalence of dg-modules $f: M\xrightarrow{\sim} A$.
Perturbation techniques produce a new complex $(T^c(\Sigma M),\partial)$
connected to $B(A) = (T^c(\Sigma A),\partial)$
by a homotopy equivalence of dg-coalgebras $\phi: (T^c(\Sigma M),\partial)\xrightarrow{\sim} (T^c(\Sigma A),\partial)$.

Tornike Kadeishvili realized an intellectual breakthrough in~\cite{Kadeishvili1}
when he observed that:

\begin{externalthm}[T. Kadeishvili~\cite{Kadeishvili1}]\hspace*{2mm}
\begin{enumerate}
\item
A differential on a tensor coalgebra $T^c(\Sigma M)$,
where $M$ is any dg-module,
is equivalent to an $A_\infty$-algebra structure on $M$.
\item
A homotopy equivalence of dg-coalgebras $f: (T^c(\Sigma M),\partial)\xrightarrow{\sim}(T^c(\Sigma N),\partial)$
is equivalent to an $A_\infty$-equivalence between $M$ and $N$.
\end{enumerate}
\end{externalthm}

Thus the earlier perturbation techniques could be applied
to prove dg-analogues for $A_\infty$-algebras
of the homotopy invariance results of Boardman-Vogt.
First applications of this idea have been developed in a series of papers~\cite{Kadeishvili2,Kadeishvili3,Kadeishvili4} by Kadeishvili,
in a joint work~\cite{HuebschmannKadeishvili}
of J. Huebschmann and Kadeishvili
and in a series of papers~\cite{GugenheimStasheff,GugenheimLambeStasheff1,GugenheimLambeStasheff2}
by V. Gugenheim, L. Lambe and J. Stasheff.
One result of Kadeishvili asserts that:

\begin{externalthm}[T. Kadeishvili~\cite{Kadeishvili1}, see also J. Huebschmann and T. Kadeishvili~\cite{HuebschmannKadeishvili}]
Let $A$ be a dg-algebra whose homology $H_*(A)$ is projective as a $\kk$-module.
The homology $H_*(A)$ carries an $A_\infty$-structure
such that $H_*(A)$ forms an $A_\infty$-algebra
equivalent to $A$.
\end{externalthm}

He calls this $A_\infty$-structure
the minimal model of~$A$.

\medskip
Though operads were introduced in topology in the late 60's,
the interest in structures defined by operads in algebra
has arisen in the mid 90's only.
Since then,
several variants of the notion of operad and prop (cyclic operads~\cite{GetzlerKapranovCyclic},
modular operads~\cite{GetzlerKapranovModular},
properads~\cite{ValletteThesis,Vallette}, wheeled props~\cite{Merkulov})
have been introduced.
The interest in structures defined by props in algebra has similarly grown up since the mid 90's
with new motivations coming from quantum field theory.

The dg-analogue of Boardman-Vogt homotopy invariance theorems has been established by M. Markl for algebras over operads in~\cite{Markl}.
The case of algebras over modular operads is addressed by J. Chuang and A. Lazarev
in~\cite{ChuangLazarevDualFeynman,ChuangLazarevMinimalModels}.
The case of algebras over properads can be deduced from techniques of S. Merkulov and B. Vallette~\cite{MerkulovVallette}.
Markl uses arguments close to Boardman-Vogt' topological constructions.
Merkulov-Vallette prove that the structure of an algebra over a properad
is equivalent to a solution of a Maurer-Cartan equation in an $L_\infty$-algebra.
The invariance of such solutions under $L_\infty$-equivalences
can be used to solve the transfer problem for algebras over a properad.
The articles~\cite{ChuangLazarevDualFeynman,ChuangLazarevMinimalModels} by Chuang-Lazarev
gives an elegant construction of minimal models (generalizing Kadeishvili's minimal models)
for algebras over (modular) operads.
Chuang-Lazarev also use the representation of an algebra structure
by a solution of Maurer-Cartan' equation.

The article~\cite{Hess} of K. Hess extends perturbation techniques to polygebras,
a relaxed version of the structure of algebra over a prop.

\medskip
We now have different theories to transfer structures within the category of topological spaces
and within the category of dg-modules.
A natural problem is to give universal arguments
which apply to any of these contexts,
and to the general case of algebras over a prop.
A nice axiomatic background
to do homotopical algebra
is given by the notion of model category.
Applications of model categories have been extensively studied
in the context of operads (see for instance~\cite{BergerMoerdijk,HinichHomotopy,Spitzweck}),
but not in the context of props yet.
It has been observed that the homotopy versions of usual algebraic structures,
like $A_\infty$-algebras,
are associated to cofibrant replacements in the category of operads.
In the light of model categories,
the homotopy invariance of algebraic structures can be formulated
as an equivalence of mapping spaces in the category of operads (see~\cite{Rezk})
and such equivalences can be used to transport structures in the context of operads (see~\cite{BergerMoerdijk,Rezk}).

\part*{Introduction}

The purpose of this paper is to study the homotopy of algebras over a prop from the viewpoint of model categories.

Our first goal is to prove that categories of props inherit a model structure.
A big part of the paper is devoted to the setting up of a general argument which applies to any category of props in a symmetric monoidal model category.
We need this general argument to include the case of props in topological spaces and the case of props in dg-modules over an arbitrary ring
within the scope of our methods.
We give another less technical argument which applies to the category of props in simplicial sets, simplicial modules,
and dg-modules over a ring of characteristic $0$.

We also aim at giving a conceptual proof of the existence of transferred structures in the prop context.
The arguments work in any ambient symmetric monoidal model category $\C$
provided certain fibrations are stable under tensor products (we say that $\C$ satisfies the monoid limit axioms,
see~\S\ref{PropAlgebras:LimitMonoidAxioms} for an explicit statement).
This includes: the category of topological spaces,
the category of simplicial sets,
the category of dg-modules over a ring,
the category of simplicial modules over a ring,
and other categories derived from these basic examples.
The additional limit monoid axiom
is not necessary when one deals with operads.

Our results give:

\begin{mainthm}\label{Mainresult:Transfer}
Let $\C$ be a cofibrantly generated symmetric monoidal model category
whose tensor product satisfies the limit monoid axioms.
Let $\POp$
be a cofibrant prop in $\C$.

Let $(Y,\psi)$ denote a $\POp$-algebra in $\C$.
The notation $Y$ refers to the underlying object in $\C$ of this $\POp$-algebra
and $\psi$ refers to the $\POp$-action on it.
Suppose we have a weak equivalence in $\C$
\begin{equation*}
f: X\xrightarrow{\sim} Y,
\end{equation*}
such that $X$, like $Y$, is both cofibrant and fibrant as an object of~$\C$.

The object $X$ can be equipped with a $\POp$-action so that we have a $\POp$-algebra $(X,\phi)$
connected to the original one $(Y,\psi)$
by weak equivalences of $\POp$-algebras:
\begin{equation*}
(X,\phi)\xrightarrow{\sim}\cdot\xleftarrow{\sim}\cdot\xrightarrow{\sim}(Y,\psi).
\end{equation*}
\end{mainthm}

As a corollary,
we obtain:

\begin{mainthm}\label{Mainresult:MinimalModel}
Let $\C$ be the category of nonnegatively lower graded differential modules over a ring $\kk$.
Let $(A,\phi)$ be an algebra in $\C$ over a cofibrant prop $\POp$.
If both the $\POp$-algebra $A$ and its homology $H_*(A)$ are projective as $\kk$-modules,
then $H_*(A)$ inherits a $\POp$-action so that we have a $\POp$-algebra $(H_*(A),\theta)$
connected to the original one $(A,\phi)$
by weak equivalences of $\POp$-algebras:
\begin{equation*}
(H_*(A),\theta)\xrightarrow{\sim}\cdot\xleftarrow{\sim}\cdot\xrightarrow{\sim}(A,\phi).
\end{equation*}
\end{mainthm}

\begin{proof}
Pick simply representatives of homology classes to form a quasiisomorphism $\gamma: H_*(A)\xrightarrow{\sim} A$
and apply Theorem~\ref{Mainresult:Transfer}.
\end{proof}

This theorem gives a generalization of Kadeishvili's minimal model
in the context of algebras over a prop.

\medskip
In these statement,
a weak equivalence of $\POp$-algebras refers obviously to a morphism of $\POp$-algebras
which forms a weak equivalence in the underlying category.
Other notions of homotopy equivalences
between algebras over a prop
have been defined in the literature.
In a sense,
we prove that a refined notion of homotopy equivalence of $\POp$-algebras
is equivalent to the naive one.

To be explicit,
recall that a $\POp$-action on an object $X\in\C$
is equivalent to a morphism of props $\phi: \POp\rightarrow\End_X$,
where $\End_X$ is a universal prop acting on $X$,
the endomorphism prop of~$X$.
For $\POp$-algebras $(X,\phi)$ and $(X,\psi)$
with the same underlying object $X\in\C$,
a notion of homotopy equivalence between $(X,\phi)$ and $(X,\psi)$
is defined by setting that $\phi,\psi: \POp\rightarrow\End_X$
are homotopic as prop morphisms.
This notion of homotopy equivalence is close to the one used by Chuang-Lazarev in the context of operads in dg-modules~\cite{ChuangLazarevMinimalModels},
and to Boardman-Vogt' notion of homotopy $W(\POp)$-equivalence~\cite{BoardmanVogt}.

The equivalence part of Theorem~\ref{Mainresult:Transfer}
arises as a consequence of the following statement:

\begin{mainthm}\label{Mainresult:Equivalence}
Let $\C$ be a cofibrantly generated symmetric monoidal model category
whose tensor product preserves fibrations.
Let $\POp$
be a cofibrant prop in $\C$.
Let $X\in\C$.
Suppose $X$ is both cofibrant and fibrant as an object of~$\C$.

The $\POp$-algebras $(X,\phi)$ and $(X,\psi)$
associated to homotopic prop morphisms $\phi,\psi: \POp\rightarrow\End_X$
are connected by weak equivalences of $\POp$-algebras:
\begin{equation*}
(X,\phi)\xrightarrow{\sim}\cdot\xleftarrow{\sim}\cdot\xrightarrow{\sim}(X,\psi).
\end{equation*}
\end{mainthm}

In the case where $\POp$ is an operad,
this result is a corollary of a more precise theorem of~\cite{Rezk}.
In the context of operads in dg-modules,
the theorem can be deduced from an extension of the bar duality of operads~\cite{FresseCylinder}.
The results of~\cite[\S IV.6]{BoardmanVogt}
also give an analogue of Theorem~\ref{Mainresult:Equivalence}
in Boardman-Vogt' setting (thus, in the context of operads in topological spaces).

\medskip
Since the publication of this paper, a further study of applications of homotopical algebra to props has appeared in~\cite{JohnsonYau}.
The authors of this new article,
M. Johnson and D. Yau,
feature the use of colored props in view towards applications in quantum field theory.

\subsection*{Overall plan}

In the first part of the paper, \emph{Homotopical algebra of props},
we study the homotopy of the category of props itself.
Our objective is to provide the category of props
with a model structure.
This part includes a review of definitions.

The second part of the paper, \emph{Homotopy invariance of structures},
is devoted to the proof of Theorem~\ref{Mainresult:Transfer}
and Theorem~\ref{Mainresult:Equivalence}.

The appendix is devoted to technical verifications required by the proof of the axioms of model categories for the category of props.
These verifications can be avoided if we deal with props in dg-modules over a field of characteristic $0$,
with props in simplicial modules over a ring or with props in simplicial sets.

\subsection*{Acknowledgements}
This research has been discussed at several internal seminars (Le Wast, La Wangenbourg)
of the project ``ANR JCJC06 OBTH''.
I am grateful to the participants of these seminars, David Chataur, Thomas Gire, Eric Hoffbeck, Joan Mill\`es and Bruno Vallette,
for many intensive and stimulating discussions on the matter of this article.
I am also grateful to Mark Johnson and Donald Yau for their interest in this work.
I thank Teimuraz Pirashvili, the referee, and the editorial board of the Georgian Mathematical Journal
for their careful reading of my manuscript and accurate remarks on the article content.

It is a great pleasure to express my intellectual debts
towards Tornike\break Kadeishvili.
My current work on the bar complex~\cite{FresseIteratedBar,FresseEinfinityBar}
has mainly arisen from stimulating discussions with Tornike Kadeishvili and Francis Sergeraert:
I have started this research in 2002,
at a workshop at the Razmadze Mathemathical Institute, Tbilissi,
which I visited at the invitation of Tornike Kadeishvili and Samson Saneblidze.
I thank Tornike Kadeishvili, Samson Saneblidze, and the other participants of the workshop
for this travel which has been so decisive for my own researches.

\part{Homotopical algebra of props}\label{PropHomotopy}
The main purpose of this part is to prove that props in a cofibrantly generated symmetric monoidal model category
inherit a model structure.
This result allows us to include the homotopical study of props in dg-modules, in topological spaces, in simplicial sets, \dots
in the same axiomatic framework.

To begin with, in~\S\ref{PropBackground},
we recall the definition of a prop in a symmetric monoidal category,
at least to fix conventions.
In~\S\ref{HomotopyBackground},
we review the definition of a cofibrantly generated symmetric monoidal model category
to fix the background of our constructions.

The axioms of model categories are fully satisfied by the category of props
in dg-modules over a ring of characteristic zero,
in topological spaces,
in simplicial sets,
but are not fully satisfied by the category of props in dg-modules over a field of positive characteristic
(the same observation holds in the context of operads).
In general,
we only have a semi-model structure in the sense of~\cite{HoveySemiModel},
but this structure is sufficient to do homotopical algebra.
The notion of semi-model category is reviewed in~\S\ref{SemiModelCategories}.
The semi-model structure of the category of props is defined in~\S\ref{PropSemiModel}.
The proof of the axioms reduces to a technical verification
which is postponed to the appendix.
The particular case of props in simplicial sets, topological spaces and dg-modules over a ring of characteristic zero
is studied in~\S\ref{PropPathObjectArgument}.
An alternate and less technical proof of the axioms, which apply to these particular instances of base categories, is given.

\section{Props in symmetric monoidal categories}\label{PropBackground}
The first purpose of this section is to review the definition of a prop
in a symmetric monoidal category.

The category of props comes equipped with free objects.
Free props have an explicit construction which is reviewed in Appendix~\ref{Graphs},
because we need some insight into the structure of free objects
in order to prove the axioms of semi-model categories.
For the moment,
we only recall the abstract definition, which suffices for our immediate needs.

In the appendix,
we review briefly a graphical representation
which helps to make an intuition of the structure of a prop.
The reader who is not aware of this representation can glance at this recollection in~\S\ref{Graphs}
after reading this section.

\subsubsection{Symmetric monoidal categories}\label{PropBackground:SymmetricMonoidalCategories}
The notion of symmetric monoidal category gives the categorical structure
in which props and algebras over props can be defined.

Recall briefly that a symmetric monoidal category is a category $\C$
equipped with a tensor product $\otimes: \C\times\C\rightarrow\C$
that satisfies unit, associativity and symmetry relations (see~\cite{MacLaneCategories}).
The notation $\unit$
is used to refer to the unit object of an abstract symmetric monoidal category $\C$.
The category of topological spaces $\C = \Top$
and the category of simplicial sets $\C = \Simp$
form symmetric monoidal categories
with respect to the cartesian product $\otimes := \times$.
The category of $\kk$-modules $\C = \kk\Mod$
forms a symmetric monoidal category
with respect to the tensor product over the ground ring $\otimes = \otimes_{\kk}$.
The category of dg-modules $\C = \dg\kk\Mod$
and the category of simplicial modules $\C = \s\kk\Mod$
inherit a tensor product from $\kk$-modules
and form symmetric monoidal categories as well.

The notation $\C$ is adopted throughout the paper to refer to a base symmetric monoidal category.
We tacitly assume that every small colimit, every small limit, exists in $\C$.
We assume moreover that the tensor product preserves colimits on both sides
and $\C$ comes equipped with a homomorphism bifunctor $\Hom_{\C}: \C^{op}\times\C\rightarrow\C$
such that
\begin{equation*}
\Mor_{\C}(A\otimes B,C) = \Mor_{\C}(A,\Hom_{\C}(B,C)).
\end{equation*}
This assumption is fulfilled by the category of topological spaces $\C = \Top$
(provided that we take a good notion of topological space, for instance compactly generated spaces),
simplicial sets $\C = \Simp$, by categories of modules $\C = \kk\Mod$,
dg-modules $\C = \dg\kk\Mod$ and simplicial modules $\C = \s\kk\Mod$.

In the context of an abstract base category~$\C$,
we use the additive notation $0$ (respectively, $\oplus$) to refer to the initial object (respectively, coproduct) of~$\C$,
but we do not necessarily assume that $0$ is a zero object,
nor that $\oplus$ is a biproduct.
For the product we use the standard notation $\times$.
The notation $*$ will refer to the final object of~$\C$.

\subsubsection{Diagram categories}\label{PropBackground:DiagramCategories}
In the paper we use several categories of diagrams in $\C$,
which are simply functors $K: \ICat\rightarrow\C$
on a fixed small category $\ICat$.
The category of $\ICat$-diagrams is denoted by $\C^{\ICat}$.

The category of $\ICat$-diagrams is tensored and enriched over~$\C$.
The tensor product $\otimes: \C\otimes\C^{\ICat}\rightarrow\C^{\ICat}$
is defined by the obvious pointwise formula $(C\otimes K)(\alpha) = C\otimes K(\alpha)$,
for any $C\in\C$ and $K\in\C^{\ICat}$.
The homomorphism bifunctor $\Hom_{\C^{\ICat}}(-,-): (\C^{\ICat})^{op}\times(\C^{\ICat})\rightarrow\C$
is given by the usual end formula
\begin{equation*}
\Hom_{\C^{\ICat}}(K,L) = \int_{\alpha\in\ICat}\Hom_{\C}(K(\alpha),L(\alpha))
\end{equation*}
for any $K,L\in\C^{\ICat}$.
The adjunction relation
\begin{equation*}
\Mor_{\C^{\ICat}}(C\otimes K,L)\simeq\Mor_{\C}(C,\Hom_{\C^{\ICat}}(K,L)),
\end{equation*}
where $C\in\C$ and $K,L\in\C^{\ICat}$,
is an immediate consequence of the definition of an end.

Any object $\alpha\in\ICat$
determines an enriched Yoneda functor $G_{\alpha}: \ICat\rightarrow\C$
such that
\begin{equation*}
\Hom_{\C^{\ICat}}(G_{\alpha},K)\simeq K(\alpha),
\end{equation*}
for any $K\in\C^{\ICat}$.
Let $\unit[S]$ denote the sum over a set $S$ of copies of the unit object $\unit\in\C$.
The image of~$\beta\in\ICat$
under $G_{\alpha}$ is given by:
\begin{equation*}
G_{\alpha}(\beta) = \unit[\Mor_{\ICat}(\alpha,\beta)].
\end{equation*}

For any functor $\phi: \ICat\rightarrow\JCat$,
we have an obvious restriction functor $\phi^*: \C^{\JCat}\rightarrow\C^{\ICat}$
naturally associated to $\phi$.
This functor has a left adjoint $\phi_!: \C^{\ICat}\rightarrow\C^{\JCat}$
given by the usual coend formula
\begin{equation*}
(\phi_! K)(\beta) = \int^{\alpha\in\ICat} \Mor_{\JCat}(\phi(\alpha),\beta)\otimes K(\alpha),
\end{equation*}
where $S\otimes C$
is another notation for the sum over a set $S$
of copies of an object $C\in\C$.

\subsubsection{The notion of prop}\label{PropBackground:PropDefinition}
A prop in a symmetric monoidal category $\C$
is a symmetric monoidal category $\POp$,
enriched over $\C$,
which has the nonnegative integers $n\in\NN$ as objects
and whose tensor product is given by the addition law $m\otimes n = m + n$
on objects $(m,n)\in\NN^2$.
From this definition,
it appears that the structure of a prop is fully determined by the structure
of the double sequence of hom-objects $\POp(m,n)\in\C$, $(m,n)\in\NN^2$.
Each hom-object $\POp(m,n)$ comes equipped with
\begin{itemize}
\item
a right action of the symmetric group $\Sigma_m$,
which reflects the action of $\Sigma_m$ on $m = 1^{\otimes m}$
at the homomorphism level,
\item
and a left action of the symmetric group $\Sigma_n$,
which reflects the action of $\Sigma_n$
on $n = 1^{\otimes n}$.
\end{itemize}
The composition structures of homomorphisms
consist of:
\begin{itemize}
\item
horizontal products
\begin{equation*}
\circ_h: \POp(k,m)\otimes\POp(l,n)\rightarrow\POp(k+l,m+n),\quad\text{for $k,l,m,n\in\NN$},
\end{equation*}
which define the tensor product of homomorphisms,
\item
vertical composition products
\begin{equation*}
\circ_v: \POp(k,m)\otimes\POp(n,k)\rightarrow\POp(n,m)\quad\text{for $n,k,m\in\NN$},
\end{equation*}
which define the composition of homomorphisms in the enriched category~$\POp$,
\item
and units
\begin{equation*}
\eta: \unit\rightarrow\POp(n,n),
\end{equation*}
which represent the identity morphisms of the objects $n\in\NN$
in $\POp$.
\end{itemize}
These operations are assumed to satisfy relations which reflect the axioms of symmetric monoidal categories
(see~\cite{EnriquezEtingof,MacLaneProps} for a comprehensive description of the relations).
In standard references,
the action of symmetric groups is encoded in extensions $\eta: \unit[\Sigma_n]\rightarrow\POp(n,n)$
of the unit morphism $\eta: \unit\rightarrow\POp(n,n)$.
Our presentation is obviously equivalent.
The naming of the horizontal and vertical products in a prop is motivated by the graphical representation of~\S\ref{Graphs}.

Naturally,
a morphism of props consists of a symmetric monoidal functor $\phi: \POp\rightarrow\QOp$,
which is fully determined by a collection of morphisms $\phi: \POp(m,n)\rightarrow\QOp(m,n)$
preserving structures.
The category of props is denoted by~$\Prop$.
The definition of limits and colimits in the category of props is postponed to~\S\ref{PropSemiModel}.
We will adopt the base-set notation $\vee$ to represent coproducts and pushouts in~$\Prop$,
the standard notation $\times$ to represent products and pullbacks.
The initial object of the category of props is denoted by $I$
and the final object by $*$.

\subsubsection{The diagram categories underlying props}\label{PropBackground:DiagramUnderlyingProps}
In the paper,
we use two categories of diagrams underlying the structure of a prop:
a first one, defined by forgetting multiplicative structures,
and another one,
defined by forgetting multiplicative and symmetric structures.
Formally,
these categories are defined as follows.

Let $\ACat = \NN\times\NN$ be the discrete category of pairs $(m,n)\in\NN\times\NN$.
Let $\BCat$ be the category formed by pairs $(m,n)\in\NN\times\NN$ as objects
together with morphism sets such that:
\begin{equation*}
\Mor_{\BCat}((m,n),(p,q)) = \begin{cases} \Sigma_m^{op}\times\Sigma_n, & \text{if $m = p$ and $n = q$}, \\
\emptyset, & \text{otherwise}. \end{cases}
\end{equation*}

Form the categories of diagrams associated to these small categories.
The $\ACat$-diagrams are just double sequences of objects $K(m,n)\in\C$.
The $\BCat$-diagrams, more usually called $\Sigma_*$-biobjects,
consist of collections
\begin{equation*}
M = \{M(m,n)\in\C\}_{(m,n)\in\NN\times\NN}
\end{equation*}
whose term $M(m,n)$ is equipped with a right $\Sigma_m$-action
and a left $\Sigma_n$-action (commuting with the right $\Sigma_m$-action).

There is an obvious functor $\phi: \ACat\rightarrow\BCat$
which reduces to the identity map on objects.
The associated restriction functor $\phi^*: \C^{\BCat}\rightarrow\C^{\ACat}$
simply forgets symmetric group actions.
The extension functor $\phi_!: \C^{\ACat}\rightarrow\C^{\BCat}$
maps an $\ACat$-diagram $K\in\C^{\ACat}$
to the $\Sigma_*$-biobject
such that
\begin{equation*}
(\phi_! K)(m,n) = \unit[\Sigma_n\times\Sigma_m^{op}]\otimes K(m,n),
\end{equation*}
for $(m,n)\in\NN\times\NN$.

The enriched Yoneda diagram $G_{(m,n)}\in\C^{\ACat}$
associated to a pair $(m,n)\in\NN\times\NN$
is given by the simple formula:
\begin{equation*}
G_{(m,n)}(p,q) = \begin{cases} \unit, & \text{if $(p,q) = (m,n)$}, \\ 0, & \text{otherwise}. \end{cases}
\end{equation*}
To avoid confusion,
we only use the notation $G_{(m,n)}$ to refer to the Yoneda diagram
associated to $(m,n)\in\NN\times\NN$
in $\C^{\ACat}$.
The Yoneda diagram associated to $(m,n)$ in $\C^{\BCat}$
can be identified with the image of $G_{(m,n)}\in\C^{\ACat}$
under the extension functor $\phi_!: \C^{\ACat}\rightarrow\C^{\BCat}$.
For this reason,
we adopt the notation $\phi_! G_{(m,n)}$
to refer to this diagram unambiguously.

\subsubsection{Forgetful functors and free props}\label{PropBackground:FreeProps}
The category of props is endowed with an obvious forgetful functor $U: \Prop\rightarrow\C^{\BCat}$.
This forgetful functor has a left adjoint $\Free: \C^{\BCat}\rightarrow\Prop$
which maps a $\Sigma_*$-biobject $M\in\C^{\BCat}$
to a corresponding free object $\Free(M)\in\Prop$.

To simplify notation,
the forgetful functor $U$ can be omitted
if the context makes clear that we deal with the $\Sigma_*$-biobject underlying a prop.

The free prop $\Free(M)\in\Prop$ associated to a $\Sigma_*$-biobject $M$
can be equivalently defined by a standard universal property.
The free prop $\Free(M)\in\Prop$
is endowed with a canonical morphism of $\Sigma_*$-biobjects $\eta: M\rightarrow\Free(M)$,
which represents the unit of the adjunction $\Free: \C^{\BCat}\rightleftarrows\Prop :U$.
Any morphism $\phi: M\rightarrow\QOp$,
where $\QOp$ is a prop,
has a unique factorization
\begin{equation*}
\xymatrix{ M\ar[rr]^{\phi}\ar[dr]_{\eta} && \QOp \\ & \Free(M)\ar@{.>}[ur]_{\exists!\tilde{\phi}} & }
\end{equation*}
such that $\tilde{\phi}: \Free(M)\rightarrow\QOp$ is a prop morphism.

The existence of free props follows easily from the adjoint functor theorem (as long as mild set-theoretic assumptions hold in the base category $\C$).
In~\S\ref{Graphs},
we review the effective construction of~\cite{EnriquezEtingof,ValletteThesis}
based on the graphical representation of props.

\section{Homotopical background}\label{HomotopyBackground}
The main purpose of this section is to specify the assumptions that we require on our base category~$\C$
to do homotopical algebra.
To summarize:
we assume that $\C$ is a cofibrantly generated closed symmetric monoidal model category.

Since we use a refinement of the notion of model category in the next section,
we prefer to review completely the definition of a model category, of a cofibrantly generated model category,
and the definition of a symmetric monoidal model category.
We refer to the books~\cite{Hirschhorn,Hovey}
for more details on the recollections of~\S\S\ref{HomotopyBackground:BasicDefinitions}-\ref{HomotopyBackground:CofibrantlyGenerated}.

\subsubsection{Basic definitions}\label{HomotopyBackground:BasicDefinitions}
We adopt the usual terminologies of homotopical algebra.
Recall briefly that a morphism $i: A\rightarrow B$
satisfies the left lifting property with respect to another morphism $p: X\rightarrow Y$,
and $p$ satisfies the right lifting property with respect to $i$,
if the fill-in morphism
exists in any diagram of the form
\begin{equation*}
\xymatrix{ A\ar[r]\ar[d]_{i} & X\ar[d]^{p} \\ B\ar[r]\ar@{.>}[ur] & Y }.
\end{equation*}
Recall also that a morphism $f: A\rightarrow B$
is a retract of another morphism $g: X\rightarrow Y$
if $f$ and $g$ fit in a commutative diagram
\begin{equation*}
\xymatrix{ A\ar[r]^{s}\ar[d]_{f} & X\ar[d]_{g}\ar[r]^{p} & A\ar[d]^{f} \\ B\ar[r]_{t} & Y\ar[r]_{q} & B }
\end{equation*}
such that $p s = \id$ and $q t = \id$.

\subsubsection{The axioms of model categories}\label{HomotopyBackground:ModelCategoryAxioms}
A model category is a category $\C$
equipped with three classes of distinguished morphisms,
called weak equivalences ($\xrightarrow{\sim}$),
fibrations ($\twoheadrightarrow$)
and cofibrations ($\rightarrowtail$),
so that the following axioms hold:
\begin{enumerate}\renewcommand{\theenumi}{M\arabic{enumi}}
\item\label{M:Completeness}\emph{completeness axiom}:
Limits and colimits exist in $\C$.
\item\label{M:TwoOutofThree}\emph{two-out-of-three axiom}:
Let $f$ and $g$ be composable morphisms.
If any two among $f$, $g$ and $f g$ are weak equivalences,
then so is the third.
\item\label{M:Retract}\emph{retract axiom}:
Suppose $f$ is a retract of $g$.
If $g$ is a weak equivalence (respectively a cofibration, a fibration),
then so is $f$.
\item\label{M:Lifting}\emph{lifting axioms}:
\begin{enumerate}
\item
The cofibrations have the left lifting property with respect to acyclic fibrations,
where an \emph{acyclic fibration}
refers to a morphism which is both a weak equivalence and a fibration.
\item
The fibrations have the right lifting property with respect to acyclic cofibrations,
where an \emph{acyclic cofibration}
refers to a morphism which is both a weak equivalence and a cofibration.
\end{enumerate}
\item\label{M:Factorization} (\emph{factorization axioms}):
\begin{enumerate}
\item
Any morphism has a factorization $f = p i$ such that $i$ is a cofibration and $p$ is an acyclic fibration.
\item
Any morphism has a factorization $f = q j$ such that $j$ is an acyclic cofibration and $q$ is a fibration.
\end{enumerate}
\end{enumerate}
By convention,
an object $X$ is \emph{cofibrant} if the initial morphism $0\rightarrow X$
is a cofibration,
\emph{fibrant} if the terminal morphism $X\rightarrow *$
is a fibration.

Recall that the class of fibrations (respectively, the class of acyclic fibrations)
in a model category is fully characterized by the right lifting property
with respect to acyclic cofibrations (respectively, cofibrations).
Dually,
the class of cofibrations (respectively, the class of acyclic cofibrations)
is fully characterized by the left lifting property
with respect to acyclic fibrations (respectively, fibrations).

\subsubsection{Homotopy classes in model categories}\label{HomotopyBackground:Homotopies}
Let $A$ be an object of~$\C$.
A cylinder object of~$A$ is a diagram
\begin{equation*}
\xymatrix{ A\oplus A\ar@{>->}[]!R+<4pt,0pt>;[r]_(0.35){(d^0,d^1)} & \widetilde{A}\ar[r]^{\sim}_{s^0} & A }
\end{equation*}
such that $s^0$ is a weak equivalence,
the pair $(d^0,d^1)$ defines a cofibration,
and $s^0 d^0 = s^0 d^1 = \id$.
Any object $A\in\C$
has a cylinder object which can be produced by a factorization~(\ref{M:Factorization}.i)
of the codiagonal $(\id,\id): A\oplus A\rightarrow A$.

A homotopy between morphisms $f,g: A\rightarrow X$
is defined by an extension of the morphism~$(f,g): A\oplus A\rightarrow X$
to a cylinder object:
\begin{equation*}
\xymatrix{ A\oplus A\ar[r]^{(f,g)}\ar[d]_{(d^0,d^1)} & X \\ \widetilde{A}\ar@{.>}[ur]_{h} & }.
\end{equation*}
The homotopy relation is an equivalence relation if $A$ is cofibrant.
Actually, this definition gives the notion of left homotopy.
There is also a symmetrical definition for right homotopy,
but this notion is not used explicitly in the paper.

The set of homotopy classes $[A,X]$ of morphisms $f: A\rightarrow X$
is the quotient of the morphism set $\Mor_{\C}(A,X)$
by the homotopy relation.
If $A$ is cofibrant, then $[A,X]$ defines a functor in $X$
that maps weak equivalences between fibrant objects to bijections.

\subsubsection{Relative cell complexes}\label{HomotopyBackground:RelativeCellComplexes}
Let $\K$ be a class of morphisms in $\C$.
A morphism $f: K\rightarrow L$
forms a relative $\K$-cell complex
if $f$ splits as a composite
\begin{equation*}
K = L_0\rightarrow\dots
\rightarrow L_{\lambda-1}\xrightarrow{j_{\lambda}} L_{\lambda}\rightarrow\cdots
\rightarrow\colim_{\lambda<\mu} L_{\lambda} = L
\end{equation*}
such that each $j_{\lambda}$
is obtained by a pushout of the form
\begin{equation*}
\xymatrix{ \bigoplus_{\alpha} C_\alpha\ar[d]_{(i_\alpha)}\ar[r] &
L_{\lambda-1}\ar@{.>}[d]^{j_{\lambda}} \\
\bigoplus_{\alpha} D_\alpha\ar@{.>}[r] & L_{\lambda} },
\end{equation*}
where $i_\alpha\in\K$.
The composite runs over a (possibly transfinite) ordinal $\mu$.
The set of morphisms $\K$ permits the small object argument
if we have an ordinal $\mu$
such that the following property holds for the domain $C$ of each morphism $i\in \K$:
\begin{enumerate}\renewcommand{\theenumi}{S}
\item\label{SO:SmallObject}\emph{small object axiom}:
any morphism $u: C\rightarrow K$
towards a $\mu$-transfinite composite of relative $\K$-cell complexes
\begin{equation*}
K_0\rightarrow\dots
\rightarrow K_{\lambda-1}\xrightarrow{f_{\lambda}} K_{\lambda}\rightarrow\cdots
\rightarrow\colim_{\lambda<\mu} K_{\lambda} = K
\end{equation*}
factors through a term $K_\lambda$,
for some $\lambda<\mu$.
\end{enumerate}
The small object argument produces a factorization $f = p i$
such that $i$ is a relative $\K$-cell complex
and $p$ has the right lifting property with respect to morphisms of~$\K$.

\subsubsection{Cofibrantly generated model categories}\label{HomotopyBackground:CofibrantlyGenerated}
A model category $\C$
is cofibrantly generated if $\C$
has a set of cofibrations $\I$, called generating cofibrations,
and a set of acyclic cofibrations $\J$, called generating acyclic cofibrations,
such that:
\begin{enumerate}\renewcommand{\theenumi}{CG\arabic{enumi}}
\item\label{CG:Smallness}\emph{smallness axiom}:
The sets $\I$ and $\J$ permit the small object argument.
\item\label{CG:Generation}\emph{generation axioms}:
\begin{enumerate}
\item
The fibrations are characterized by the right lifting property
with respect to the generating acyclic cofibrations $j\in\J$.
\item
The acyclic fibrations are characterized by the right lifting property
with respect to the generating cofibrations $i\in\I$.
\end{enumerate}
\end{enumerate}
Under these conditions,
the small object argument can be applied to $\I$ (respectively, $\J$)
to produce the factorizations $f = p i$ (respectively, $f = q j$)
demanded by axiom M5 of model categories.
Simply note that any relative $\I$-cell (respectively, $\J$-cell) complex
forms a cofibration (respectively, an acyclic cofibration).

From the construction of the small object argument,
one deduces that any (acyclic) cofibration in $\C$
forms a retract of a relative $\I$-cell (respectively, $\J$-cell) complex.
Since fibrations are characterized by axiom~\ref{CG:Generation},
the weak equivalences, generating cofibrations, and generating acyclic cofibrations
are sufficient to characterize the full structure
of a cofibrantly generated model category.

\subsubsection{Examples}\label{HomotopyBackground:Examples}
Our basic examples of model categories
are:
\begin{itemize}
\item
the category of (compactly generated) topological spaces $\Top$ (see~\cite[\S\S 2.4.21-25]{Hovey} and further references therein);
\item
the category of simplicial sets $\Simp$ (see for instance~\cite{GoerssJardine} and \cite[\S 3]{Hovey});
\item
the category of (possibly unbounded) dg-modules $\dg\kk\Mod$, where $\kk$ denotes any fixed ground ring (see for instance~\cite[\S 2.3]{Hovey});
\item
the category of simplicial $\kk$-modules $\s\kk\Mod$ (see for instance~\cite{QuillenHomotopy}).
\end{itemize}
These categories are all cofibrantly generated.
The explicit definitions of these model structures
are given in \emph{loc. cit.}.

\subsubsection{Axioms of symmetric monoidal model categories}\label{HomotopyBackground:SymmetricMonoidalModelCategories}
The notion of symmetric monoidal model category gives the general background of homotopical algebra
for operads and props.
Basically,
a symmetric monoidal model category is a symmetric monoidal category $\C$
equipped with a model structure together with additional axioms
which ensure that cofibrations assemble properly in tensor products:
\begin{enumerate}\renewcommand{\theenumi}{MM\arabic{enumi}}\setcounter{enumi}{-1}
\item\label{MM:Unit}\emph{unit axiom}:
The unit object $\unit$
is cofibrant in $\C$.
\item\label{MM:PushoutProduct}\emph{pushout-product axiom}:
If $i: A\rightarrowtail B$ and $j: C\rightarrowtail D$ are cofibrations,
then so is the natural morphism
\begin{equation*}
(i_*,j_*): A\otimes D\bigoplus_{A\otimes C} B\otimes C\rightarrow B\otimes D
\end{equation*}
arising from the diagram
\begin{equation*}
\xymatrix{ A\otimes C\ar[d]_{A\otimes j}\ar[r]^{i\otimes C} & B\otimes C\ar@{.>}[d]\ar@/^6mm/[ddr]^{B\otimes j} & \\
A\otimes D\ar@/_6mm/[drr]_{i\otimes D}\ar@{.>}[r] & A\otimes D\bigoplus_{A\otimes C} B\otimes C\ar@{.>}[dr]^{(i_*,j_*)} & \\
&& B\otimes D }.
\end{equation*}
This cofibration $(i_*,j_*)$ is also acyclic
whenever $i$ or $j$ is.
\end{enumerate}
These axioms are borrowed from~\cite[\S 4]{Hovey}.

By adjunction,
the pushout-product axiom is equivalent to:
\begin{enumerate}\renewcommand{\theenumi}{MM\arabic{enumi}${}^{\sharp}$}
\item\label{MM:PullbackHom}
If $i: A\rightarrowtail B$ is a cofibration and $p: X\twoheadrightarrow Y$ is a fibration,
then the natural morphism
\begin{equation*}
(i^*,p_*): \Hom_{\C}(B,X)\rightarrow\Hom_{\C}(A,X)\times_{\Hom_{\C}(A,Y)}\Hom_{\C}(B,Y)
\end{equation*}
arising from the diagram
\begin{equation*}
\xymatrix{ \Hom_{\C}(B,X)\ar@{.>}[dr]^{(i^*,p_*)}\ar@/^6mm/[drr]^{p_*}\ar@/_6mm/[ddr]_{i^*} && \\
& \Hom_{\C}(A,X)\times_{\Hom_{\C}(A,Y)}\Hom_{\C}(B,Y)\ar@{.>}[r]\ar@{.>}[d] & \Hom_{\C}(B,Y)\ar[d]^{i^*} \\
& \Hom_{\C}(A,X)\ar[r]_{p_*} & \Hom_{\C}(A,Y) }
\end{equation*}
forms a fibration,
which is also acyclic whenever $i$ or $p$ is.
\end{enumerate}

The categories of topological spaces $\Top$,
simplicial sets $\Simp$,
dg-modules $\dg\kk\Mod$
and simplicial modules $\s\kk\Mod$
all form symmetric monoidal model categories in the sense of this paragraph (see~\cite[\S 4.2]{Hovey}).

\section{Semi-model categories}\label{SemiModelCategories}
The purpose of this section is to review the axioms of semi-model categories
and the definition of semi-model structures by adjunction
from a given model category.

In general,
we are given an adjunction $F: \X\rightleftarrows\A :U$
where $\X$ is a cofibrantly generated model category.
The category $\A$ is supposed to have all colimits and all limits.
For our purpose,
we assume moreover that the functor $U: \A\rightarrow\X$ preserves colimits over non-empty ordinals.
Note that $F(0)$ defines the initial object of $\A$
since $F$ preserves colimits.

In our constructions,
we use morphisms $f: A\rightarrow B$
such that $U(f)$ defines a cofibration in $\X$.
Usually,
we say that a morphism of $\A$
is an $\X$-cofibration
if its image under the functor $U: \A\rightarrow\X$
defines a cofibration in $\X$,
an object $A\in\A$
is $\X$-cofibrant
if the initial morphism $\eta: F(0)\rightarrow A$
forms an $\X$-cofibration.
Note that the initial object $F(0)\in\A$ is $\X$-cofibrant (as the identity morphism $\id: F(0)\rightarrow F(0)$ is a cofibration),
but we do not assume necessarily that $U F(0)$
forms a cofibrant object in $\X$.

\subsubsection{The axioms of semi-model categories}\label{SemiModelCategories:Axioms}
Define classes of weak equivalences, cofibrations and fibrations
in $\A$ by assuming that:
\begin{itemize}
\item
a morphism forms a weak equivalence (respectively, a fibration) in $\A$
if and only if its image under the functor $U: \A\rightarrow\X$
forms a weak equivalence (respectively, a fibration) in $\X$
(we say that the functor $U: \A\rightarrow\X$ creates weak equivalences and fibrations);
\item
cofibrations have the left lifting property with respect to acyclic fibrations.
\end{itemize}
In good cases,
these definitions provide $\A$ with a full model structure.
But in the case of props
we do not have all lifting and factorization properties.
For that reason,
we have to introduce weaker forms of the lifting and factorization axioms
of model categories.

The category $\A$ is said to form a semi-model category
if the following properties hold:
\begin{enumerate}\renewcommand{\theenumi}{M\arabic{enumi}'}\setcounter{enumi}{3}
\item\label{Mprime:Lifting}\emph{lifting axioms}:
\begin{enumerate}
\item
The fibrations have the right lifting property
with respect to the acyclic cofibrations $i: A\rightarrow B$
such that $A$ is cofibrant.
\item
The acyclic fibrations have the right lifting property
with respect to the cofibrations $i: A\rightarrow B$
such that $A$ is cofibrant.
\end{enumerate}
\item\label{Mprime:Factorization}\emph{factorization axioms}:
\begin{enumerate}
\item
Any morphism $f: A\rightarrow B$
has a factorization $f = p i$,
where $i$ is a cofibration and $p$ is an acyclic fibration,
provided that $A$ is cofibrant.
\item
Any morphism $f: A\rightarrow B$
has a factorization $f = q j$,
where $j$ is an acyclic cofibration and $q$ is a fibration,
provided that $A$ is cofibrant.
\end{enumerate}
\end{enumerate}
The two-out-of-three and the retract axioms~\ref{M:TwoOutofThree}-\ref{M:Retract}
of model categories hold automatically and, by assumption, so does the completeness axiom~\ref{M:Completeness}.

The definition of weak equivalences, cofibrations, and fibrations in $\A$
implies readily:
\begin{enumerate}\renewcommand{\theenumi}{M\arabic{enumi}'}\setcounter{enumi}{-1}
\item\label{Mprime:Initial}\emph{initial object axiom}: The initial object of $\A$ is cofibrant.
\setcounter{enumi}{5}
\item\label{Mprime:Fibration}\emph{fibration axioms}:
\begin{enumerate}
\item
The class of (acyclic) fibrations in $\A$ is stable under (possibly transfinite) composites.
\item
The class of (acyclic) fibrations in $\A$ is stable under pullbacks.
\end{enumerate}
\end{enumerate}
Note that these properties are not implied by axioms~\ref{Mprime:Lifting}-\ref{Mprime:Factorization},
because these axioms are not sufficient to characterize the class of (acyclic) cofibrations of~$\A$
from the class of (acyclic) fibrations.
For that reason, properties~\ref{Mprime:Initial} and~\ref{Mprime:Fibration}
should be added to the axioms of semi-model categories.

Properties~\ref{Mprime:Initial},
\ref{M:TwoOutofThree}-\ref{M:Retract},
\ref{Mprime:Lifting}-\ref{Mprime:Factorization}
and~\ref{Mprime:Fibration}
are sufficient to do homotopical algebra.
In particular,
the definition and the usual properties of homotopy classes~$[A,B]$
have literal generalizations in the setting of semi-model categories (see~\cite{Spitzweck}).

\subsubsection{Generating (acyclic) cofibrations in semi-model categories}\label{SemiModelCategories:CofibrantlyGenerated}
In usual cases,
the model category $\X$ comes equipped with a set of generating cofibrations $\I$
and a set of generating acyclic cofibrations $\J$.

Form the sets of morphisms $F\I = \{F(i), i\in\I\}$ and $F\J = \{F(j), j\in\J\}$.
The adjunction relation implies readily that a morphism $p$
is a fibration (respectively, an acyclic fibration) in $\A$
if and only if $p$ has the right lifting property with respect to morphisms $F(j)\in F\J$
(respectively, $F(i)\in F\I$).
Hence,
a natural idea is to use $F\I$ (respectively, $F\J$)
as a set of generating cofibrations (respectively, acyclic cofibrations)
to prove axioms~\ref{Mprime:Lifting}-\ref{Mprime:Factorization}.

In the context of a semi-model category,
we only need to apply the small object argument to morphisms $f: A\rightarrow B$
such that $A$ is cofibrant.
Therefore,
we say that a set of morphisms $F\K = \{F(i),\ i\in\K\}$
permits the small object argument
if we have an ordinal $\mu$
such that the following property holds for the domain $F(C)$ of each morphism $F(i)\in F\K$:
\begin{enumerate}\renewcommand{\theenumi}{S'}
\item\label{SOprime:SmallObject}
any morphism $u: F(C)\rightarrow A$
towards a $\mu$-transfinite composite of relative $F\K$-cell complexes
\begin{equation*}
A_0\rightarrow\dots
\rightarrow A_{\lambda-1}\xrightarrow{f_{\lambda}} A_{\lambda}\rightarrow\cdots
\rightarrow\colim_{\lambda<\mu} A_{\lambda} = A
\end{equation*}
such that $A_0$ is cofibrant in $\A$
factors through a term $A_\lambda$,
for some $\lambda<\mu$.
\end{enumerate}
Under this assumption,
the small object argument can be applied to produce,
for any morphism $f: A\rightarrow B$ with a cofibrant domain $A$,
a factorization $f = p i$ such that $i$ is a relative $F\K$-cell complex
and $p$ has the right lifting property with respect to the morphisms of $F\K$.

The category $\A$ is cofibrantly generated as a semi-model category with $F\I$ (respectively, $F\J$)
as generating cofibrations (respectively, acyclic cofibrations)
if axiom~\ref{SOprime:SmallObject} holds for $F\I$ and $F\J$.

The next theorem gives sufficient conditions
to ensure this property
and to prove the axioms of semi-model categories.
For technical reason,
we consider the class $\X_c$ of cofibrations of $\X$
and the associated class of morphisms $F\X_c = \{F(i),\ i\in\X_c\}$ in $\A$.

\begin{thm}\label{SemiModelCategories:AdjointStructure}
Suppose that:
\begin{enumerate}\renewcommand{\theenumi}{*}
\item\label{AdjointStructure:Requirement}
for any pushout
\begin{equation*}
\xymatrix{ F(X)\ar[r]\ar[d]_{F(i)} & A\ar@{.>}[d]^{f} \\
F(Y)\ar@{.>}[r] & B }
\end{equation*}
such that $A$ is an $F\X_c$-cell complex in $\A$,
the morphism $U(f)$
forms a cofibration (respectively an acyclic cofibration) in $\X$
whenever $i$ is a cofibration (respectively an acyclic cofibration) in $\X$.
\end{enumerate}
Then we have:
\begin{enumerate}
\item
the sets $F\I$ and $F\J$ permit the small object argument,
\item
the lifting and factorization axioms~\ref{Mprime:Lifting}-\ref{Mprime:Factorization} of~\S\ref{HomotopyBackground:ModelCategoryAxioms}
hold in $A$,
\end{enumerate}
so that $\A$ forms a cofibrantly generated semi-model category over $\X$,
and moreover:
\begin{enumerate}\setcounter{enumi}{2}
\item
the functor $U: \A\rightarrow\X$ maps cofibrations with a cofibrant domain to cofibrations.
\end{enumerate}
\end{thm}

The proof is divided into a series of verifications.
Theorem~\ref{SemiModelCategories:AdjointStructure} is a variant of~\cite[Theorem 12.1.4]{FresseModuleBook}.
Therefore we will give references to~\emph{loc. cit.} for the details of some verifications.

\begin{lemm}\label{SemiModelCategories:CofibrantRetracts}
The cofibrant objects of $\A$ form retracts of relative $F\I$-cell complexes.
\end{lemm}

\begin{proof}
Let
\begin{equation*}
A_0\rightarrow\dots
\rightarrow A_{\lambda-1}\xrightarrow{f_{\lambda}} A_{\lambda}\rightarrow\cdots
\rightarrow\colim_{\lambda<\mu} A_{\lambda} = A
\end{equation*}
be a $\mu$-transfinite composite of relative $F\I$-cell complexes
such that $A_0 = F(0)$.
For an object $X\in\X$,
we have the adjunction relation
\begin{equation*}
\Mor_{\A}(F(X),\colim_{\lambda} A_{\lambda}) = \Mor_{\A}(X,U(\colim_{\lambda} A_{\lambda})).
\end{equation*}
Recall that $U$ is supposed to preserve colimits over non-empty ordinals.
Therefore we have moreover $U(\colim_{\lambda} A_{\lambda}) = \colim_{\lambda} U(A_{\lambda})$.
The assumption implies by induction
that each $U(f_\lambda): U(A_{\lambda-1})\rightarrow U(A_{\lambda})$
forms a cofibration in $\X$.
Thus,
the morphism $U(f_\lambda)$ is not necessarily a relative $\I$-cell complex,
but forms at least a retract of a relative $\I$-cell complex.
This assertion implies that the morphism $X\rightarrow\colim_{\lambda} U(A_{\lambda})$
factors through a term $U(A_{\lambda})$ for some $\lambda<\mu$
if $X$ satisfies the smallness axiom~\ref{SOprime:SmallObject}
with respect to relative $\I$-cell complexes in $\X$ (see~\cite[Proof of proposition 11.1.14]{FresseModuleBook}).
By adjunction,
we obtain that $F(X)\rightarrow\colim_{\lambda} A_{\lambda}$
factors through $A_{\lambda}$.

According to these observations,
we can apply the small object argument to produce a factorization
$F(0)\xrightarrow{i} B\xrightarrow{p} A$
of the initial morphism $\eta: F(0)\rightarrow A$
of any object $A\in\A$,
so that:
\begin{itemize}
\item the morphism $i$ is a relative $F\I$-cell complex,
\item the morphism $p$ has the right lifting property with respect to morphisms of the form $F(i)$, $i\in\I$.
\end{itemize}
By adjunction,
we deduce that $U(p)$ has the right lifting property with respect to all morphisms $i\in\I$,
and hence forms a fibration in $\A$.
If $A$ is cofibrant in $\A$,
then $\eta$ has the left lifting property with respect to $p$,
from which we deduce that $A$ forms a retract of $B$,
and hence of a relative $F\I$-cell complex.
\end{proof}

\begin{lemm}\label{SemiModelCategories:CofibrantPushout}
The assumption of Theorem~\ref{SemiModelCategories:AdjointStructure}
holds whenever $A$ is a cofibrant object in $\A$:
for any pushout of the form
\begin{equation*}
\xymatrix{ F(X)\ar[r]\ar[d]_{F(i)} & A\ar@{.>}[d]^{f} \\
F(Y)\ar@{.>}[r] & B },
\end{equation*}
the morphism $U(f)$ forms a cofibration (respectively an acyclic cofibration) in $\X$
if $i$ does.
\end{lemm}

\begin{proof}
The assumption of Theorem~\ref{SemiModelCategories:AdjointStructure}
implies by induction that any $F\I$-cell complex
is $\X$-cofibrant.
By Lemma~\ref{SemiModelCategories:CofibrantRetracts},
we have an $F\I$-cell complex $C$ such that $A$ is a retract of $C$.
Explicitly,
we have morphisms $A\xrightarrow{j} C\xrightarrow{s} A$
such that $s j = \id$.
Form the diagram of pushouts
\begin{equation*}
\xymatrix{ F(X)\ar[r]\ar[d]_{F(i)} & A\ar@{.>}[d]^{f}\ar[r]^{j} & C\ar@{.>}[d]^{g}\ar[r]^{s} & A\ar@{.>}[d]^{f} \\
F(Y)\ar@{.>}[r] & B\ar@{.>}[r] & D\ar@{.>}[r] & B }.
\end{equation*}
The assumption of Theorem~\ref{SemiModelCategories:AdjointStructure}
asserts that $U(g)$ forms a cofibration in $\X$.
The commutativity of the diagram implies
that $U(f)$ is a retract of $U(g)$,
and hence forms a cofibration as well.
\end{proof}

\begin{lemm}\label{SemiModelCategories:SmallObjectArgument}
The sets $F\I$ and $F\J$ permit the small object argument.
\end{lemm}

\begin{proof}
The verification of this claim involves immediate generalizations of arguments
of Lemma~\ref{SemiModelCategories:CofibrantRetracts}.
Let
\begin{equation*}
A_0\rightarrow\dots
\rightarrow A_{\lambda-1}\xrightarrow{f_{\lambda}} A_{\lambda}\rightarrow\cdots
\rightarrow\colim_{\lambda<\mu} A_{\lambda} = A
\end{equation*}
be a $\mu$-transfinite composite of relative $F\I$-cell complexes
such that $A$ is cofibrant.

For an object $X\in\X$,
we have
\begin{equation*}
\Mor_{\A}(F(X),\colim_{\lambda} A_{\lambda})
= \Mor_{\A}(X,U(\colim_{\lambda} A_{\lambda}))
= \Mor_{\A}(X,\colim_{\lambda} U(A_{\lambda})).
\end{equation*}
The assertion of Lemma~\ref{SemiModelCategories:CofibrantRetracts}
implies by induction
that each $U(f_\lambda): U(A_{\lambda-1})\rightarrow U(A_{\lambda})$
forms a cofibration in $\X$.
Thus,
the morphism $U(f_\lambda)$ is not necessarily a relative $\I$-cell complex,
but forms at least a retract of a relative $\I$-cell complex.
Again
this assertion implies that the morphism $X\rightarrow\colim_{\lambda} U(A_{\lambda})$
factors through a term $U(A_{\lambda})$ for some $\lambda<\mu$
if $X$ satisfies the smallness axiom~(\ref{SO:SmallObject})
with respect to relative $\I$-cell complexes in $\X$.
The conclusion follows readily.
\end{proof}

\begin{prop}\label{SemiModelCategories:LiftingFactorization}
The lifting and factorization axioms~\ref{Mprime:Lifting}-\ref{Mprime:Factorization} hold in $\A$.
\end{prop}

\begin{proof}
Let $f$ be a morphism with a cofibrant domain.
The small object argument produces a factorization $f = p i$
where $i$ is a relative $F\K$-cell complex
and $p$ has the right lifting property with respect to the morphisms of $F\K$.
For $\K = \I$ (respectively, $\K = \J$),
this gives the factorization required by axiom~\ref{Mprime:Factorization}.i (respectively, \ref{Mprime:Factorization}.ii).

Axiom~\ref{Mprime:Lifting}.i
is tautologically satisfied in $\A$.
Let $i$ be any acyclic cofibration.
In order to prove axiom~\ref{Mprime:Lifting}.ii,
we use the small object argument to produce a factorization $i = q j$,
where $j$ is a relative $F\J$-cell complex
and $q$ has the right lifting property with respect to the morphisms of $F\J$.
This morphism $q$ forms a fibration
and a weak equivalence by the two-out-of-three axiom.
Hence $p$ forms an acyclic fibration.
From the lifting axiom~\ref{Mprime:Lifting}.i,
we deduce readily that $i$ is a retract of $j$.
Note that relative $F\J$-cell complexes
have the left lifting property with respect to fibrations.
Use that retracts inherit lifting properties to conclude
that $i$ has the left lifting property with respect to fibrations,
as required.
\end{proof}

\begin{prop}\label{SemiModelCategories:Cofibrations}
The cofibrations of $\A$ with a cofibrant domain are $\X$-cofibrations.
\end{prop}

\begin{proof}
The requirement of Theorem~\ref{SemiModelCategories:AdjointStructure}
implies by induction that any relative $F\I$-cell complex $f: A\rightarrow B$
such that $A$ is $\X$-cofibrant
forms an $\X$-cofibration.
So does any cofibration with a cofibrant domain
since we have proved that any of these cofibrations
arises as a retract of a relative $F\I$-cell complex.
\end{proof}

This verification achieves the proof of Theorem~\ref{SemiModelCategories:AdjointStructure}.\qed

\section{The semi-model category of props}\label{PropSemiModel}
The goal of this section is to define a semi-model structure on the category of props.
For this aim,
we apply Theorem~\ref{SemiModelCategories:AdjointStructure} to the composite adjunction
\begin{equation*}
\xymatrix{ \C^{\ACat}\ar@<+2pt>[r]^{\phi_!} & \C^{\BCat}\ar@<+2pt>[r]^{F}\ar@<+2pt>[l]^{\phi^*} & \Prop\ar@<+2pt>[l]^{U} },
\end{equation*}
between the category of props and $\C^{\ACat}$.
The technical verification of the sufficient condition of Theorem~\ref{SemiModelCategories:AdjointStructure}
is postponed to the appendix part.
Note that we really apply Theorem~\ref{SemiModelCategories:AdjointStructure}
to the composite adjunction between $\C^{\ACat}$ and $\Prop$,
but the definition of the free prop is more natural (and more usual) on the category of $\Sigma_*$-biobjects $\C^{\BCat}$.
For this reason,
we also use the intermediate adjunction between $\C^{\BCat}$ and $\Prop$.

In fact,
we cannot even have a semi-model structure on the whole category of props,
because of symmetries inherent in the horizontal composition product.
To fix this problem,
we restrict ourselves to the full subcategory of props
with non-empty inputs (outputs).
This category is defined in the course
of our construction.

First of all,
we review the overall definition of a cofibrantly generated model structure
on a category of diagrams:

\begin{prop}[{see~\cite[\S 11.6]{Hirschhorn}}]\label{PropSemiModel:DiagramModelCategory}
Let $\ICat$ be any small category.
The category of diagrams $\C^{\ICat}$ is equipped with a full model structure
so that:
\begin{itemize}
\item
a morphism $f: K\rightarrow L$ forms a weak equivalence (respectively, a fibration) in $\C^{\ICat}$
if and only if its components $f: K(\alpha)\rightarrow L(\alpha)$, $\alpha\in\ICat$, are all weak equivalences (respectively, fibrations)
in $\C$;
\item
the cofibrations are characterized by the left lifting property with respect to acyclic fibrations.
\end{itemize}
This model structure has generating (acyclic) cofibrations
defined by morphisms of the form
\begin{equation*}
i\otimes G_{\alpha}: C\otimes G_{\alpha}\rightarrow D\otimes G_{\alpha},
\end{equation*}
where $i: C\rightarrow D$ ranges over the generating (acyclic) cofibrations of~$\C$
and $G_{\alpha}$ is the enriched Yoneda diagram associated to any object $\alpha\in\ICat$.
\end{prop}

\begin{proof}
This statement is a variation of the result of~\cite[\S 11.6]{Hirschhorn}
in the case where $\C$ forms a closed symmetric monoidal model category.
\end{proof}

Note that:

\begin{obsv}\label{PropSemiModel:DiagramExtensionRestriction}
Let $\phi: \ICat\rightarrow\JCat$ be any functor between small categories.
The extension and restriction functors $\phi_!: \C^{\ICat}\rightleftarrows\C^{\JCat} :\phi^*$
define a Quillen pair of adjoint functors.
In the case of the functor $\phi: \ACat\rightarrow\BCat$ of~\S\ref{PropBackground:DiagramCategories},
we have moreover:
\begin{itemize}
\item
the functor $\phi^*: \C^{\BCat}\rightarrow\C^{\ACat}$
creates fibrations, creates weak equivalences,
and maps cofibrations to cofibrations;
\item
the functor $\phi_!: \C^{\ACat}\rightarrow\C^{\BCat}$
defines a one-to-one correspondence between generating (acyclic) cofibrations.
\end{itemize}
\end{obsv}

According to the definition of Proposition~\ref{PropSemiModel:DiagramModelCategory},
the generating cofibrations of the category of $\ACat$-diagrams (respectively, $\BCat$-diagrams)
are the tensor products $i\otimes G: C\otimes G\rightarrow D\otimes G$
where $i: C\rightarrow D$ ranges over the generating (acyclic) cofibrations of~$\C$
and $G$ ranges over the enriched Yoneda diagrams $G = G_{(m,n)}$ (respectively, $G = \phi_! G_{(m,n)}$).
Assertion (2) follows from the formal identity $\phi_!(i\otimes G) = i\otimes\phi_! G$.

\medskip
The next stage of our verification consists in proving the existence of limits and colimits
in the category of props.
The next assertion is standard for objects equipped with a multiplicative structure
and is proved by the usual line of argument:

\begin{prop}\label{PropSemiModel:Limits}
The forgetful functor $U: \Prop\rightarrow\C^{\BCat}$ creates limits in the category of props $\Prop$.\qed
\end{prop}

Recall that tensor powers preserve certain special colimits,
namely reflexive coequalizers and filtered colimits.
Since the free prop is made of colimits of tensor products (see~\S\S\ref{Graphs:FreePropsConstruction}-\ref{Graphs:FreePropDefinition}),
we obtain:

\begin{obsv}\label{PropSemiModel:FreeColimits}
The free prop functor $\Free: \C^{\BCat}\rightarrow\Prop$ preserves reflexive coequalizers and filtered colimits.\qed
\end{obsv}

From this observation,
we obtain readily:

\begin{prop}\label{PropSemiModel:CreatedColimits}
The forgetful functor $U: \Prop\rightarrow\C^{\BCat}$ creates reflexive coequalizers and filtered colimits
in the category of props $\Prop$.\qed
\end{prop}

The next proposition is a standard consequence of the previous result:

\begin{prop}\label{PropSemiModel:Colimits}
The colimit of any diagram $i\mapsto\POp_i$
in the category of props $\Prop$
can be realized by a reflexive coequalizer of the form
\begin{gather*}\\
\xymatrix{ \Free(\colim_i U \Free U(\POp_i))\ar@<+2pt>[r]\ar@<-2pt>[r] & \Free(\colim_i U(\POp_i))\ar@{.>}[r]\ar@/_8mm/[l] &
\colim_i\POp_i },
\end{gather*}
where we perform the colimit of the diagrams $i\mapsto U(\POp_i)$ and $i\mapsto U \Free U(\POp_i)$ in~$\C^{\BCat}$.\qed
\end{prop}

Note that the initial object of the category of props
has an easy description:

\begin{fact}\label{PropSemiModel:InitialObject}
The $\Sigma_*$-biobject $I$
such that
\begin{equation*}
I(m,n) = \begin{cases} \unit[\Sigma_n], & \text{if $m = n$}, \\
0, & \text{otherwise},
\end{cases}
\end{equation*}
is equipped with an obvious prop structure
and forms the initial object in the category of props.
\end{fact}

Observe that $I$ does not form a cofibrant $\Sigma_*$-biobject in general.

\subsubsection{Props with non-empty inputs (outputs)}\label{PropSemiModel:NonEmptyInputs}
The horizontal composition product $\circ_h: \POp(k,m)\otimes\POp(l,n)\rightarrow\POp(k+l,m+n)$
defines a symmetric operation
on components such that $k = l = m = n = 0$.
This is not a problem if the base category is the category of topological spaces, simplicial sets,
simplicial modules, or dg-modules over a ring of characteristic $0$,
but in general this symmetry relation gives an obstruction to the definition
of a model structure (see Remark~\ref{PropSemiModel:DGModuleCase}).
To fix this problem,
we restrict ourselves to props $\POp$
such that:
\begin{equation*}
\POp(0,n) = \begin{cases} \unit, & \text{if $n = 0$}, \\ 0, & \text{otherwise}. \end{cases}
\end{equation*}
If this condition is satisfied, then we say that $\POp$ is a prop with non-empty inputs.

Of course,
we have a symmetrical notion of prop with non-empty outputs.

Many usual categories of algebras, coalgebras, and bialgebras, are associated to props
with either non-empty inputs or non-empty outputs.
But the prop associated to the category of unital augmented Hopf algebras has both components with empty inputs
and components with empty outputs
since the structure of a unital augmented Hopf algebra includes a unit operation $\eta: \kk\rightarrow A$
and an augmentation operation $\epsilon: A\rightarrow\kk$.
This observation implies that the homotopical algebra of unital augmented Hopf algebras
cannot be modelled by props in positive characteristic.

Let $\Prop_0$ be the full subcategory of~$\Prop$
formed by props with non-empty inputs.
The embedding $\Prop_0\hookrightarrow\Prop$ has a left adjoint which maps a prop $\POp\in\Prop$
to the reduced prop $\overline{\POp}\in\Prop_0$
such that:
\begin{equation*}
\overline{\POp}(m,n) = \begin{cases} \unit, & \text{if $m = n = 0$}, \\
0, & \text{if $m = 0$ and $n\not=0$}, \\
\POp(m,n), & \text{otherwise}. \end{cases}
\end{equation*}

The adjunction relations
\begin{equation*}
\xymatrix{ \C^{\ACat}\ar@<+2pt>[r]^{\phi_!} & \C^{\BCat}\ar@<+2pt>[r]^{F}\ar@<+2pt>[l]^{\phi^*} & \Prop\ar@<+2pt>[l]^{U} },
\end{equation*}
restrict to adjunction relations
\begin{equation*}
\xymatrix{ \C^{\ACat}_0\ar@<+2pt>[r]^{\phi_!} & \C^{\BCat}_0\ar@<+2pt>[r]^{F}\ar@<+2pt>[l]^{\phi^*} & \Prop_0\ar@<+2pt>[l]^{U} },
\end{equation*}
where we take subcategories of diagrams $K$
such that $K(m,n) = 0$ if $m = 0$.
Simply remove the component $\POp(0,0) = \unit$ from props $\POp\in\Prop_0$
to obtain a forgetful functor $U: \Prop_0\rightarrow\C^{\BCat}_0$.

The assertions of~\S\S\ref{PropSemiModel:Limits}-\ref{PropSemiModel:InitialObject}
remain obviously valid if we replace the category of props by props with non-empty inputs,
and the diagram categories~$\C^{\ACat}$ and~$\C^{\BCat}$
by their subcategories~$\C^{\ACat}_0$ and~$\C^{\BCat}_0$.

\medskip
We prove the following theorem:

\begin{thm}\label{PropSemiModel:Result}
The category of props with non-empty inputs $\Prop_0$
in any cofibrantly generated symmetric monoidal model category $\C$
inherits a semi-model structure
so that:
\begin{itemize}
\item
the forgetful functor $\phi^* U: \Prop_0\rightarrow\C^{\ACat}_0$
creates weak equivalences, creates fibrations,
and maps the cofibrations with a cofibrant domain to cofibrations;
\item
the morphisms of free props $\Free(i): \Free(K)\rightarrow\Free(L)$,
where $i: K\rightarrow L$
ranges over (acyclic) cofibrations of $\Sigma_*$-biobjects with non-empty inputs,
form generating (acyclic) cofibrations of the category of props.\qed
\end{itemize}

The same result holds in the symmetrical case of props with non-empty outputs.
\end{thm}

We check that the adjunction $\Free: \C^{\ACat}_0\rightleftarrows\Prop_0 :U$
satisfies the requirement of Theorem~\ref{SemiModelCategories:AdjointStructure}
and we deduce our result from that statement.
We postpone this technical verification to the appendix.

\subsubsection{Remark: the case of dg-modules over a ring}\label{PropSemiModel:DGModuleCase}
\begin{figure}
\begin{equation*}
\xymatrix@!C=0.5em@!R=0.5em@M=2pt{ \ar@{.}[r] & 1\ar[d]\ar@{.}[r] & \\ & M(1,p)\ar[dl]\ar[dr] & \\
L(1,1)\ar[dr] & \cdots & L(1,1)\ar[dl] \\
& M(p,1)\ar[d] & \\ \ar@{.}[r] & 1\ar@{.}[r] & }
\end{equation*}
\caption{A symmetric summand of the free prop}\label{Fig:NonAcyclicSummand}
\end{figure}
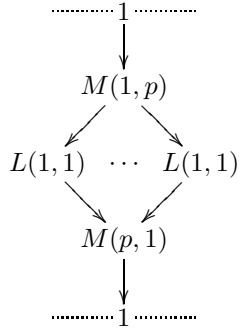

If the ground ring $\kk$ is a field of positive characteristic,
then we cannot obtain a full model structure
on the category of props with non-empty inputs (or outputs),
because pushouts along morphisms $\Free(i): \Free(K)\rightarrow\Free(L)$
such that $i$ is an acyclic cofibration do not necessarily
form weak equivalences (in contradiction with a standard property of model categories).
To see easily the existence of obstructions,
we can assume that $K = 0$ and $L$ is an acyclic dg-module.
Take $\POp = \Free(M)$, so that $\POp\vee\Free(L) = \Free(M\oplus L)$.
Observe that the expansion of $\Free(M\oplus L)$
includes a summand of the form
\begin{equation*}
M(p,1)\otimes_{\Sigma_p}(L(1,1)^{\otimes p})\otimes_{\Sigma_p} M(1,p),
\end{equation*}
where $p = \Char\kk$.
In the graphical representation of~\S\ref{Graphs},
this summand consists of composites represented by the graph of Figure~\ref{Fig:NonAcyclicSummand}).
This summand is not acyclic if $M$ does not form a cofibrant object in $\Sigma_*$-bimodules.

In addition, we cannot have a semi-model structure on the whole category of props $\Prop$,
because the free prop $\Free(M)$ includes symmetric tensor products as summands.
Indeed,
for the $\Sigma_*$-biobject $M = C\otimes G_{0 0}$,
which has
\begin{equation*}
C\otimes G_{0 0}(m,n) = \begin{cases} C, & \text{if $m=n=0$}, \\ 0, & \text{otherwise}, \end{cases}
\end{equation*}
we obtain
\begin{equation*}
\Free(C\otimes G_{0 0})(m,n) = \begin{cases} \Sym(C), & \text{if $m=n=0$}, \\ 0, & \text{otherwise}, \end{cases}
\end{equation*}
where $\Sym(C) = \bigoplus_{n=0}^{\infty} (C^{\otimes n})_{\Sigma_n}$
is the symmetric algebra generated by $C$.
In the graphical representation of~\S\ref{Graphs},
the free prop $\Free(C\otimes G_{0 0})$
is just a disjoint union of vertices labeled by $C$.
The morphism $\Sym(i): \Sym(C)\rightarrow\Sym(D)$
induced by a generating acyclic cofibration of dg-modules $i: C\trivialcofib D$
is not a weak equivalence in positive characteristic
and hence, neither is the morphism of free props
\begin{equation*}
\Free(i\otimes G_{0 0}): \Free(C\otimes G_{0 0})\rightarrow\Free(D\otimes G_{0 0}),
\end{equation*}
though this morphism has the left lifting property with respect to fibrations.

\subsubsection{Remark: difficulties specific to props}\label{PropSemiModel:Difficulties}
The result of Theorem~\ref{PropSemiModel:Result} can be improved for the category of operads (respectively, properads)
in the sense that the lifting and factorization axioms
hold for cofibrations $i: \POp\rightarrow\QOp$
such that $\POp$ is cofibrant in the underlying category of $\Sigma_*$-objects (see~\cite{Spitzweck}).
This is not the case for props, because the horizontal composition product
makes the description of coproducts more complicated (see~\S\ref{PropPushoutHomotopy:CanonicalCoproductForm}).

\section{Path objects and applications for particular base categories}\label{PropPathObjectArgument}
The example of dg-modules over a field of positive characteristic
shows that we have not always a full model structure on the category of props.
Moreover, we have to restrict ourselves to props with non-empty inputs (or non-empty outputs)
to have proper homotopical properties.

The purpose of this section is to study the case of props in a category of dg-modules over a ring $\kk$ such that $\QQ\subset\kk$,
in a category of simplicial modules,
in the category of simplicial sets,
and in the category of topological spaces.
For these instances of symmetric monoidal model categories,
we prove that the whole category of props inherits a full model structure.
Besides,
the lifting and factorization axioms have a simpler verification
in these examples than in the general case.

First of all,
we review the structures particular to dg-modules over a field of characteristic zero,
simplicial modules, simplicial sets and topological spaces
which permit this simplification.
We give the alternate proof of the axioms of model categories
afterwards.

\medskip
To begin with,
the axiom of~\S\ref{HomotopyBackground:RelativeCellComplexes}
gives the general setting of the small object argument,
but better properties hold in the usual examples of this section:

\begin{fact}[{see for instance~\cite{Hovey}}]\hspace*{2mm}\label{PropPathObjectArgument:SmallObjectArgument}
\begin{enumerate}
\item\label{EverythingSmallObjectArgument}
The domain of generating (acyclic) cofibrations of dg-modules, simplicial modules, and simplicial sets
is small with respect to every countable composite
of morphisms,
and not only with respect to composites of relative cell complexes.
\item\label{TopologicalSmallObjectArgument}
The domain of generating (acyclic) cofibrations of topological spaces
is small with respect to composites
of topological inclusions of spaces (see~\cite[Lemma 2.4.1]{Hovey}).
\end{enumerate}
\end{fact}

By convention,
a fibrant replacement of an object $X$
in a model category $\C$
is a fibrant object $R_X\in\C$
together with a weak equivalence $X\xrightarrow{\sim} R_X$.
Equivalently,
a fibrant replacement of $X$ sits in a factorization $X\xrightarrow{\sim} R_X\twoheadrightarrow *$
of the final morphism $X\rightarrow *$.
Factorization axiom~(\ref{M:Factorization}.ii)
implies that any object $X\in\C$
has a fibrant replacement.
In the usual model categories studied in this section, we have better:

\begin{fact}\hspace*{2mm}\label{PropPathObjectArgument:FibrantReplacement}
\begin{enumerate}
\item\label{TrivialFibrantReplacement}
Each object of the category of dg-modules, simplicial modules, and topological spaces,
is fibrant.
\item\label{MonoidalFibrantReplacement}
In the category of simplicial sets, not all objects are fibrant,
but we have fibrant replacements $R: X\mapsto R_X$
which are functorial in $X\in\C$
and lax symmetric monoidal
in the sense that we have natural weak equivalences
\begin{equation*}
\xymatrix@R=1em@C=1em{ && R_X\otimes R_Y\ar@{.>}[dd] \\
X\otimes Y\ar[drr]!UL_{\sim}\ar[urr]!DL^{\sim} && \\
&& R_{X\otimes Y} }
\end{equation*}
which preserve the unit, associativity and symmetry of cartesian products.
The category of simplicial sets is also right proper: pullbacks of weak equivalences along fibrations
give weak equivalences.
\end{enumerate}
\end{fact}

Kan's functor $R_X = \Ex^{\infty}(X)$
or the composite $R_X = S(|X|)$
of the singular complex of the topological realization of $X\in\Simp$
give examples of functorial fibrant replacements (see~\cite{GoerssJardine})
that satisfy the condition of~(\ref{MonoidalFibrantReplacement}).
The existence of a symmetric monoidal transformation
$R_X\otimes R_Y\rightarrow R_{X\otimes Y}$
implies that the fibrant replacements $R_{\POp} = \{R_{\POp(m,n)}\}$
of the components of a prop~$\POp$
inherit a prop structure and define a fibrant replacement of~$\POp$
in the category of props.

Recall that a category all of whose objects are fibrant is automatically right proper (see~\cite[Corollary 13.1.3]{Hirschhorn}).

\medskip
Recall that a path object of an object $X$
in a model category $\C$
is an object $P(X)\in\C$
together with morphisms
\begin{equation*}
\xymatrix{ X\ar[r]^(0.4){s_0} & P(X)\ar@<+2pt>[r]^(0.6){d_0}\ar@<-2pt>[r]_(0.6){d_1} & X }
\end{equation*}
such that $d_0 s_0 = d_1 s_0 = \id$, the morphism $s_0$ is a weak equivalence,
and $(d_0,d_1): P(X)\rightarrow X\times X$ is a fibration.
If $X$ is fibrant, then $d_0$ and $d_1$ are acyclic fibrations.

The next assertion
holds in categories of dg-modules over a ring of characteristic $0$,
in any category of simplicial modules,
in the category of simplicial sets,
in the category of topological spaces,
but not in categories of dg-modules over a field of positive characteristic.

\begin{fact}\label{PropPathObjectArgument:PathObjects}
In the model categories of dg-modules over a ring of characteristic $0$,
of simplicial modules over any ring,
of simplicial sets,
and topological spaces,
we have an object $P(X)$ naturally associated to any $X\in\C$
together with natural transformations $s_0: X\rightarrow P(X)$ and $d_0,d_1: P(X)\rightarrow X$
such that $P(X)$ defines a path object of~$X$
when $X$ is fibrant in~$\C$.
Moreover, this natural path object $P(X)$
defines a lax monoidal functor in the sense that we have natural transformations
\begin{equation*}
\xymatrix@R=0.8em@C=0.8em{ && P(X)\otimes P(Y)\ar@<-2pt>[]!R;[drr]\ar@<+2pt>[]!R;[drr]\ar@{.>}[dd] && \\
X\otimes Y\ar[drr]!L\ar[urr]!L &&&& X\otimes Y \\
&& P(X\otimes Y)\ar@<-2pt>[]!R;[urr]\ar@<+2pt>[]!R;[urr] && }
\end{equation*}
which preserve the unit, associativity and symmetry of tensor products.
\end{fact}

The simplicial (topological) mapping spaces $P(X) = X^{\Delta^1}$
work in the context of simplicial sets, simplicial $\kk$-modules,
and topological spaces.
In the case of dg-modules over a ring $\kk$
such that $\QQ\subset\kk$,
we can take the tensor product $P(C) = C\otimes\Omega^*(\Delta^1)$,
where $\Omega^*(\Delta^1)$ is Sullivan's dg-algebra of differentials
of the simplicial interval $\Delta^1$.
The natural transformation $P(C\otimes D)\rightarrow P(C)\otimes P(D)$
is yielded by the product of~$\Omega^*(\Delta^1)$.
This natural transformation is unital, associative and commutative,
because so is the dg-algebra~$\Omega^*(\Delta^1)$.

The existence of a symmetric monoidal transformation
$P(X)\otimes P(Y)\rightarrow P(X\otimes Y)$
implies that the path objects $P(\POp(m,n))$
of the components of a prop $\POp$
inherit a prop structure and define a path object of $\POp$
in the category of props (at least if $\POp$ is fibrant).

\medskip
The structures introduced in facts~\ref{PropPathObjectArgument:SmallObjectArgument}-\ref{PropPathObjectArgument:PathObjects}
are already sufficient to obtain:

\begin{lemm}\label{PropPathObjectArgument:FibrationFactorization}
In the situation of facts~\ref{PropPathObjectArgument:SmallObjectArgument}-\ref{PropPathObjectArgument:PathObjects},
any morphism $g: \POp\rightarrow\QOp$ in the category of props~$\Prop$
has a factorization $g = r u$
such that $r$ is a fibration and $u$ is a weak equivalence.
\end{lemm}

\begin{proof}
Form the pullback diagram
\begin{equation*}
\xymatrix{ R_{\POp}\ar[r]^{R_g}\ar@/_8mm/[ddr]_{=}\ar@{.>}[dr]^(0.6){(\id,s_0 R_g)} & R_{\QOp}\ar@/^/[dr]^{s_0} & \\
& R_{\POp}\times_{R_{\QOp}} P(R_{\QOp})\ar@{.>}[r]\ar@{.>>}[d]_{\sim} & P(R_{\QOp})\ar@{->>}[d]_{\sim}^{d_1} \\
& R_{\POp}\ar[r]_{R_g} & R_{\QOp} }.
\end{equation*}
By the universal definition of pullbacks,
we have a natural morphism
$(\id,s_0 R_g): R_{\POp}\rightarrow R_{\POp}\times_{R_{\QOp}} P(R_{\QOp})$.
This morphism forms a weak equivalence by the two-out-of-three axiom.
The object $R_{\POp}\times_{R_{\QOp}} P(R_{\QOp})$
also fits in a pullback
of the form:
\begin{equation*}
\xymatrix{ R_{\POp}\times_{R_{\QOp}} P(R_{\QOp})\ar@{.>>}[d]_{(\id,d_0)}\ar@{.>}[r] &
P(R_{\QOp})\ar@{->>}[d]^{(d_0,d_1)} \\
R_{\POp}\times R_{\QOp}\ar[r]_{(R_g,\id)} & R_{\QOp}\times R_{\QOp} }.
\end{equation*}

Form the pullback diagram
\begin{equation*}
\xymatrix{ \POp\ar@/_8mm/[ddr]_{(\eta,g)}\ar[r]^{\sim}_{\eta}\ar@{.>}[dr]^(0.6){\sim} & R_{\POp}\ar@/^/[dr]^{\sim} & \\
& T\ar@{.>}[r]^(0.3){\sim}\ar@{.>>}[d] & R_{\POp}\times_{R_{\QOp}} P(R_{\QOp})\ar@{->>}[d]^{(\id,d_0)} \\
& R_{\POp}\times\QOp\ar[r]_{\id\times\eta}^{\sim} & R_{\POp}\times R_{\QOp} }.
\end{equation*}
The morphism $\id\times\eta$ is a weak equivalence by right properness (of the underlying category)
and so is its pullback along $(\id,d_0)$.
By the universal definition of pullbacks,
we have again a natural morphism $\POp\rightarrow T$.
This morphism forms a weak equivalence by the two-out-of-three axiom.
Since $R_{\POp}$ is fibrant,
the natural morphism $(*,\id): R_{\POp}\times\QOp\rightarrow\QOp$
forms a fibration.

Finally,
our construction gives a factorization of $g: \POp\rightarrow\QOp$
of the form
\begin{equation*}
\POp\xrightarrow{\sim} T\twoheadrightarrow R_{\POp}\times\QOp\twoheadrightarrow\QOp.
\end{equation*}
The conclusion follows.
\end{proof}

This lemma, together with assertion~(\ref{EverythingSmallObjectArgument})
of Fact~\ref{PropPathObjectArgument:SmallObjectArgument},
is the main ingredient of the proof of:

\begin{thm}\label{PropPathObjectArgument:SimplicialDGCase}
If the base category $\C$
is either a category of dg-modules over a ring $\kk$ such that $\QQ\subset\kk$,
or a category of simplicial modules over a ring $\kk$,
or the category of simplicial sets $\Simp$,
then the definition of Theorem~\ref{PropSemiModel:Result}
returns a full model structure on the whole category of props $\Prop$.
\end{thm}

In the proof of this theorem,
we use the adjunction between $\Sigma_*$-biobjects and props $F: \C^{\BCat}\rightleftarrows\Prop :U$
rather than the composite adjunction between $F\phi_!: \C^{\ACat}\rightleftarrows\Prop :\phi^* U$.
Recall that the restriction functor $\phi^*: \C^{\BCat}\rightarrow\C^{\ACat}$
creates weak equivalences and fibrations.
Therefore the definition of Theorem~\ref{PropSemiModel:Result}
amounts to saying that weak equivalences and fibrations
of props are created in the category of $\Sigma_*$-biobjects.

\begin{proof}
Recall that the category of props~$\Prop$
inherits automatically axioms~\ref{M:Completeness}-\ref{M:Retract}
from the category of $\Sigma_*$-biobjects~$\C^{\BCat}$.
The lifting axiom~\ref{M:Lifting}.i
is also tautological from the definition of cofibrations in~$\Prop$.

Let $\I$ (respectively, $\J$)
refer to the generating cofibrations (respectively, acyclic cofibrations) of~$\C^{\BCat}$.
Fact~\ref{PropPathObjectArgument:SmallObjectArgument}
implies by adjunction that the domain $\Free(C)$
of any morphism of free props $\Free(i)$, where $i\in\I$ (respectively, $i\in\J$),
is small with respect to every composite of morphisms in $\C$.
Therefore the set $\Free\I$, and similarly $\Free\J$,
permits the small object argument,
from which we deduce that any morphism $f$ has a factorization $f = p i$ (respectively, $f = q j$)
such that $i$ is a relative $\Free\I$-cell (respectively, $\Free\J$-cell) complex
and $p$ (respectively, $q$)
has the right lifting property with respect to the morphisms of~$\Free\I$ (respectively, $\Free\J$).

The morphism $p$ (respectively, $q$) in the factorization $f = p i$ (respectively, $f = q j$)
forms a fibration (respectively, an acyclic fibration)
since $U(p)$ has the right lifting property
with respect to generating cofibrations (respectively, acyclic cofibrations) by adjunction.
The morphisms of $\Free\I$ (respectively, $\Free\J$)
have the left lifting property
with respect to fibrations (respectively, acyclic fibrations)
by adjunction
and so do the retracts of relative $\Free\I$-cell (respectively, $\Free\J$-cell) complexes
since any retract~$f$ of a morphism~$g$
inherits the lifting properties of~$g$.
From these observations,
we deduce that the small object argument, applied to $\Free\I$,
produces the factorization $f = p i$
required by axiom~\ref{M:Factorization}.i,
but we have still to prove that a relative $\Free\J$-cell complex forms a weak equivalence
to prove the existence of the factorization required by axiom~\ref{M:Factorization}.ii.

To get round difficulties,
we use the fibrant replacement functors $R_{\POp}$ and the path objects $P(R_{\POp})$
supplied by facts~\ref{PropPathObjectArgument:FibrantReplacement}-\ref{PropPathObjectArgument:PathObjects}.
By Lemma~\ref{PropPathObjectArgument:FibrationFactorization},
any morphism has a factorization $f = r u$
such that $r$ is a fibration and $u$ is a weak equivalence.
Apply the first factorization axiom~\ref{M:Factorization}.i to $g$
to obtain a factorization $u = p i$,
where $p$ is an acyclic fibration and $i$ is a cofibration.
By the two-out-of-three axiom~\ref{M:TwoOutofThree},
the morphism $i$ is an acyclic cofibration as well.
Hence the identity $f = r u = r p i$
gives a factorization $f = q j$
such that $q = r p$ is a fibration
and $j = i$ is an acyclic cofibration.

We can now prove that the morphisms $i: A\rightarrow B$
which have the left lifting property with respect to cofibrations
are acyclic cofibrations:
apply~\ref{M:Factorization}.ii to get a factorization $i = q j$,
where $q$ is a fibration and $j$ is an acyclic cofibration,
pick a lifting
\begin{equation*}
\xymatrix{ A\ar[d]_{i}\ar@{>->}[]!R+<4pt,0pt>;[r]^{j} & X\ar@{->>}[d]^{q} \\ B\ar[r]_{=}\ar@{.>}[ur] & B }
\end{equation*}
to check that $i$ forms a retract of $j$,
and use~\ref{M:Retract} to conclude.
We deduce from this argument that retracts of relative $\Free\J$-cell complexes
are acyclic cofibrations.

It remains to prove the second lifting axiom~\ref{M:Lifting}.ii.
Let $i: A\rightarrow B$ be an acyclic cofibration.
Apply the small object argument to obtain a factorization $i = q j$,
where $q$ is a fibration and $j$ is a relative $\Free\J$-cell complex.
Recall that $j$ has the left lifting property
with respect to fibrations.
Since we just proved that $j$ forms an acyclic cofibration,
the two-out-of-three axiom~\ref{M:TwoOutofThree}
implies that $q$ forms an acyclic fibration as well.
As in the previous verification,
we apply the first lifting axiom~\ref{M:Lifting}.i
to check that $i$ forms a retract of~$j$.
Then we conclude that $i$ inherits
the left lifting property with respect to fibrations
from $j$.

This verification achieves the proof of Theorem~\ref{PropPathObjectArgument:SimplicialDGCase}.
\end{proof}

The case of topological spaces is more difficult since the domain of generating (acyclic) cofibrations
is not small with respect to every kind of composites.
Therefore, in this case, we cannot completely avoid an analysis
of relative cell complexes in the spirit of~\S\ref{PropPushoutDecomposition}.
The result of this analysis gives:

\begin{lemm}\label{PropPathObjectArgument:TopologicalPushouts}
Let $\K$ be the class of generating (acyclic) cofibrations of $\Sigma_*$-biobjects
in topological spaces.
The components $g: \POp(m,n)\rightarrow\QOp(m,n)$
of a relative $F\K$-cell complex of props in topological spaces
are topological inclusions.
\end{lemm}

\begin{proof}See~\S\ref{PropPushoutDecomposition:Topology}.\end{proof}

Once this lemma is proved,
the line of argument of Theorem~\ref{PropPathObjectArgument:SimplicialDGCase}
can be applied and we obtain:

\begin{thm}\label{PropPathObjectArgument:TopologicalCase}
The definition of Theorem~\ref{PropSemiModel:Result}
also returns a full model structure on the whole category of props
in topological spaces.\qed
\end{thm}

\part{Homotopy invariance of structures}\label{HomotopyEquivalences}
The objective of this part is to prove Theorem~\ref{Mainresult:Transfer} (transfer of structures)
and Theorem~\ref{Mainresult:Equivalence} (equivalences of structures).

First of all,
we review the definition of the endomorphism prop~$\End_X$
of an object~$X$ in a symmetric monoidal category~$\C$.
The structure of an algebra over a prop $\POp$
is defined by a pair $(A,\phi)$,
where $A$ is an object of~$\C$ and $\phi$ is a prop morphism $\phi: \POp\rightarrow\End_A$.
From this definition,
it is tautological that the endomorphism prop $\End_X$ of an object $X\in\C$
is the universal prop acting on $X$.
The definition of the endomorphism prop has a natural generalization
for diagrams of objects $\{X_i\}_{i\in\ICat}\in\C^{\ICat}$.
In the context of diagrams,
the existence of a morphism $\phi: \POp\rightarrow\End_{\{X_i\}_{i\in\ICat}}$
implies that each object $X_i\in\C$ is provided with a $\POp$-algebra structure
in such a way that the map $i\mapsto X_i$
defines a diagram in the category of $\POp$-algebras.
The idea of the proof of Theorem~\ref{Mainresult:Transfer} and Theorem~\ref{Mainresult:Equivalence}
is to use such endomorphism props to construct the demanded diagrams of $\POp$-algebras.
We introduce new axioms of symmetric monoidal model categories, the limit monoid axioms,
which ensure that endomorphism props satisfy certain homotopy invariance properties
necessary to apply lifting arguments.

The background is explained in~\S\ref{PropAlgebras}.
In~\S\ref{StructureTransfer}, we prove the first assertion of Theorem~\ref{Mainresult:Transfer}.
In~\S\ref{PropAlgebraEquivalences}, we prove Theorem~\ref{Mainresult:Equivalence}
and we achieve the proof of Theorem~\ref{Mainresult:Transfer}.

\section{The category of algebras over a prop}\label{PropAlgebras}
The purpose of this section is to give the background of the proof of Theorem~\ref{Mainresult:Transfer} and Theorem~\ref{Mainresult:Equivalence}:
we recall the definition of the endomorphism prop of an object $X$ in a symmetric monoidal category $\C$,
we review the definition of an algebra over a prop,
we define the endomorphism prop of a diagram,
and we introduce the limit monoid axioms
which ensure that endomorphism props satisfy good homotopy invariance properties
when $\C$ is equipped with a model structure.

\subsubsection{Endomorphism props}\label{PropAlgebras:EndomorphismProps}
For any object $X\in\C$,
we form a prop $\End_X$
such that
\begin{equation*}
\End_X(m,n) = \Hom_{\C}(X^{\otimes m},X^{\otimes n}).
\end{equation*}
The symmetric group $\Sigma_m$ (respectively, $\Sigma_n$) operates on $\End_X(m,n)$
by tensor permutations on the source (respectively, on the target).
The horizontal composition product,
respectively the vertical composition product of $\End_X$,
is the natural extension to homomorphisms
of the tensor product of morphisms
\begin{equation*}
\Hom_{\C}(X^{\otimes k},X^{\otimes m})\otimes\Hom_{\C}(X^{\otimes l},X^{\otimes n})
\xrightarrow{\circ_h = \otimes}\Hom_{\C}(X^{\otimes(k+l)},X^{\otimes(m+n)}),
\end{equation*}
respectively the composition product of morphisms
\begin{equation*}
\Hom_{\C}(X^{\otimes k},X^{\otimes n})\otimes\Hom_{\C}(X^{\otimes m},X^{\otimes k})
\xrightarrow{\circ_v = \circ}\Hom_{\C}(X^{\otimes m},X^{\otimes n}).
\end{equation*}
The identities of the objects $X^{\otimes n}$
are represented by morphisms $\eta: \unit\rightarrow\Hom_{\C}(X^{\otimes n},X^{\otimes n})$.
These morphisms determine the unit of the prop $\End_X$.
This prop $\End_X$ is the endomorphism prop of $X$.

\subsubsection{Algebras over a prop}\label{PropAlgebras:Definition}
The structure of an algebra over a prop~$\POp$ in a category~$\C$
is defined by a pair $(A,\phi)$,
where $A$ is an object of~$\C$ and $\phi$ is a prop morphism $\phi: \POp\rightarrow\End_A$.
The morphism $\phi: \POp\rightarrow\End_A$
is equivalent to a collection of morphisms $\phi: \POp(m,n)\rightarrow\Hom_{\C}(A^{\otimes m},A^{\otimes n})$
which, in a point-set context, associates an actual operation $p: A^{\otimes m}\rightarrow A^{\otimes n}$
to any element $p\in\POp(m,n)$.
In the sequel,
we say that  $\phi: \POp\rightarrow\End_A$ defines an action of the prop~$\POp$ on~$A$.

In the literature,
an algebra over a prop is often defined as a functor of enriched symmetric monoidal categories $\underline{A}: \POp\rightarrow\C$.
This abstract definition is equivalent to ours:
since $\POp$ has the non-negative integers as objects and since we have $n = 1^{\otimes n}$, for every $n\in\NN$,
such a functor is determined on objects by tensor powers $\underline{A}(n) = A^{\otimes n}$,
where we set $A = \underline{A}(1)\in\C$,
and the functor $\underline{A}: \POp\rightarrow\C$
is determined on homomorphisms by a collection of morphisms
\begin{equation*}
\POp(m,n)\xrightarrow{\underline{A}}\Hom_{\C}(\underline{A}(m),\underline{A}(n))\xrightarrow{=}\Hom_{\C}(A^{\otimes m},A^{\otimes n})
\end{equation*}
that satisfy relations equivalent to the ones of a prop morphism.

\subsubsection{Endomorphism props of diagrams}\label{PropAlgebras:DiagramEndomorphisms}
Let $\{X_i\}_{i\in\ICat}$
be an $\ICat$-diagram in $\C$.
Form the double sequence of ends
\begin{equation*}
\End_{\{X_i\}_{i\in\ICat}}(m,n) = \int_{i\in\ICat}\Hom_{\C}(X_i^{\otimes m},X_i^{\otimes n}),\quad\text{$m,n\in\NN$}.
\end{equation*}
By definition of an end,
the object $\End_{\{X_i\}_{i\in\ICat}}(m,n)$ can also be defined by a coreflexive equalizer
\begin{equation*}
\xymatrix{ \End_{\{X_i\}_{i\in\ICat}}(m,n)\ar@{.>}[r]^(0.4){\eta} &
\prod_{i\in\ICat}\Hom_{\C}(X_i^{\otimes m},X_i^{\otimes n})\ar@<+2pt>[r]^(0.4){d^0}\ar@<-2pt>[r]_(0.4){d^1} &
\prod_{u: i\rightarrow j}\Hom_{\C}(X_i^{\otimes m},X_j^{\otimes n})\ar@/_2em/[l]_{s^0} },
\end{equation*}
where $u: i\rightarrow j$ ranges over morphisms of~$\ICat$.

The morphism $d^0$ is the cartesian product of the morphisms
\begin{equation*}
\Hom_{\C}(X_i^{\otimes m},X_i^{\otimes n})\xrightarrow{u_*}\Hom_{\C}(X_i^{\otimes m},X_j^{\otimes n})
\end{equation*}
induced by the morphisms $u: i\rightarrow j$ of~$\ICat$.
The morphism $d^1$ is the cartesian product of the morphisms
\begin{equation*}
\Hom_{\C}(X_j^{\otimes m},X_j^{\otimes n})\xrightarrow{u^*}\Hom_{\C}(X_i^{\otimes m},X_j^{\otimes n}).
\end{equation*}
The morphism
\begin{equation*}
\prod_{u: i\rightarrow j}\Hom_{\C}(X_i^{\otimes m},X_j^{\otimes n})\xrightarrow{s^0}\prod_{i\in\ICat}\Hom_{\C}(X_i^{\otimes m},X_i^{\otimes n})
\end{equation*}
is the projection onto the summands associated to identity morphisms $\id: i\rightarrow i$, $i\in\ICat$.

The next proposition is the easy generalization in the context of props of~\cite[Proposition 4.1.2]{Rezk}:

\begin{prop}\hspace*{2mm}\label{PropAlgebras:DiagramEndomorphismProp}
\begin{enumerate}
\item\label{DiagramEndomorphismProp:Structure}
The $\Sigma_*$-biobject $\End_{\{X_i\}_{i\in\ICat}} = \{\End_{\{X_i\}_{i\in\ICat}}(m,n)_{(m,n)\in\NN\times\NN}\}$
inherits a prop structure so that the universal morphism
\begin{equation*}
\eta: \End_{\{X_i\}_{i\in\ICat}}\rightarrow\prod_i\End_{X_i}
\end{equation*}
is a prop morphism.
\item\label{DiagramEndomorphismProp:Application}
Let $\POp$ be a prop.
Suppose we have a collection of morphisms $\phi_i: \POp\rightarrow\End_{X_i}$, $i\in\ICat$,
which provides each object $X_i\in\C$
with the structure of a $\POp$-algebra.
The product morphism $(\phi_i)_i: \POp\rightarrow\prod_i\End_{X_i}$
factors through $\End_{\{X_i\}_{i\in\ICat}}$
if and only if each morphism $u_*: X_i\rightarrow X_j$
defines a morphism of $\POp$-algebras
\begin{equation*}
u_*: (X_i,\phi_i)\rightarrow(X_j,\phi_j).\qed
\end{equation*}
\end{enumerate}
\end{prop}

The prop $\End_{\{X_i\}_{i\in\ICat}}$
is the endomorphism prop of the diagram $\{X_i\}_{i\in\ICat}$.
In the situation of~(\ref{DiagramEndomorphismProp:Application}),
we obtain that the $\POp$-algebras $(X_i,\phi_i)$, $i\in\ICat$,
form an $\ICat$-diagram in the category of $\POp$-algebras.

Note that endomorphism props of diagrams are natural with respect to indexing categories:

\begin{obsv}\label{PropAlgebras:DiagramEndomorphismPropFunctoriality}
Let $\alpha: \ICat\rightarrow\JCat$ be a functor between small categories.
For any $\JCat$-diagram $\{X_j\}_{j\in\JCat}$,
we have a morphism of props
\begin{equation*}
\alpha^*: \End_{\{X_j\}_{j\in\JCat}}\rightarrow\End_{\{X_{\alpha(i)}\}_{i\in\ICat}}
\end{equation*}
such that the diagram
\begin{equation*}
\xymatrix{ \End_{\{X_j\}_{j\in\JCat}}\ar@{.>}[d]_{\alpha^*}\ar[r] & \prod_{j\in\JCat}\End_{X_j}\ar[d]^{\alpha^*} \\
\End_{\{X_{\alpha(i)}\}_{i\in\ICat}}\ar[r] & \prod_{i\in\ICat}\End_{X_{\alpha(i)}} }
\end{equation*}
commutes,
where the right-hand side vertical morphism is given by the projection onto factors
associated to indices $j = \alpha(i)$, $i\in\ICat$.

The map $\alpha\mapsto\alpha^*$
satisfies the functoriality relations $\id^* = \id$ and $(\alpha\beta)^{*} = \beta^*\alpha^*$.
\end{obsv}

\subsubsection{Limit monoid axioms}\label{PropAlgebras:LimitMonoidAxioms}
In the sequel,
we use the notation $\Hom_{X Y}$
to refer to the double sequence of hom-objects
\begin{equation*}
\Hom_{X Y}(m,n) = \Hom_{\C}(X^{\otimes m},Y^{\otimes n}),
\end{equation*}
for any $X,Y\in\C$.
The map $(X,Y)\mapsto\Hom_{X Y}$
defines clearly a functor in $X,Y\in\C$.
We have by definition $\End_X = \Hom_{X X}$, for any $X\in\C$.

In the next sections,
our main task is to analyze the homotopy invariance of certain endomorphism props of diagrams.
For this purpose,
we arrange the end of the definition to decompose the endomorphism prop $\End_{\{X_i\}_{i\in\ICat}}$
into manageable pullbacks formed from objects of the form $\Hom_{X Y}$.
The idea is to apply the adjoint form of the pushout-product axiom
to prove that certain functors between small categories $\alpha: \ICat\rightarrow\JCat$
induce (acyclic) fibrations on endomorphism props.
But the usual axioms of symmetric monoidal model categories
are not sufficient to prove that the target of~$\Hom_{X Y}(m,n) = \Hom_{\C}(X^{\otimes m},Y^{\otimes n})$
is a functor in $Y$ which preserves fibrations,
because this is not the case of the tensor power $Y^{\otimes n}$ in general.
Therefore,
we introduce additional axioms
to ensure such properties:
\begin{enumerate}\renewcommand{\theenumi}{LM\arabic{enumi}}
\item (\emph{final monoid axiom}):
The natural morphism $*\otimes *\rightarrow *$, where $*$ denotes the terminal object of $\C$,
is an isomorphism.
\item (\emph{cartesian monoid axiom}):
If we have a fibration of the form $(f,g): S\rightarrow X\times_B Y$,
then the morphism
$(f\otimes Z,g\otimes Z): S\otimes Z\rightarrow X\otimes Z\times_{B\otimes Z} Y\otimes Z$
arising from the diagram
\begin{equation*}
\xymatrix{ S\otimes Z\ar@{.>}[dr]^(0.6){(f\otimes Z,g\otimes Z)}\ar@/_6mm/[ddr]_{f\otimes Z}\ar@/^6mm/[rrd]^{g\otimes Z} && \\
& X\otimes Z\times_{B\otimes Z} Y\otimes Z\ar@{.>}[d]\ar@{.>}[r] & Y\otimes Z\ar[d] \\
& X\otimes Z\ar[r] & B\otimes Z }
\end{equation*}
forms a fibration too,
for any fibrant object $Z\in\C$.
\end{enumerate}
We say that $\C$ satisfies the limit monoid axioms when these requirements are satisfied.

The category of topological spaces, simplicial sets,
and, more generally, any symmetric monoidal category whose tensor product is given by the cartesian product,
satisfy the limit monoid axioms.
So do the categories of dg-modules
and the categories of simplicial modules over a ring,
because tensor products preserve cokernels
and, in the case $\C = \dg\kk\Mod,\s\kk\Mod$,
fibrations are morphisms which are surjective (over connected components).

\medskip
The next proposition follows from an easy inspection:

\begin{prop}\label{PropAlgebras:FibrantTensors}
The next properties hold in any symmetric monoidal model category which satisfies the limit monoid axioms.
\begin{enumerate}
\item
If $X$ is fibrant, then so is the object $X^{\otimes n}$, for every $n\in\NN$.
\item
If $p: S\rightarrow X$ is a fibration and $X$ is fibrant,
then $f^{\otimes n}: S^{\otimes n}\rightarrow X^{\otimes n}$
is a fibration as well.
\item
If $p: S\rightarrow X\times Y$ is a fibration and $X$ and $Y$ are fibrant objects,
then $f^{\otimes n}: S^{\otimes n}\rightarrow X^{\otimes n}\times Y^{\otimes n}$
is a fibration as well.\qed
\end{enumerate}
\end{prop}

\subsubsection{Conventions}
From now on,
we assume that $\C$ is a symmetric monoidal model category
which satisfies the limit monoid axioms.
The letter $\POp$
refers to a cofibrant prop in~$\C$.

We restrict ourselves to the full subcategory of props with non-empty inputs (outputs)
if this is necessary to have a semi-model structure on props.
The components $\End_{X}(0,n)$ of an endomorphism prop are nontrivial,
but endomorphism props occur always as the targets of morphisms $\phi: \POp\rightarrow\End_{X}(0,n)$
and for that reason can be replaced by the associated prop with non-empty inputs $\overline{\End}_{X}$
if necessary.
To simplify,
we do not mark this replacement in our notation.

\section{Transfer of structures}\label{StructureTransfer}
The purpose of this section is to prove the first assertion of Theorem~\ref{Mainresult:Transfer}.
The idea is to use the endomorphism prop of the diagram
\begin{equation*}
\{A\xrightarrow{f} B\}
\end{equation*}
formed by a single morphism $f$.
Let
\begin{equation*}
\End_{A}\xleftarrow{d_1}\End_{\{A\rightarrow B\}}\xrightarrow{d_0}\End_{B}
\end{equation*}
denote the components of the universal morphism $\eta: \End_{\{A\rightarrow B\}}\rightarrow\End_A\times\End_B$.
We prove that $(d_0,d_1)$ satisfies good homotopy properties
in order to apply lifting arguments in the category of props
and to transport a prop action from~$B$ to~$A$.
We use the following presentation of the endomorphism prop $\End_{\{A\rightarrow B\}}$
to obtain our result:

\begin{prop}\label{StructureTransfer:MorphismProp}
The endomorphism prop $\End_{\{A\rightarrow B\}}$ of a morphism $f: A\rightarrow B$
fits in a pullback:
\begin{equation*}
\xymatrix{ \End_{\{A\rightarrow B\}}\ar@{.>}[r]^{d_0}\ar@{.>}[d]_{d_1} & \End_{B}\ar[d]^{f^*} \\
\End_{A}\ar[r]_{f_*} & \Hom_{A B} }
\end{equation*}
where $d_0,d_1$ are prop morphisms.
\end{prop}

\begin{proof}Exercise.\end{proof}

This proposition is an obvious generalization of~\cite[Proposition 4.1.4]{Rezk}.

\begin{lemm}\label{StructureTransfer:MorphismPropHomotopy}
Suppose $A$ and $B$ in the construction of Proposition~\ref{StructureTransfer:MorphismProp} are both cofibrant and fibrant.
Then the endomorphism props $\End_A$ and $\End_B$
are fibrant.
Moreover:
\begin{enumerate}
\item\label{MorphismPropHomotopy:AcyclicFibration}
If $f$ is an acyclic fibration,
then so is $d_0$.
\item\label{MorphismPropHomotopy:AcyclicCofibration}
If $f$ is an acyclic cofibration,
then $d_0$ is a weak equivalence and $d_1$ is an acyclic fibration.
\end{enumerate}
\end{lemm}

This lemma is a generalization of~\cite[Proposition 4.1.7-8]{Rezk}
in the context of props.

\begin{proof}[Proof of assertion (\ref{MorphismPropHomotopy:AcyclicFibration})]
The limit monoid axiom implies that $A^{\otimes n}$ and $B^{\otimes n}$
are both fibrant by Proposition~\ref{PropAlgebras:FibrantTensors}.
The pushout-product axiom of symmetric monoidal model categories implies that $A^{\otimes m}$ and $B^{\otimes m}$
are both cofibrant.
The adjoint axiom~\ref{MM:PullbackHom}
implies that the objects $\Hom_{\C}(A^{\otimes m},A^{\otimes n})$
(respectively, $\Hom_{\C}(B^{\otimes m},B^{\otimes n})$)
are fibrant for every $(m,n)\in\NN\times\NN$,
from which we conclude that $\End_{A}$ (respectively, $\End_{B}$)
is a fibrant object in the category of props.

If $f$ is a fibration,
then so is $f^{\otimes n}: A^{\otimes n}\rightarrow B^{\otimes n}$
by Proposition~\ref{PropAlgebras:FibrantTensors}.
If $f$ is a weak equivalence and the objects $A$ and $B$ are cofibrant,
then $f^{\otimes n}: A^{\otimes n}\rightarrow B^{\otimes n}$
forms a weak equivalence too,
because the pushout-product axiom
and Brown's lemma (see~\cite[Lemma 1.1.12]{Hovey})
imply that the tensor product with a cofibrant object preserves weak equivalences
between cofibrant objects.
If $f$ is an acyclic fibration,
then we obtain that $f^{\otimes n}: A^{\otimes n}\rightarrow B^{\otimes n}$
forms an acyclic fibration as well.

Since we have already observed that $A^{\otimes m}$ is cofibrant,
axiom~\ref{MM:PullbackHom}
implies that $(f^{\otimes n})_*: \Hom_{\C}(A^{\otimes m},A^{\otimes n})\rightarrow\Hom_{\C}(A^{\otimes m},B^{\otimes n})$
forms an (acyclic) fibration
if $f^{\otimes n}$ does.
Since this assertion holds for every pair $(m,n)\in\NN\times\NN$,
we deduce that $f_*: \Hom_{A A}\rightarrow\Hom_{A B}$
forms a fibration in the category of $\Sigma_*$-biobjects.
Since the class of acyclic fibrations is closed under pullbacks,
we obtain that $d_0$ forms an acyclic fibration as well.
This proves assertion~(\ref{MorphismPropHomotopy:AcyclicFibration}).
Just recall that the forgetful functor from props to $\Sigma_*$-biobjects
creates weak equivalences and fibrations.
\end{proof}

\begin{proof}[Proof of assertion (\ref{MorphismPropHomotopy:AcyclicCofibration})]
We have already observed that $A^{\otimes n}$ and $B^{\otimes n}$
are both cofibrant and fibrant for every $n\in\NN$.
Axiom~\ref{MM:PullbackHom}
and the Brown lemma imply that the functor $\Hom_{\C}(A^{\otimes m},-)$
preserves weak equivalences between fibrant objects.
Use again that $f^{\otimes n}: A^{\otimes n}\rightarrow B^{\otimes n}$
forms a weak equivalence
to deduce that $(f^{\otimes n})_*$
is still a weak equivalence.
Since this assertion holds for every pair $(m,n)\in\NN\times\NN$,
we deduce that $f_*: \Hom_{A A}\rightarrow\Hom_{A B}$
forms a weak equivalence in the category of $\Sigma_*$-biobjects.

Since $B^{\otimes n}$ is fibrant,
axiom~\ref{MM:PullbackHom},
implies that $(f^{\otimes m})^*: \Hom_{\C}(B^{\otimes m},B^{\otimes n})\rightarrow\Hom_{\C}(A^{\otimes m},B^{\otimes n})$
is an acyclic fibration for every $(m,n)\in\NN\times\NN$.
Hence, the morphism $f^*: \Hom_{B B}\rightarrow\Hom_{B A}$
forms an acyclic fibration of $\Sigma_*$-biobjects.
Since the class of acyclic fibrations is closed under pullbacks,
we obtain that $d_1$ forms an acyclic fibration as well.

Thus the morphisms $d_1,f_*,f^*$ are all weak equivalences.
We conclude from the two-out-of-three axiom
that $d_0$ is a weak equivalence as well.
\end{proof}

\begin{thm}\label{StructureTransfer:TransferConstruction}
Let $\POp$ be a cofibrant prop.
Let $f: A\rightarrow B$
be a morphism in $\C$
such that $A$ and $B$ are both cofibrant and fibrant
in $\C$.
Suppose we have a prop morphism $\psi: \POp\rightarrow\End_B$
which provides $B$ with a $\POp$-algebra structure.
\begin{enumerate}
\item\label{TransferConstruction:AcyclicFibration}
If $f$ is an acyclic fibration,
then we have a prop morphism $\phi: \POp\rightarrow\End_A$
such that $f$ defines a morphism of $\POp$-algebras
\begin{equation*}
(A,\phi)\xrightarrow{f}(B,\psi).
\end{equation*}
\item\label{TransferConstruction:AcyclicCofibration}
If $f$ is an acyclic cofibration,
then we have prop morphisms $\phi: \POp\rightarrow\End_A$ and $\theta: \POp\rightarrow\End_B$
such that $f$ defines a morphism of $\POp$-algebras
\begin{equation*}
(A,\phi)\xrightarrow{f}(B,\theta)
\end{equation*}
and $\theta$ is homotopic to $\psi$ in the semi-model category of props.
\end{enumerate}
\end{thm}

This theorem is a generalization of~\cite[Theorem 3.5]{BergerMoerdijk}
in the context of props.

\begin{proof}[Proof of assertion (\ref{TransferConstruction:AcyclicFibration})]
We apply Lemma~\ref{StructureTransfer:MorphismPropHomotopy} to the morphism $f: A\rightarrow B$.

If $f$ is an acyclic fibration, then we have a diagram
\begin{equation*}
\xymatrix{ & \End_{\{A\rightarrow B\}}\ar@{.>}[r]^{d_1}\ar@{.>>}[d]_{\sim}^{d_0} & \End_{A}\ar@{->>}[d]_{\sim}^{f_*} \\
\POp\ar[r]_{\phi}\ar@{.>}[ur]^{\tilde{\phi}} & \End_{B}\ar[r]_{f^*} & \Hom_{A B} }
\end{equation*}
such that $d_0$ is an acyclic fibration.
Hence, by the lifting axiom of semi-model categories,
we have a morphism $\tilde{\phi}$ such that $d_0\tilde{\phi} = \phi$,
from which we deduce the existence of a diagram
\begin{equation*}
(A,\phi)\xrightarrow{f}(B,\psi)
\end{equation*}
in the category of $\POp$-algebras,
where we set $\phi = d_1\tilde{\phi}$
to get the $\POp$-algebra structure of~$A$.
This achieves the proof of assertion~(\ref{TransferConstruction:AcyclicFibration}).
\end{proof}

\begin{proof}[Proof of assertion (\ref{TransferConstruction:AcyclicCofibration})]
If $f$ is weak equivalence, then we have a diagram
\begin{equation*}
\xymatrix{ \End_{\{A\rightarrow B\}}\ar@{.>}[r]^{d_0}_{\sim}\ar@{.>>}[d]_{\sim}^{d_1} & \End_{B}\ar@{->>}[d]_{\sim}^{f^*} \\
\End_{A}\ar[r]^{f_*}_{\sim} & \Hom_{A B} }
\end{equation*}
such that $d_0$ is a weak equivalence and $d_1$ is an acyclic fibration.
Recall also that $\End_A$ and $\End_B$
are both fibrant,
and as a byproduct, so is the endomorphism prop $\End_{\{A\rightarrow B\}}$.
As $\POp$ is also cofibrant,
we obtain that $d_0$ induces a bijection between homotopy classes
of prop morphisms:
\begin{equation*}
[\POp,\End_{\{A\rightarrow B\}}]\xrightarrow[d_0]{\simeq}[\POp,\End_{B}].
\end{equation*}
Consequently,
we have a prop morphism $\tilde{\phi}: \POp\rightarrow\End_{\{A\rightarrow B\}}$
such that $\theta = d_0\tilde{\phi}$
is left homotopic to $\psi$.
Take the diagram of $\POp$-algebras determined by this morphism
to reach the conclusion of assertion~(\ref{TransferConstruction:AcyclicCofibration}).
\end{proof}

For a weak equivalence $f: A\xrightarrow{\sim} B$,
which is not necessarily an acyclic cofibration nor an acyclic fibration,
we apply the factorization axiom and the two-out-of-three axiom of model categories
to obtain a factorization $f = p i$
such that $i: A\rightarrow Z$ is an acyclic cofibration and $p: Z\rightarrow B$ is an acyclic fibration.
Theorem~\ref{StructureTransfer:TransferConstruction}
applied to $p$ and $i$
gives the following result:

\begin{prop}\label{StructureTransfer:TransferDecomposition}
We have prop morphisms $\phi: \POp\rightarrow\End_A$ and $\rho,\sigma: \POp\rightarrow\End_Z$
such that $\rho$ is homotopic to $\sigma$,
the morphism $i$ defines a weak equivalence of $\POp$-algebras $i: (A,\phi)\xrightarrow{\sim}(Z,\rho)$
and $p$ defines a weak equivalence of $\POp$-algebras $p: (Z,\sigma)\xrightarrow{\sim}(B,\psi)$.
\end{prop}

Hence to achieve the proof of Theorem~\ref{Mainresult:Transfer}
it remains to prove that $\POp$-algebras $(Z,\rho)$ and $(Z,\sigma)$
associated to homotopic prop morphisms $\rho,\sigma: \POp\rightrightarrows\End_Z$
are connected by weak equivalences of~$\POp$-algebras
$(Z,\rho)\xleftarrow{\sim}\cdot\xrightarrow{\sim}(Z,\sigma)$.
This verification is the goal of the next section.

\section{Prop homotopies and naive equivalences of algebras over props}\label{PropAlgebraEquivalences}
The purpose of this section is to prove Theorem~\ref{Mainresult:Equivalence}
and to achieve the proof of Theorem~\ref{Mainresult:Transfer}.
Suppose we have homotopic prop morphisms
$\phi^0,\phi^1: \POp\rightrightarrows\End_A$
that provide $A$ with a $\POp$-algebra structure.

By a standard assertion of the theory of model categories,
a left homotopy in the category of props, determined by a morphism on a cylinder object,
is equivalent to a right homotopy, determined by a morphism towards a path object.
Our idea is to form a path object $A\rightarrow Z\rightrightarrows A$
in the underlying category and to use the endomorphism prop of this diagram
as a good approximation of the path object of $\End_A$.
Then we adapt the standard construction of the equivalence between right and left homotopies
to form a diagram of~$\POp$-algebras
$(A,\phi^0)\xleftarrow{\sim}(Z,\psi)\xrightarrow{\sim}(A,\phi^1)$
as asserted by Theorem~\ref{Mainresult:Equivalence}.

\subsubsection{Path-objects, the $\YCat$ and $\VCat$ diagrams}\label{PropAlgebraEquivalences:Diagrams}
From now on,
we assume that $A$ is a cofibrant and fibrant object of $\C$.
We form a factorization $\displaystyle A\cofib^{\sim} Z\fib A\times A$
of the diagonal $A\xrightarrow{(\id,\id)} A\times A$.
We adopt the notation
\begin{equation*}
\SCopy_{= A}\quad\cofib^{s_0}_{\sim}\quad Z\fib^{(d_0,d_1)}\TzeroCopy_{= A}\times\ToneCopy_{= A}
\end{equation*}
to refer unambiguously to the copies of~$A$
which sit in this factorization.
Since $T_0 = T_1 = A$ is fibrant, we obtain that $d_0: Z\rightarrow T_0$ and $d_1: Z\rightarrow T_1$
are both acyclic fibrations.

Consider the $\YCat$-diagram
\begin{equation*}
\YCat = \left\{\vcenter{\xymatrix@!@M=2pt@H=6pt@W=6pt@R=4pt@C=12pt{ && T_0 \\
S\ar@{>->}[]!R+<4pt,0pt>;[r]_(0.4){\sim}^(0.4){s_0} & Z\ar@{->>}[ur]_{\sim}^{d_0}\ar@{->>}[dr]^{\sim}_{d_1} & \\
&& T_1 }}\right\}
\end{equation*}
and the $\VCat$-diagrams
\begin{equation*}
\{T_0\leftarrow Z\rightarrow T_1\} = \left\{\vcenter{\xymatrix@!@M=2pt@H=6pt@W=6pt@R=4pt@C=12pt{ & T_0 \\
Z\ar@{->>}[ur]_{\sim}^{d_0}\ar@{->>}[dr]^{\sim}_{d_1} & \\ & T_1 }}\right\}
\quad\text{and}
\quad\{T_0 = S = T_1\} = \left\{\vcenter{\xymatrix@!@M=2pt@H=6pt@W=6pt@R=4pt@C=12pt{ & T_0 \\
S\ar[ur]^{\id}\ar[dr]_{\id} & \\ & T_1 }}\right\}.
\end{equation*}
Our constructions involve the diagram embeddings
\begin{equation*}
\xymatrix@!0@M=2pt@R=1.5cm@C=3cm{ \{S\}\ar@{^{(}->}[]!R+<4pt,0pt>;[r]^(0.22){i} &
\{T_0\xleftarrow{=} S\xrightarrow{=} T_1\}\ar@{^{(}->}[]!R+<4pt,0pt>;[r]^(0.45){v} &
\YCat & \\
&
\{T_0,T_1\}\ar@{_{(}->}[]!U+<0pt,4pt>;[u]^(0.4){t}\ar@{^{(}->}[]!R+<4pt,0pt>;[r]_(0.22){u} &
\{T_0\leftarrow Z\rightarrow T_1\}\ar@{_{(}->}[]!U+<0pt,4pt>;[u]^(0.35){w} &
\{Z\}\ar@{_{(}->}[]!L-<4pt,0pt>;[l]^(0.22){j} }
\end{equation*}
and the induced prop morphisms
\begin{equation*}
\xymatrix@!0@M=2pt@R=1.5cm@C=3cm{ \End_{S}\ar@/_8pt/@{.>}[dr]_{t^* (i^*)^{-1}} &
\ar[l]_(0.55){i^*}^(0.55){\simeq}\End_{\{T_0 = S = T_1\}}\ar[d]_(0.4){t^*} &
\ar[l]_(0.35){v^*}\End_{\YCat}\ar[d]^(0.4){w^*} & \\
&
\End_{\{T_0,T_1\}} &
\ar[l]^(0.52){u^*}\End_{\{T_0\leftarrow Z\rightarrow T_1\}}\ar[r]_(0.62){j^*} &
\End_{Z} }
\end{equation*}
Note that $\End_{\{T_0,T_1\}} = \End_{T_0}\times\End_{T_1} = \End_A\times\End_A$,
the morphism $i^*$ defines an isomorphism $\End_{\{T_0 = S = T_1\}}\simeq\End_{S} = \End_A$,
and the composite $t^* (i^*)^{-1}$
is identified with the diagonal morphism $\Delta: \End_A\rightarrow\End_A\times\End_A$
in the category of props.

\begin{lemm}\hspace*{2mm}\label{PropAlgebraEquivalences:VDiagramHomotopy}
\begin{enumerate}
\item
The endomorphism prop of the $\VCat$-diagram $\{T_0\leftarrow Z\rightarrow T_1\}$
fits in a pullback:
\begin{equation*}
\xymatrix{ \End_{\{T_0\leftarrow Z\rightarrow T_1\}}\ar@{.>}[d]_{u^*}\ar@{.>}[rrrr]^{j^*} &&&& \End_{Z}\ar[d]^{((d_0)_*,(d_1)_*)} \\
\End_{\{T_0,T_1\}}\ar[rr]_(0.45){=} && \End_{T_0}\times\End_{T_1}\ar[rr]_(0.45){(d_0)^*\times(d_1)^*} && \Hom_{Z T_0}\times\Hom_{Z T_1} }
\end{equation*}
\item
The morphism $u^*: \End_{\{T_0\leftarrow Z\rightarrow T_1\}}\rightarrow\End_{\{T_0,T_1\}}$
is a fibration.
\end{enumerate}
\end{lemm}

\begin{proof}
The first assertion follows from an easy exercise.

Observe that $Z$ is cofibrant if $S = A$ is.
The pushout-product axiom implies that $Z^{\otimes m}$
forms a cofibrant object.
The limit monoid axioms imply that $((d_0)_*,(d_1)_*): Z^{\otimes n}\rightarrow T_0^{\otimes n}\times T_1^{\otimes n}$
forms a fibration.
By~(\ref{MM:PullbackHom}),
we obtain that the morphism
\begin{multline*}
\Hom_{\C}(Z^{\otimes m},Z^{\otimes n})
\xrightarrow{((d_0)_*,(d_1)_*)}\Hom_{\C}(Z^{\otimes m},T_0^{\otimes n}\times T_1^{\otimes n})\\
= \Hom_{\C}(Z^{\otimes m},T_0^{\otimes n})\times\Hom_{\C}(Z^{\otimes m},T_1^{\otimes n})
\end{multline*}
forms a fibration, for every $(m,n)\in\NN\times\NN$.
Hence the morphism $((d_0)_*,(d_1)_*): \End_{Z}\rightarrow\Hom_{Z T_0}\times\Hom_{Z T_1}$
forms a fibration of $\Sigma_*$-biobjects.
Use that the class of fibrations is stable under pullbacks to conclude.
\end{proof}

\begin{lemm}\hspace*{2mm}\label{PropAlgebraEquivalences:YDiagramHomotopy}
\begin{enumerate}
\item
The endomorphism prop of the $\YCat$-diagram
is given by a cartesian product of the form:
\begin{equation*}
\End_{\YCat}
= (\End_{S})\bigtimes_{(\Hom_{S Z})}(\End_Z)\bigtimes_{(\Hom_{Z T_0}\times\Hom_{Z T_1})}
(\End_{T_0}\times\End_{T_1})
\end{equation*}
\item
The morphism $v^*: \End_{\YCat}\rightarrow\End_{\{T_0 = S = T_1\}}$
is an acyclic fibration.
\end{enumerate}
\end{lemm}

\begin{proof}
The identity of assertion (1) follows from an easy exercise.

The argument of Lemma~\ref{PropAlgebraEquivalences:VDiagramHomotopy},
applied to the $\VCat$-diagram $\{T_0\xleftarrow{=} S\xrightarrow{=} T_1\}$,
also gives an identity
\begin{equation*}
\End_{\{T_0 = S = T_1\}} = (\End_{S})\bigtimes_{(\Hom_{S T_0}\times\Hom_{S T_1})}
(\End_{T_0}\times\End_{T_1}).
\end{equation*}
The morphism $v^*: \End_{\YCat}\rightarrow\End_{\{T_0=S=T_1\}}$
can be identified with the base extension
\begin{equation*}
(\End_{S})\bigtimes_{(\Hom_{S Z})}\mathbf{\Phi}\bigtimes_{(\Hom_{Z T_0}\times\Hom_{Z T_1})}
(\End_{T_0}\times\End_{T_1})
\end{equation*}
of the natural morphism
\begin{equation*}
\End_{Z}\xrightarrow{\mathbf{\Phi}}(\Hom_{S Z})\bigtimes_{(\Hom_{S T_0}\times\Hom_{S T_1})}
(\Hom_{Z T_0}\times\End_{Z T_1})
\end{equation*}
yielded by the diagram
\begin{equation*}
\xymatrix@!0@M=4pt@R=12mm@C=24mm{\End_Z\ar@{.>}[dr]^{\mathbf{\Phi}}\ar@/^/[drrr]!U^{((d_0)_*,(d_1)_*)}\ar@/_4mm/[ddr]!L_{(s_0)^*} &&& \\
& *+<4pt>{\text{pull}}\ar@{.>}[rr]\ar@{.>}[d] && \Hom_{Z T_0}\times\End_{Z T_1}\ar[d]^{(s_0)^*} \\
& \Hom_{S Z}\ar[rr]_(0.4){((d_0)_*,(d_1)_*)} && \Hom_{S T_0}\times\Hom_{S T_1} }.
\end{equation*}
Since the class of acyclic fibrations is stable under base extension,
it is sufficient to prove that $\Phi$ is one to reach the conclusion of the lemma.

Observe that
\begin{align*}
& \Hom_{S T_0}(m,n)\times\Hom_{S T_1}(m,n) = \Hom_{\C}(S^{\otimes m},T_0^{\otimes n}\times T_1^{\otimes n}) \\
\text{and}\quad & \Hom_{Z T_0}(m,n)\times\Hom_{Z T_1}(m,n) = \Hom_{\C}(Z^{\otimes m},T_0^{\otimes n}\times T_1^{\otimes n}).
\end{align*}
The components of~$\mathbf{\Phi}$
are yielded by diagrams of the form
\begin{equation*}
\xymatrix@!0@M=4pt@R=12mm@C=28mm{\Hom(Z^{\otimes m},Z^{\otimes n})\ar@{.>}[dr]^{\mathbf{\Phi}}
\ar@/^/[drrr]^{(d_0,d_1)_*}\ar@/_4mm/[ddr]!L_{(s_0)^*} &&& \\
& *+<4pt>{\text{pull}}\ar@{.>}[rr]\ar@{.>}[d] && \Hom_{\C}(Z^{\otimes m},T_0^{\otimes n}\times T_1^{\otimes n})\ar[d]^{(s_0)^*} \\
& \Hom_{\C}(S^{\otimes m},Z^{\otimes n})\ar[rr]_{(d_0,d_1)_*} && \Hom_{\C}(S^{\otimes m},T_0^{\otimes n}\times T_1^{\otimes n}) }.
\end{equation*}
Since $S = A$ is supposed to be cofibrant and $s_0: S\rightarrow Z$
is an acyclic cofibration,
the pushout-product axiom implies that $(s_0)^{\otimes n}$
forms an acyclic cofibration as well.
The limit monoid axioms imply that $((d_0)_*,(d_1)_*): Z^{\otimes n}\rightarrow T_0^{\otimes n}\times T_1^{\otimes n}$
forms a fibration.
The assertion that $\mathbf{\Phi}$ forms an acyclic fibration
follows from axiom~(\ref{MM:PullbackHom}).
This verification achieves the proof of the lemma.
\end{proof}

\begin{thm}\label{PropAlgebraEquivalences:Construction}
Let $\POp$ be a cofibrant prop in $\C$.
Suppose $A$ is a cofibrant and fibrant object of $\C$.
Let $\phi^0,\phi^1: \POp\rightarrow\End_A$ be prop morphisms which provide $A$
with a $\POp$-algebra structure each.
If $(\phi^0,\phi^1)$ are left homotopic,
then we have a $\POp$-algebra $(Z,\psi)$
together with weak equivalences
\begin{equation*}
(A,\phi^0)\xleftarrow[d_0]{\sim}(Z,\psi)\xrightarrow[d_1]{\sim}(A,\phi^1)
\end{equation*}
in the category of $\POp$-algebras.
\end{thm}

\begin{proof}
Let
\begin{equation*}
\xymatrix{ \POp\ar@<+2pt>[r]^{d^0}\ar@<-2pt>[r]_{d^1} & \widetilde{\POp}\ar[r]^{s^0}_{\sim} & \POp }
\end{equation*}
be a cylinder object of~$\POp$.
Recall that $(d^0,d^1)$
is supposed to form a cofibration.
The morphisms $d^0$ and $d^1$ taken individually are cofibrations too,
because $\POp$
is supposed to be cofibrant.
These cofibrations $d^0$ and $d^1$
are also acyclic by the two-out-of-three axiom.
By assumption,
we have a lifting in the diagram
\begin{equation*}
\xymatrix{ \POp\vee\POp\ar@{>->}[]!D+<0pt,-4pt>;[d]_(0.35){(d^0,d^1)}\ar[rr]^{(\phi^0,\phi^1)} && \End_{A}\ar[d] \\
\widetilde{\POp}\ar@{.>}[urr]^{h}\ar[rr] && *+<4pt>{*} }
\end{equation*}

Apply (M4.i)
to get a lifting in the diagram
\begin{equation*}
\xymatrix{ I\ar@{>->}[]!D+<0pt,-4pt>;[d]\ar[rr] && \End_{\YCat}\ar@{->>}[d]^(0.4){v^*}_{\sim}\ar[r]^(0.4){w^*} &
\End_{\{T_0\leftarrow Z\rightarrow T_1\}}\ar@{->>}[d]^(0.4){u^*} \\
\POp\ar@{.>}[urr]^{k}\ar[r]_(0.4){\phi^0} & \End_{S} & \End_{\{T_0 = S = T_1\}}\ar[l]_(0.6){\simeq}^(0.6){i^*}\ar[r]_(0.55){t^*} &
\End_{\{T_0,T_1\}} }.
\end{equation*}
The bottom horizontal composite
can be identified with the diagonal morphism
\begin{equation*}
\POp\xrightarrow{(\phi^0,\phi^0)}\End_{A}\times\End_{A} = \End_{\{T_0,T_1\}}.
\end{equation*}
Form now the product morphism
\begin{equation*}
\widetilde{\POp}\xrightarrow{(\phi^0 s^0,h)}\End_{A}\times\End_{A} = \End_{\{T_0,T_1\}}.
\end{equation*}
Check that the solid frame
commutes in the diagram
\begin{equation*}
\xymatrix{ \POp\ar@{>->}[]!D+<0pt,-4pt>;[d]^(0.35){\sim}_(0.35){d^0}\ar[r]^{k} & \End_{\YCat}\ar[r]^(0.4){w^*} &
\End_{\{T_0\leftarrow Z\rightarrow T_1\}}\ar@{->>}[d]^(0.4){u^*} \\
\widetilde{\POp}\ar@{.>}[urr]^{\ell}\ar[rr]_{(\phi^0 s^0,h)} && \End_{\{T_0,T_1\}} },
\end{equation*}
and apply (M4.ii)
to obtain the existence of a lifting $\ell$.

Form the morphism
\begin{equation*}
\POp\xrightarrow{\ell d^1}\End_{\{T_0\leftarrow Z\rightarrow T_1\}}
\end{equation*}
and take $\psi = j^*\ell d^1$
to provide the object $Z$ with a $\POp$-algebra structure.
It is easy to check that the composite
\begin{equation*}
\POp\xrightarrow{\ell d^1}\End_{\{T_0\leftarrow Z\rightarrow T_1\}}\xrightarrow{u^*}\End_{\{T_0,T_1\}} = \End_A\times\End_A
\end{equation*}
agrees with the morphism $(\phi^0,\phi^1)$.
Hence the existence of the morphism $\ell d^1$
implies the existence of a diagram of $\POp$-algebras
\begin{equation*}
(A,\phi^0)\xleftarrow[d_0]{\sim}(Z,\psi)\xrightarrow[d_1]{\sim}(A,\phi^1),
\end{equation*}
from which we draw the conclusion of the theorem.
\end{proof}

This result achieves the proof of Theorem~\ref{Mainresult:Equivalence}
and, together with Proposition~\ref{StructureTransfer:TransferDecomposition},
achieves the proof of Theorem~\ref{Mainresult:Transfer}.\qed

\renewcommand{\thesection}{\Alph{section}}
\setcounter{section}{0}
\part*{Appendix: analysis of pushouts in the category of props}

Recall that we have a composite adjunction
\begin{equation*}
\xymatrix{ \C^{\ACat}\ar@<+2pt>[r]^{\phi_!} & \C^{\BCat}\ar@<+2pt>[r]^{F}\ar@<+2pt>[l]^{\phi^*} & \Prop\ar@<+2pt>[l]^{U} }.
\end{equation*}
The first forgetful functor $U: \Prop\rightarrow\C^{\BCat}$
retains only symmetric group actions on the underlying collection $\{\POp(m,n)\}_{(m,n)}$ of a prop $\POp$.
The restriction functor $\phi^*: \C^{\BCat}\rightarrow\C^{\ACat}$ forgets these actions
and retains only the $\NN\times\NN$-grading.

These adjunctions restrict to subcategories of objects with non-empty inputs
\begin{equation*}
\xymatrix{ \C^{\ACat}_0\ar@<+2pt>[r]^{\phi_!} & \C^{\BCat}_0\ar@<+2pt>[r]^{F}\ar@<+2pt>[l]^{\phi^*} & \Prop_0\ar@<+2pt>[l]^{U} }.
\end{equation*}
The goal of this appendix is to prove that the composite adjunction
$F\phi_!: \C^{\ACat}_0\rightleftarrows\Prop_0 :\phi^* U$
fulfils the requirement of Theorem~\ref{SemiModelCategories:AdjointStructure}.
Then we can conclude that the category of props with non-empty inputs inherits a semi-model structure
as asserted by Theorem~\ref{PropSemiModel:Result}.
The case of props with non-empty outputs can be addressed similarly.

Recall that $\phi_!: \C^{\ACat}\rightleftarrows\C^{\BCat} :\phi^*$
forms a Quillen adjunction
and so does its restriction to objects with non-empty inputs.
Thus the functor $\phi_!$ maps (acyclic) cofibrations to (acyclic) cofibrations.
Since the definition of the free prop is more natural on $\Sigma_*$-biobjects,
we prove the following statement rather than the requirement of Theorem~\ref{SemiModelCategories:AdjointStructure}
in its original form:

\begin{mainlemm}\label{MainResult:PropPushouts}
For any pushout of props
\begin{equation*}
\mathrm{(*)}
\qquad\qquad\vcenter{\xymatrix{ \Free(K)\ar[r]^{u}\ar[d]_{\Free(i)} & \POp\ar@{.>}[d]^{f} \\
\Free(L)\ar@{.>}[r]_(0.33){v} & \Free(L)\bigvee_{\Free(K)}\POp }}\qquad\qquad
\end{equation*}
such that $\POp$ is an $\Free(\C^{\BCat})_c$-cell complex in $\Prop_0$,
the morphism $\phi^* U(f)$
forms a cofibration (respectively an acyclic cofibration) in $\C^{\ACat}_0$
whenever $i$ is a cofibration (respectively an acyclic cofibration) in $\C^{\BCat}_0$.
\end{mainlemm}

The proof of this lemma requires a good description of pushouts $\Free(L)\bigvee_{\Free(K)}\POp$.
The idea is to rearrange the construction of the pushout
in order to decompose the morphism
\begin{equation*}
\POp\xrightarrow{f}\Free(L)\bigvee_{\Free(K)}\POp
\end{equation*}
into a sequence of manageable pushouts in the category of $\Sigma_*$-biobjects.
This desired decomposition cannot be performed within the category of props
and we have to introduce new objects
in order to handle the structure of~$\Free(L)\bigvee_{\Free(K)}\POp$.

Pushouts of operads have a canonical representation in the underlying category.
This is not the case for pushouts in the category of props (see~\S\ref{PropPushoutHomotopy:CanonicalCoproductForm}).
For that reason,
the proof of Lemma~\ref{MainResult:PropPushouts} in the context of props
differs significantly from the case of operads
addressed in~\cite{HinichHomotopy,Spitzweck}.
Here is our plan:
\begin{enumerate}
\item\label{PushoutAnalysisPlan:Coproducts}
We observe that a coproduct of the form $\POp\vee\Free(M)$ is in some sense a functor of symmetric tensors in~$M$
with coefficients in a symmetric sequence $\Env_{\Prop}(\POp)$,
naturally associated to $\POp$,
and with $\POp$ as a leading term.
We study the morphism of symmetric sequences $\Env_{\Prop}(f): \Env_{\Prop}(\POp)\rightarrow\Env_{\Prop}(\Free(L)\bigvee_{\Free(K)}\POp)$
induced by $f: \POp\rightarrow\Free(L)\bigvee_{\Free(K)}\POp$
rather than $f$ itself,
because the object $\Env_{\Prop}(\Free(L)\bigvee_{\Free(K)}\POp)$
can be obtained from $\Env_{\Prop}(\POp)$ by an easy formula
which reflects coproduct identities at the functor level.
We have just to take the leading term of~$\Env_{\Prop}(f)$
to retrieve~$f$.
\item\label{PushoutAnalysisPlan:Decomposition}
We prove the existence of a decomposition
\begin{equation*}
\Env_{\Prop}(\Free(L)\bigvee_{\Free(K)}\POp) = \colim_r\Bigl\{\Env_{\Prop}^r(\Free(L)\bigvee_{\Free(K)}\POp)\Bigr\}
\end{equation*}
so that each term $\Env_{\Prop}^r(\Free(L)\bigvee_{\Free(K)}\POp)$
is built from the previous one by a pushout involving
an $r$-fold pushout-product
\begin{equation*}
\lambda: \colim\{L\otimes\dots\otimes K\otimes\dots\otimes L\}\rightarrow L^{\otimes r}.
\end{equation*}
\item\label{PushoutAnalysisPlan:Patching}
We apply usual patching techniques to prove that
the morphism of symmetric sequences
\begin{equation*}
\Env_{\Prop}(f): \Env_{\Prop}(\POp)\rightarrow\Env_{\Prop}(\Free(L)\bigvee_{\Free(K)}\POp)
\end{equation*}
consists of cofibrations (respectively, acyclic cofibrations) of $\Sigma_*$-biobjects.
We use that $f$ forms the leading term of the morphism of symmetric sequences $\Env_{\Prop}(f)$
to conclude that $f$ forms itself a cofibration (respectively, an acyclic cofibration) of $\Sigma_*$-biobjects.
\end{enumerate}

In the context of algebras over an operad,
a version of these results is proved in~\cite{BergerMoerdijk,BergerMoerdijkRepresentations} (assuming the existence of model structures)
and in~\cite[\S\S 19-20]{FresseModuleBook} (using a direct analysis of pushouts).
In both cases,
the obtained results apply to algebras over operads satisfying technical cofibration requirements.
In a sense,
props are algebras over a certain colored operad
and we adapt the arguments of the second-mentioned reference~\cite[\S\S 19-20]{FresseModuleBook}
to this operad.

The symmetric sequences,
which define the coefficients of functors of symmetric tensors on the category of $\Sigma_*$-biobjects,
are called symmetric $\Sigma_*$-multiobjects
and are introduced in~\S\ref{PropFunctor}.
The symmetric $\Sigma_*$-multiobject $\Env_{\Prop}(\POp)$ which represents the functor $\POp\vee\Free(-)$
is defined in that section.

Sections~\ref{PropPushoutDecomposition} and~\ref{PropPushoutHomotopy}
are devoted to parts~(\ref{PushoutAnalysisPlan:Decomposition}) and~(\ref{PushoutAnalysisPlan:Patching})
of the proof of Lemma~\ref{MainResult:PropPushouts}.

First of all,
we review the explicit description of free props
arising from a graphical interpretation of the operations of a prop.

\section{The language of graphs and free props}\label{Graphs}
In the point-set context,
the homomorphisms $p\in\POp(m,n)$ of a prop~$\POp$
are usually represented by a labeled box with $m$ inputs arranged on a horizontal upper line,
and $n$ outputs arranged on a horizontal lower line:
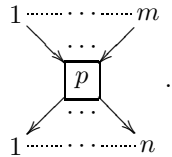
\begin{equation*}
\vcenter{\xymatrix@H=6pt@W=3pt@M=2pt@!R=1pt@!C=1pt{ 1\ar[dr]\ar@{.}[r] & \cdots\ar@{.}[r]\ar@{}[d]|{\displaystyle{\cdots}} &
m\ar[dl] \\
& *+<8pt>[F]{p}\ar[dl]\ar[dr]\ar@{}[d]|{\displaystyle{\cdots}} & \\
1\ar@{.}[r] & \cdots\ar@{.}[r] & n }}.
\end{equation*}
This graphical interpretation reflects the interpretation of~$p\in\POp(m,n)$
as an operation $p: 1^{\otimes m}\rightarrow 1^{\otimes n}$.

Tensor products of homomorphisms $p\in\POp(k,m)$ and $q\in\POp(l,n)$
are represented by horizontal concatenations
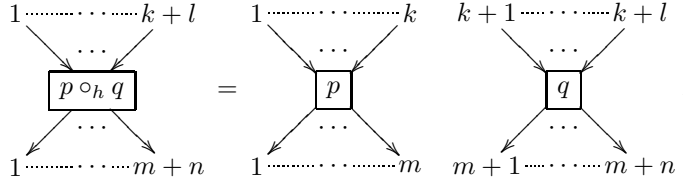
\begin{equation*}
\vcenter{\xymatrix@!R=0.5em@!C=0.5em@M=2pt{ 1\ar[dr]\ar@{.}[r] & \cdots\ar@{.}[r]\ar@{}[d]|{\displaystyle{\cdots}} & k+l\ar[dl] \\
& *+<8pt>[F]{p\circ_h q}\ar[dl]\ar[dr]\ar@{}[d]|{\displaystyle{\cdots}} & \\
1\ar@{.}[r] & \cdots\ar@{.}[r] & m+n }}
= \vcenter{\xymatrix@!R=0.5em@!C=0.5em@M=2pt{ 1\ar[dr]\ar@{.}[r] & \cdots\ar@{.}[r]\ar@{}[d]|{\displaystyle{\cdots}} & k\ar[dl] &
k+1\ar[dr]\ar@{.}[r] & \cdots\ar@{.}[r]\ar@{}[d]|{\displaystyle{\cdots}} & k+l\ar[dl] \\
& *+<8pt>[F]{p}\ar[dl]\ar[dr]\ar@{}[d]|{\displaystyle{\cdots}} &&& *+<8pt>[F]{q}\ar[dl]\ar[dr]\ar@{}[d]|{\displaystyle{\cdots}} & \\
1\ar@{.}[r] & \cdots\ar@{.}[r] & m & m+1\ar@{.}[r] & \cdots\ar@{.}[r] & m+n }}.
\end{equation*}
Composites of homomorphisms $q\in\POp(k,n)$ and $p\in\POp(m,k)$
are represented by vertical concatenations
so that outputs of~$p$ are plugged into inputs of~$q$:
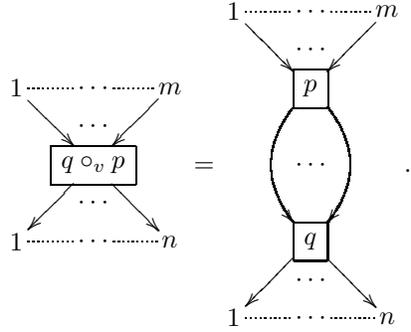
\begin{equation*}
\vcenter{\xymatrix@!R=0.5em@!C=0.5em@M=2pt{ 1\ar[dr]\ar@{.}[r] & \cdots\ar@{.}[r]\ar@{}[d]|{\displaystyle{\cdots}} & m\ar[dl] \\
& *+<8pt>[F]{q\circ_v p}\ar[dl]\ar[dr]\ar@{}[d]|{\displaystyle{\cdots}} & \\
1\ar@{.}[r] & \cdots\ar@{.}[r] & n }}
= \vcenter{\xymatrix@!R=0.5em@!C=0.5em@M=2pt{ 1\ar[dr]\ar@{.}[r] & \cdots\ar@{.}[r]\ar@{}[d]|{\displaystyle{\cdots}} & m\ar[dl] \\
& *+<8pt>[F]{p}\ar@/_1.5em/[dd]\ar@{}[dd]|{\displaystyle{\cdots}}\ar@/^1.5em/[dd] & \\
&&\\
& *+<8pt>[F]{q}\ar[dl]\ar[dr]\ar@{}[d]|{\displaystyle{\cdots}} & \\
1\ar@{.}[r] & \cdots\ar@{.}[r] & n }}.
\end{equation*}

The purpose of this section is to review an explicit description of free props
which relies on the graphical representation
of this introduction.
Intuitively,
the elements of a free prop $\Free(M)$ consist of formal composites of generating operations~$x\in M(m,n)$
which are modelled by directed graphs with these operations on vertices.
To begin with,
we define the graph structure involved in this construction of~$\Free(M)$.

\subsubsection{The language of graphs}\label{Graphs:Language}
We adapt the formalism of~\cite[\S I.2]{Serre}.
For us,
a directed graph $\Gamma$ with $m$ inputs and $n$ outputs (for short, an $(m,n)$-graph)
consists of a (finite) set of vertices $V(\Gamma)$
together with a (finite) set of edges $E(\Gamma)$,
oriented from a source $s(e)\in V(\Gamma)\amalg\{1,\dots,m\}$ to a target $t(e)\in V(\Gamma)\amalg\{1,\dots,n\}$,
such that:
\begin{enumerate}
\item
for each input $i\in\{1,\dots,m\}$,
there is one and only one edge $e\in E(\Gamma)$ such that $s(e) = i$,
\item
for each output $j\in\{1,\dots,n\}$,
there is one and only one edge $e\in E(\Gamma)$ such that $t(e) = j$,
\item
multiple edges
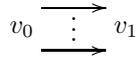
\begin{equation*}
\xymatrix{ v_0\ar@<+8pt>[r]\ar@<+3pt>@{}[r]|{\vdots}\ar@<-8pt>[r] & v_1 }
\end{equation*}
and parallel chains of edges
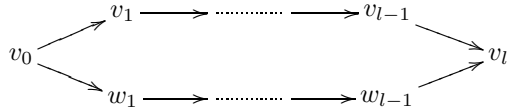
\begin{equation*}
\xymatrix@R=6pt@M=2pt{ & v_1\ar[r] & \ar@{.}[r] & \ar[r] & v_{l-1}\ar[dr] & \\
v_0\ar[ur]\ar[dr] &&&&& v_l \\
& w_1\ar[r] & \ar@{.}[r] & \ar[r] & w_{l-1}\ar[ur] & }
\end{equation*}
are allowed in $\Gamma$,
but we have no chains of edges
\begin{equation*}
\xymatrix{ v_0\ar[r] & \ar@{.}[r] & \ar[r] & v_l }
\end{equation*}
such that $v_0 = v_l$.
\end{enumerate}
For our needs,
we do not assume that a graph is necessarily connected.

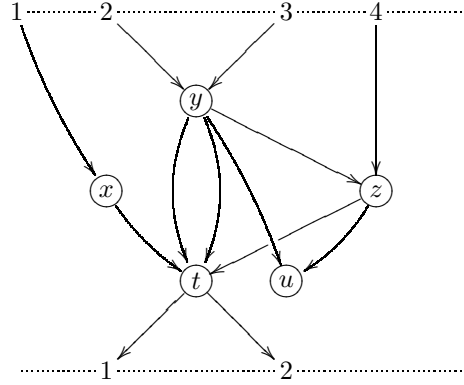
\begin{figure}[h]
\begin{equation*}
\xymatrix@!R=1em@!C=1em@M=2pt{ 1\ar@/_3pt/[ddr]\ar@{.}[r] & 2\ar[dr]\ar@{.}[rr] && 3\ar[dl]\ar@{.}[r] & 4\ar[dd]\ar@{.}[r] & \\
&& *+<8pt>[o][F]{y}\ar[drr]\ar@/^3pt/[ddr]\ar@/_9pt/[dd]\ar@/^9pt/[dd] &&& \\
& *+<8pt>[o][F]{x}\ar@/_3pt/[dr] &&& *+<8pt>[o][F]{z}\ar[dll]|(0.57){\hole}\ar@/^4pt/[dl] & \\
&& *+<8pt>[o][F]{t}\ar[dl]\ar[dr] & *+<8pt>[o][F]{u} & \\
\ar@{.}[r] & 1\ar@{.}[rr] && 2\ar@{.}[rr] && }
\end{equation*}
\caption{A directed graph with $4$ inputs and $2$ outputs}\label{Fig:DirectedGraph}
\end{figure}

An example of a graph with $4$ inputs, $2$ outputs and $5$ vertices $V(\Gamma) = \{x,y,z,t,u\}$
is displayed in Figure~\ref{Fig:DirectedGraph}.
To clarify the picture, we usually put inputs and outputs on horizontal lines.

The input set $\In_v$ of a vertex $v$
consists of the edges $e$ such that $t(e) = v$.
The output set $\Out_v$
consists of the edges $e$ such that $s(e) = v$.

Define an isomorphism of $(m,n)$-graphs $f: \Gamma\rightarrow\Delta$
to be a pair of bijections $f_V: V(\Gamma)\rightarrow V(\Delta)$
and $f_E: E(\Gamma)\rightarrow E(\Delta)$
that satisfy the obvious commutation relation with respect to the source and target of edges.
Let $\Graph(m,n)$ be the set of $(m,n)$-graphs together with this groupoid structure.

The collection of groupoids $\Graph(m,n)$, $m,n\in\NN$,
is equipped with a natural prop structure:
\begin{itemize}
\item
the symmetric groups $\Sigma_m$ and $\Sigma_n$ operate on $\Graph(m,n)$
by re-indexing inputs and outputs of $(m,n)$-graphs;
\item
the horizontal composites $\Gamma\circ_h\Delta$ are given by the disjoint union of graphs
together with an appropriate index shift on the inputs and outputs of $\Delta$
(see Figure~\ref{Fig:HorizontalComposite});
\item
the vertical composites $\Gamma\circ_v\Delta$ are defined by plugging the outputs of $\Delta$
into the corresponding inputs of $\Gamma$
(see Figure~\ref{Fig:VerticalComposite});
\item
the identity $\id\in\Graph(n,n)$, $n\in\NN$,
is represented by the parallel $(n,n)$-graph (see Figure~\ref{Fig:UnitGraph}).
\end{itemize}
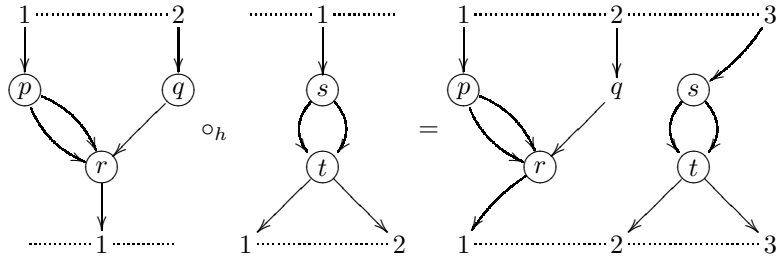
\begin{figure}[t]
\begin{equation*}
\vcenter{\xymatrix@!R=0.5em@!C=0.5em@M=2pt{ 1\ar[d]\ar@{.}[rr] && 2\ar[d] \\
*+<8pt>[o][F]{p}\ar@/_6pt/[dr]\ar@/^6pt/[dr] && *+<8pt>[o][F]{q}\ar[dl] \\
& *+<8pt>[o][F]{r}\ar[d] & \\
\ar@{.}[r] & 1\ar@{.}[r] & }}
\circ_h\vcenter{\xymatrix@!R=0.5em@!C=0.5em@M=2pt{ \ar@{.}[r] & 1\ar[d]\ar@{.}[r] & \\
& *+<8pt>[o][F]{s}\ar@/_9pt/[d]\ar@/^9pt/[d] & \\
& *+<8pt>[o][F]{t}\ar[dl]\ar[dr] & \\
1\ar@{.}[rr] && 2 }}
= \vcenter{\xymatrix@!R=0.5em@!C=0.5em@M=2pt{ 1\ar[d]\ar@{.}[rr] && 2\ar[d]\ar@{.}[rr] && 3\ar@/^3pt/[dl] \\
*+<8pt>[o][F]{p}\ar@/_6pt/[dr]\ar@/^6pt/[dr] && q\ar[dl] & *+<8pt>[o][F]{s}\ar@/_9pt/[d]\ar@/^9pt/[d] & \\
& *+<8pt>[o][F]{r}\ar@/_3pt/[dl] && *+<8pt>[o][F]{t}\ar[dl]\ar[dr] & \\
1\ar@{.}[rr] && 2\ar@{.}[rr] && 3 }}
\end{equation*}
\caption{A horizontal composite of graphs}\label{Fig:HorizontalComposite}
\end{figure}
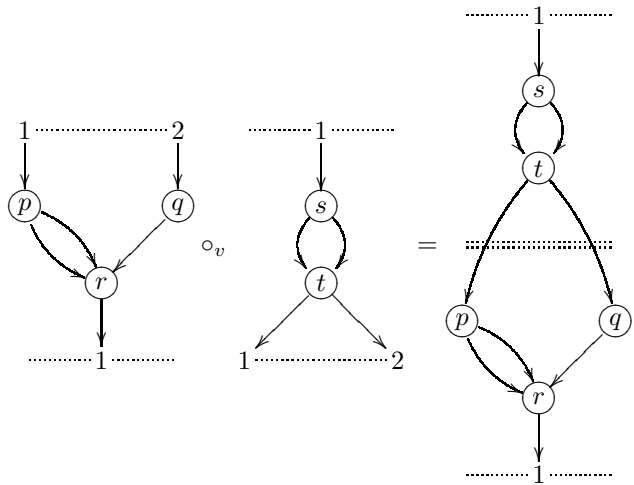
\begin{figure}[t]
\begin{equation*}
\vcenter{\xymatrix@!R=0.5em@!C=0.5em@M=2pt{ 1\ar[d]\ar@{.}[rr] && 2\ar[d] \\
*+<8pt>[o][F]{p}\ar@/_6pt/[dr]\ar@/^6pt/[dr] && *+<8pt>[o][F]{q}\ar[dl] \\
& *+<8pt>[o][F]{r}\ar[d] & \\
\ar@{.}[r] & 1\ar@{.}[r] & }}
\circ_v\vcenter{\xymatrix@!R=0.5em@!C=0.5em@M=2pt{ \ar@{.}[r] & 1\ar[d]\ar@{.}[r] & \\
& *+<8pt>[o][F]{s}\ar@/_9pt/[d]\ar@/^9pt/[d] & \\
& *+<8pt>[o][F]{t}\ar[dl]\ar[dr] & \\
1\ar@{.}[rr] && 2 }}
= \vcenter{\xymatrix@!R=0.5em@!C=0.5em@M=2pt{ \ar@{.}[r] & 1\ar[d]\ar@{.}[r] & \\
& *+<8pt>[o][F]{s}\ar@/_9pt/[d]\ar@/^9pt/[d] & \\
& *+<8pt>[o][F]{t}\ar@/_0.5em/[ddl]\ar@/^0.5em/[ddr] & \\
\ar@{:}[rr] && \\
*+<8pt>[o][F]{p}\ar@/_6pt/[dr]\ar@/^6pt/[dr] && *+<8pt>[o][F]{q}\ar[dl] \\
& *+<8pt>[o][F]{r}\ar[d] & \\
\ar@{.}[r] & 1\ar@{.}[r] & }}
\end{equation*}
\caption{A vertical composite of graphs}\label{Fig:VerticalComposite}
\end{figure}
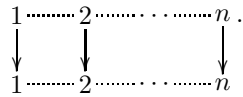
\begin{figure}[t]
\begin{equation*}
\xymatrix@!R=4pt@!C=2pt@M=2pt{ 1\ar[d]\ar@{.}[r] & 2\ar[d]\ar@{.}[r] & \cdots\ar@{.}[r] & n\ar[d] \\
1\ar@{.}[r] & 2\ar@{.}[r] & \cdots\ar@{.}[r] & n }.
\end{equation*}
\caption{The identity $(n,n)$-graph}\label{Fig:UnitGraph}
\end{figure}

\subsubsection{The explicit construction of free props}\label{Graphs:FreePropsConstruction}
Let $M$ be a $\Sigma_*$-biobject.

Let $I$ and $J$ be finite sets with $m$ and $n$ elements respectively.
For our needs,
we define objects $M(I,J)$
associated to $M$.
The idea is to reindex the inputs (respectively, outputs) of~$M(m,n)$
by $I$ (respectively, $J$).
Let $\Bij(K,L)$.
denote the set of bijections between any pair of finite sets.
Formally,
the object $M(I,J)$
is defined by the tensor product
\begin{equation*}
M(I,J) = \Bij(I,\{1,\dots,m\})\otimes_{\Sigma_m} M(m,n)\otimes_{\Sigma_n}\Bij(\{1,\dots,n\},J)
\end{equation*}
in which internal actions of permutations on $M(m,n)$
are made equal to translations of bijections.
Hence,
an element of~$M(I,J)$
correctly represents an element of~$M(m,n)$
whose inputs (respectively, outputs) are in bijection with the elements of~$I$ (respectively,~$J$).

For a graph $\Gamma$,
we form the tensor product $\Free_{\Gamma}(M) = \bigotimes_{v\in V(\Gamma)} M(\In_v,\Out_v)$.
In the point set context,
an element of $\Free_{\Gamma}(M)$
represents a labeling of the vertices $v\in V(\Gamma)$ by elements $x_v\in M(m_v,n_v)$
together with a bijection between the inputs (outputs) of $v$ and the inputs (outputs) of~$x_v$.
For instance,
the picture of Figure~\ref{Fig:DirectedGraph} can be used to represent a graph whose vertices are labeled by elements
$x\in M(1,1)$, $y\in M(2,4)$, $z\in M(2,2)$, $t\in M(4,2)$ and $u\in M(2,0)$ respectively.

The map $\Gamma\mapsto\Free_{\Gamma}(M)$
defines a functor on the groupoid of $(m,n)$-graphs.
We define the underlying $\Sigma_*$-biobject of the free prop
by the colimits:
\begin{equation*}
\Free(M)(m,n) = \colim_{\Gamma\in\Graph(m,n)}\Free_{\Gamma}(M).
\end{equation*}
The right action of a permutation $w\in\Sigma_m$ on $\Free(M)$
is yielded by the re-indexing functor $w^*: \Graph(m,n)\rightarrow\Graph(m,n)$
associated to $w$
and symmetrically as regards the left action of permutations $w\in\Sigma_n$.

For graphs $\Gamma,\Delta$,
the obvious bijections $V(\Gamma\circ_h\Delta) = V(\Gamma)\amalg V(\Delta)$
yield natural isomorphisms $\Free_{\Gamma\circ_h\Delta}(M)\simeq\Free_{\Gamma}(M)\otimes\Free_{\Delta}(M)$
which give rise to the horizontal composition operation of the free prop $\Free(M)$.
We also have a canonical bijection $V(\Gamma\circ_v\Delta) = V(\Gamma)\amalg V(\Delta)$
for vertical composites of graphs.
These bijections give natural isomorphisms $\Free_{\Gamma\circ_v\Delta}(M)\simeq\Free_{\Gamma}(M)\otimes\Free_{\Delta}(M)$
which yield the vertical composition operation of the free prop $\Free(M)$.
The unit $\eta: \unit\rightarrow\Free(M)(n,n)$
identifies the unit object $\unit$
with the summand of $\Free(M)(n,n)$
associated to the identity element of $\Graph(n,n)$.

It is easy to check that the definition of this paragraph provides $\Free(M)$
with a well-defined prop structure.

In addition,
we have a natural morphism $\eta: M(m,n)\rightarrow\Free(M)(m,n)$
which identifies $M(m,n)$
with the summand of~$\Free(M)(m,n)$
associated to the $(m,n)$-corolla:
\begin{equation*}
X_{m n} = \vcenter{\xymatrix@H=6pt@W=3pt@M=2pt@!R=1pt@!C=1pt{ 1\ar[dr]\ar@{.}[r] & \cdots\ar@{.}[r]\ar@{}[d]|{\displaystyle{\cdots}} & m\ar[dl] \\
& *+<8pt>[o][F]{x}\ar[dl]\ar[dr]\ar@{}[d]|{\displaystyle{\cdots}} & \\
1\ar@{.}[r] & \cdots\ar@{.}[r] & n }}.
\end{equation*}
One proves:

\begin{prop}[{See~\cite{EnriquezEtingof,ValletteThesis}}]\label{Graphs:FreePropDefinition}
The prop $\Free(M)$ together with the morphism $\eta: M\rightarrow\Free(M)$ defined in~\S\ref{Graphs:FreePropsConstruction}
satisfies the universal property of a free prop:
any morphism of $\Sigma_*$-biobjects $\phi: M\rightarrow\QOp$ towards a prop $\QOp$
has a unique factorization
\begin{equation*}
\xymatrix{ M\ar[rr]^{\phi}\ar[dr]_{\eta} && \QOp \\ & \Free(M)\ar@{.>}[ur]_{\exists!\tilde{\phi}} & }
\end{equation*}
such that $\tilde{\phi}: \Free(M)\rightarrow\QOp$ is a prop morphism.
\qed
\end{prop}



\subsubsection{Monad structures}\label{Graphs:PropMonad}
The identity morphism of a prop $\id: \POp\rightarrow\POp$ determines a prop morphism $\lambda: \Free(\POp)\rightarrow\POp$,
to which we refer as the total composition product of~$\POp$.
This morphism defines the augmentation of the adjunction between $\Sigma_*$-biobjects and props.

This construction applied to the free prop $\POp = \Free(M)$
gives a universal composition product $\mu: \Free(\Free(M))\rightarrow\Free(M)$
that makes $\Free: \C^{\BCat}\rightarrow\C^{\BCat}$
a monad on the category of $\Sigma_*$-biobjects (see~\cite{MacLaneCategories}).
The category of props is clearly isomorphic to the category of algebras over this monad.

\begin{figure}[t]
\begin{equation*}
\vcenter{\xymatrix@!R=0.6em@!C=0.6em@M=2pt{ \ar@{.}[r] & 1\ar@{.}[r]\ar[d] & 2\ar[d]\ar@{.}[r] & \\
& *+<8pt>[o][F]{x}\ar[d]\ar[dr]\save[]!C.[r]!C *+<14pt>[F:2pt]\frm{}\restore & *+<8pt>[o][F]{y}\ar[d]\ar[dl]|{\hole} & \\
& *+<8pt>[o][F]{y}\ar[d]\ar@{}[dr]\save[]!C.[r]!C *+<14pt>[F:2pt]\frm{}\restore & *+<8pt>[o][F]{t}\ar[d] & \\
\ar@{.}[r] & 1\ar@{.}[r] & 2\ar@{.}[r] & }}
\quad\mapsto\quad\vcenter{\xymatrix@!R=0.6em@!C=0.6em@M=2pt{ \ar@{.}[r] & 1\ar@{.}[r]\ar[d] & 2\ar[d]\ar@{.}[r] & \\
& *+<8pt>[o][F]{x}\ar[d]\ar[dr]\save[]!C.[r]!C *+<14pt>[F.:2pt]\frm{}\restore & *+<8pt>[o][F]{y}\ar[d]\ar[dl]|{\hole} & \\
& *+<8pt>[o][F]{y}\ar[d]\ar@{}[dr]\save[]!C.[r]!C *+<14pt>[F.:2pt]\frm{}\restore & *+<8pt>[o][F]{t}\ar[d] & \\
\ar@{.}[r] & 1\ar@{.}[r] & 2\ar@{.}[r] & }}
\quad\mapsto\quad\vcenter{\xymatrix@!R=0.6em@!C=0.6em@M=2pt{ \ar@{.}[r] & 1\ar@{.}[r]\ar[d] & 2\ar[d]\ar@{.}[r] & \\
& *+<8pt>[o][F]{x}\ar[d]\ar[dr] & *+<8pt>[o][F]{y}\ar[d]\ar[dl]|{\hole} & \\
& *+<8pt>[o][F]{y}\ar[d]\ar@{}[dr] & *+<8pt>[o][F]{t}\ar[d] & \\
\ar@{.}[r] & 1\ar@{.}[r] & 2\ar@{.}[r] & }}
\end{equation*}
\caption{The universal composition product}\label{Fig:UniversalComposition}
\end{figure}
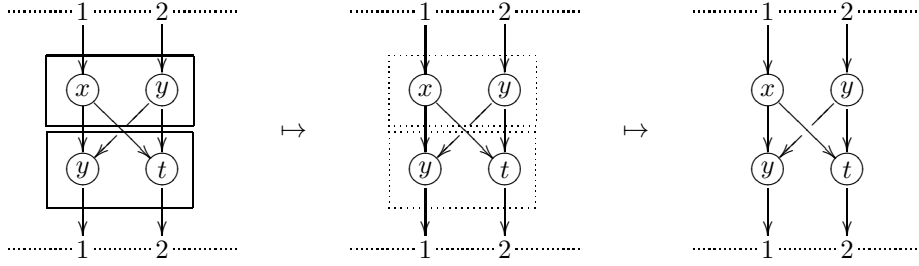

Graphically,
an element of~$\Free(\Free(M))$ is a composite $(m,n)$-graph
whose vertices are themselves $(m_v,n_v)$-graphs
together with an $M$-labeling.
The universal composition product $\mu: \Free(\Free(M))\rightarrow\Free(M)$
simply expands the inner graph structure of vertices and forget the boundary of these inner graphs
to return an $M$-labeled $(m,n)$-graph without extra structure.
This process is represented in Figure~\ref{Fig:UniversalComposition}.

\section{Symmetric $\Sigma_*$-multiobjects and functors}\label{PropFunctor}
According to the construction of the previous section,
the free prop $F(M)$ is defined by colimits of tensor products $M(a_*,b_*)^{\otimes r} = \bigotimes_{i=1}^{r} M(a_i,b_i)$.
We introduce symmetric $\Sigma_*$-multiobjects as a structure to model functors of this form
on the category of $\Sigma_*$-biobjects.
Precise definitions are given in the next paragraph.

We prove that the coproduct~$\POp\vee\Free(M)$ in the category of props~$\Prop$
is a functor associated to a certain symmetric $\Sigma_*$-multiobject,
for which we adopt the notation~$\Env_{\Prop}(\POp)$.
In the next section,
we study the symmetric $\Sigma_*$-multiobject $\Env_{\Prop}(\Free(L)\bigvee_{\Free(K)}\POp)$
associated to a pushout along a morphism of free props~$\Free(i): \Free(K)\rightarrow\Free(L)$.
For the moment,
we only observe that the symmetric $\Sigma_*$-multiobject $\Env_{\Prop}(\POp\vee\Free(M))$
associated to a coproduct with a free prop is determined from $\Env_{\Prop}(\POp)$
by a simple formula which reflects the identity~$\POp\vee\Free(M)\vee\Free(N) = \POp\vee\Free(M\oplus N)$
at the functor level.

\subsubsection{Symmetric $\Sigma_*$-multiobjects and associated functors}\label{PropFunctor:Definitions}
Recall that $\BCat$ refers to the category formed by pairs of integers $(m,n)\in\NN\times\NN$
and morphism sets such that
\begin{equation*}
\Mor_{\BCat}((m,n),(p,q)) = \begin{cases} \Sigma_m^{op}\times\Sigma_n, & \text{if $m = p$ and $n = q$}, \\
\emptyset, & \text{otherwise}. \end{cases}
\end{equation*}
Let $\Sigma_*\wr\BCat\times\BCat$
be the category formed by collections of objects
\begin{equation*}
((a_1,b_1),\dots,(a_r,b_r);(m,n))\in\BCat^{\times r}\times\BCat,\quad\text{$r\in\NN$},
\end{equation*}
whose morphisms are formal composites
\begin{multline*}
((a_1,b_1),\dots,(a_r,b_r);(m,n))
\xrightarrow{(\rho_1,\dots,\rho_r;\sigma)}
((c_1,d_1),\dots,(c_r,d_r);(p,q))\\
\xrightarrow{w^*} ((c_{w(1)},d_{w(1)}),\dots,(c_{w(r)},d_{w(r)});(p,q)),
\end{multline*}
where $w\in\Sigma_r$ and $(\rho_1,\dots,\rho_r;\sigma)\in\Mor_{\BCat}((a_*,b_*),(c_*,d_*))^{\times r}\times\Mor_{\BCat}((m,n),(p,q))$.
The composition of morphisms in $\Sigma_*\wr\BCat\times\BCat$
is given by the usual formula of wreath products.
We define a symmetric $\Sigma_*$-multiobject~$S$
as a diagram $S\in\C^{\Sigma_*\wr\BCat\times\BCat}$.

The functor on $\Sigma_*$-biobjects associated to a symmetric $\Sigma_*$-multiobject~$S$
is defined by symmetric tensor products
\begin{equation*}
S(M)(m,n) = \bigoplus_{r\in\NN}\Bigl\{\bigoplus_{(a_*,b_*)} S((a_*,b_*);(m,n))
\otimes_{\Sigma_{(a_*,b_*)}} M(a_*,b_*)^{\otimes r}\Bigr\}/\Sigma_r,
\end{equation*}
where the sum ranges over collections $(a_*,b_*) = (a_1,b_1),\dots,(a_r,b_r)$
and we use the short notation
\begin{equation*}
M(a_*,b_*)^{\otimes r} = \bigotimes_{i=1}^{r} M(a_i,b_i)
\quad\text{and}\quad\Sigma_{(a_*,b_*)} = \prod_{i=1}^{r} (\Sigma_{a_i}\times\Sigma_{b_i}).
\end{equation*}
The quotient under the action of~$\Sigma_r$
makes permutations of tensor products $M(a_*,b_*)^{\otimes r}$
agree with the internal $\Sigma_r$-action of~$S$.
Throughout the paper,
we use the notation of the symmetric $\Sigma_*$-multiobject~$S$
to represent the associated functor~$S: \C^{\BCat}\rightarrow\C^{\BCat}$.

The functor $S: \C^{\BCat}\rightarrow\C^{\BCat}$ associated to a symmetric $\Sigma_*$-multiobject $S$
preserves reflexive coequalizers and filtered colimits,
like any composite of tensor products (see~\cite[\S 1.2]{FresseModuleBook}),
and we have the pointwise identity $(\colim_i S_i)(M) = \colim_i S_i(M)$,
for any colimit of symmetric $\Sigma_*$-multiobjects $S_i$.

\subsubsection{Shifted symmetric $\Sigma_*$-multiobjects}\label{PropFunctor:ShiftedObjects}
In our arguments,
we use a shifted symmetric $\Sigma_*$-multiobject $S[M]$, associated to any $M\in\C^{\BCat}$,
formed by the partial evaluations of~$S$ on~$M$
\begin{multline*}
S[M]((a_1,b_1),\dots,(a_r,b_r);(m,n))\\
= \bigoplus_{s\in\NN}\Bigl\{\bigoplus_{(c_*,d_*)} S((a_*,b_*),(c_*,d_*);(m,n))
\otimes_{\Sigma_{(c_*,d_*)}} M(c_*,d_*)^{\otimes s}\Bigr\}/\Sigma_s,
\end{multline*}
where the sum ranges over collections $(c_*,d_*) = (c_1,d_1),\dots,(c_s,d_s)$.
We clearly have $S(M)(m,n) = S[M](\emptyset;(m,n))$.

For a colimit of symmetric $\Sigma_*$-multiobjects $S_i$,
we have an obvious pointwise identity
\begin{equation*}
(\colim_i S_i)[M] = \colim_i S_i[M],
\end{equation*}
which holds for every $M\in\C^{\BCat}$.
In contrast,
the functor $S[-]: M\mapsto S[M]$ associated to a fixed symmetric $\Sigma_*$-multiobject $S$
preserves reflexive coequalizers and filtered colimits in~$M$,
but not every colimit (like the unshifted functor associated to $S$).
Nevertheless the partial evaluation of~$S$ on a sum $M\oplus N$
is determined by an easy formula,
namely:

\begin{obsv}\label{PropFunctor:SumFunctor}
We have a natural isomorphism $S[M\oplus N]\simeq S[M][N]$, for every $M,N\in\C^{\BCat}$.
\end{obsv}

In particular, if we take $(\emptyset,(m,n))$-terms of~$S[M\oplus N]\simeq S[M][N]$,
then we obtain the relation $S(M\oplus N)\simeq S[M](N)$.
This observation, which is immediate from the definition of~\S\ref{PropFunctor:ShiftedObjects},
motivates the introduction of shifted objects.

\subsubsection{The representation of free props}\label{PropFunctor:FreePropFunctor}
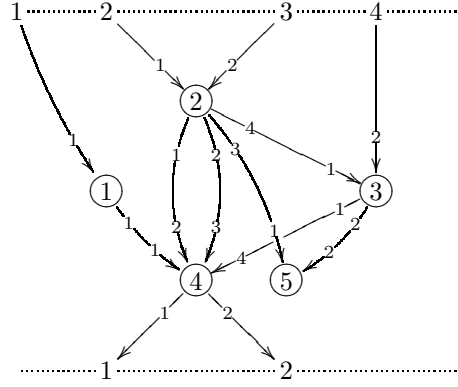
\begin{figure}[h]
\begin{equation*}
\xymatrix@!R=1em@!C=1em@M=2pt{ 1\ar@/_3pt/[ddr]|(0.7){1}\ar@{.}[r] & 2\ar[dr]|(0.6){1}\ar@{.}[rr] && 3\ar[dl]|(0.6){2}\ar@{.}[r] & 4\ar[dd]|(0.7){2}\ar@{.}[r] & \\
&& *+<8pt>[o][F]{2}\ar[drr]|(0.3){4}|(0.75){1}\ar@/^6pt/[ddr]|(0.3){3}|(0.73){\hole}|(0.75){1}\ar@/_9pt/[dd]|(0.28){\hole}|(0.3){1}|(0.7){2}\ar@/^9pt/[dd]|(0.3){2}|(0.7){3} &&& \\
& *+<8pt>[o][F]{1}\ar@/_3pt/[dr]|(0.3){1}|(0.6){1} &&& *+<8pt>[o][F]{3}\ar[dll]|(0.2){1}|(0.54){\hole}|(0.75){4}\ar@/^4pt/[dl]|(0.3){2}|(0.6){2} & \\
&& *+<8pt>[o][F]{4}\ar[dl]|(0.35){1}\ar[dr]|(0.35){2} & *+<8pt>[o][F]{5} & \\
\ar@{.}[r] & 1\ar@{.}[rr] && 2\ar@{.}[rr] && }
\end{equation*}
\caption{A directed graph equipped with a vertex numbering}\label{Fig:NumberedGraph}
\end{figure}

Let $\Graph((a_1,b_1),\dots,(a_r,b_r);(m,n))$
be the category formed by $(m,n)$-graphs $\Gamma$
with $r$ vertices $v_1,\dots,v_r$, numbered from $1$ to $r$,
such that $v_i$ has $a_i$ inputs, numbered from $1$ to $a_i$,
and $b_i$ outputs, numbered from $1$ to $b_i$.
Formally,
an object of $\Graph((a_1,b_1),\dots,(a_r,b_r);(m,n))$
consists of an $(m,n)$-graph $\Gamma$
together with a bijection $v_*: \{1,\dots,r\}\rightarrow V(\Gamma)$,
which defines the numbering of the vertices of~$\Gamma$,
and bijections $i_*: \{1,\dots,a_i\}\rightarrow\In_{v_i}$ and $j_*: \{1,\dots,b_i\}\rightarrow\Out_{v_i}$
which define the numbering of the inputs and outputs of each vertex~$v_i\in V(\Gamma)$.
Figure~\ref{Fig:NumberedGraph}
gives an example of such a graph structure.

The morphisms of $\Graph((a_1,b_1),\dots,(a_r,b_r);(m,n))$
are the isomorphisms of $(m,n)$-graphs
which preserve the numbering of vertices and the numbering of the inputs and outputs of vertices.
The objects $\Graph((a_1,b_1),\dots,(a_r,b_r);(m,n))$
form a symmetric $\Sigma_*$-multiobject in the category of groupoids,
the action of symmetric groups on this symmetric $\Sigma_*$-multiobject is simply defined by the usual renumbering process.

Let $F((a_1,b_1),\dots,(a_r,b_r);(m,n))$
be the object defined by the sum of copies of the unit object $\unit\in\C$
over the groupoid $\Graph((a_1,b_1),\dots,(a_r,b_r);(m,n))$
divided out by the relation
which identifies the summands associated to isomorphic objects.
The collection $F((a_1,b_1),\dots,(a_r,b_r);(m,n))$ inherits symmetric group actions
from the indexing category~$\Graph((a_1,b_1),\dots,(a_r,b_r);(m,n))$
and forms a symmetric $\Sigma_*$-multiobject.
The functor associated to this symmetric $\Sigma_*$-multiobject
\begin{equation*}
F(M)(m,n) = \bigoplus_{r\in\NN}\Bigl\{\bigoplus_{(a_*,b_*)} F((a_*,b_*);(m,n))
\otimes_{\Sigma_{(a_*,b_*)}} M(a_*,b_*)^{\otimes r}\Bigr\}/\Sigma_r
\end{equation*}
is formally identified with the underlying $\Sigma_*$-biobject
of the free prop $F(M)$,
because we have
\begin{multline*}
\Bigl\{F((a_*,b_*);(m,n))\otimes_{\Sigma_{(a_*,b_*)}} M(a_*,b_*)^{\otimes r}\Bigr\}/\Sigma_r\\
\simeq\colim_{\Graph((a_*,b_*);(m,n))}\unit\otimes M(a_*,b_*)^{\otimes r}/\equiv\\
\simeq\colim_{\Graph((a_*,b_*);(m,n))} M(a_*,b_*)^{\otimes r}/\equiv
\end{multline*}
by definition of $F((a_*,b_*);(m,n))$
and this colimit has the same summands and the same relations
as the construction of~\S\ref{Graphs:FreePropsConstruction}.

The following easy proposition
is a crucial ingredient of the verifications
of~\S\ref{PropPushoutHomotopy}:

\begin{prop}\label{PropFunctor:FreePropSymmetricStructure}
The subobject $\overline{F}\subset F$ generated by components
\begin{equation*}
F((a_1,b_1),\dots,(a_r,b_r);(m,n))
\quad\text{such that $m\not=0$ and $a_i\not=0$, for every $i = 1,\dots,r$,}
\end{equation*}
forms a free $\Sigma_r$-object
in the sense that $\Sigma_r$ permutes freely the summands of~$\overline{F}$.
A similar assertion holds for the subobject of $F$ generated by components
such that $n\not=0$ and $b_i\not=0$ for every $i = 1,\dots,r$.
\end{prop}

The condition amounts to considering graphs $\Gamma$
whose vertices $v$ have non-empty inputs $\In_v\not=\emptyset$
since $(a_1,\dots,a_r)$
represent the number of inputs of the vertices of~$\Gamma$.

\begin{proof}
We use the labelings of edge-paths from tree inputs to vertices to define a canonical ordering on the vertex set of each graph.
We consider the lexicographic ordering of edge-path labelings
and order vertices $v$ according to the minimum labeling of edge-paths
which abut to $v$.
In the example of Figure~\ref{Fig:NumberedGraph},
we obtain (in the lexicographic order):
\begin{align*}
& \text{minimal labeling toward vertex 1} = (\text{graph input}: 1 ;\ \text{edge labels}: 1), \\
& \text{minimal labeling toward vertex 4} = (\text{graph input}: 1 ;\ \text{edge labels}: 1, 1, 1), \\
& \text{minimal labeling toward vertex 2} = (\text{graph input}: 2 ;\ \text{edge labels}: 1), \\
& \text{minimal labeling toward vertex 5} = (\text{graph input}: 2 ;\ \text{edge labels}: 1, 3, 1), \\
& \text{minimal labeling toward vertex 3} = (\text{graph input}: 2 ;\ \text{edge labels}: 1, 4, 1).
\end{align*}

This ordering is canonical in the sense
that it is invariant under graph isomorphisms preserving input labels.
The ordering of vertices defines a canonical bijection $V(\Gamma)\simeq\{1,\dots,r\}$
and we have a one-to-one correspondence between vertex numberings $i_*: \{1,\dots,r\}\xrightarrow{\simeq} V(\Gamma)$
and permutations $\sigma: \{1,\dots,r\}\xrightarrow{\simeq}\{1,\dots,r\}$.
The conclusion follows.

To apply this argument we simply have to assume that each vertex of~$\Gamma$ has at least one input.
The case of graphs whose vertices have non-empty output is addressed symmetrically.
\end{proof}

\subsubsection{Partial evaluation of free props and composition products}\label{PropFunctor:PartialEvaluation}
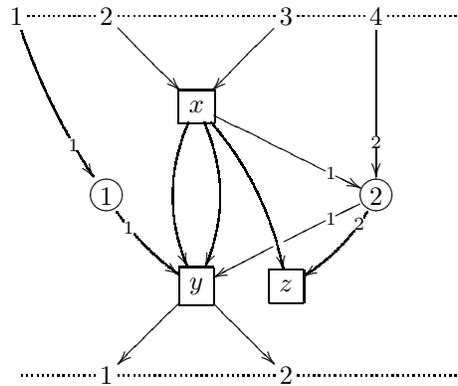
\begin{figure}[h]
\begin{equation*}
\xymatrix@!R=1em@!C=1em@M=2pt{ 1\ar@/_3pt/[ddr]|(0.69){\hole}|(0.7){1}\ar@{.}[r] & 2\ar[dr]\ar@{.}[rr] && 3\ar[dl]\ar@{.}[r] & 4\ar[dd]|(0.7){2}\ar@{.}[r] & \\
&& *+<8pt>[F]{x}\ar[drr]|(0.75){1}\ar@/^6pt/[ddr]\ar@/_9pt/[dd]\ar@/^9pt/[dd] &&& \\
& *+<8pt>[o][F]{1}\ar@/_3pt/[dr]|(0.3){1} &&& *+<8pt>[o][F]{2}\ar[dll]|(0.25){1}|(0.54){\hole}\ar@/^4pt/[dl]|(0.25){2} & \\
&& *+<8pt>[F]{y}\ar[dl]\ar[dr] & *+<8pt>[F]{z} & \\
\ar@{.}[r] & 1\ar@{.}[rr] && 2\ar@{.}[rr] && }
\end{equation*}
\caption{A directed graph equipped with a partial labeling of vertices}\label{Fig:PartiallyLabelledGraph}
\end{figure}

In the next paragraphs,
we use the shifted objects $F[M]$
defined by a partial evaluation of the symmetric $\Sigma_*$-multiobject of~\S\ref{PropFunctor:FreePropFunctor}.
Graphically,
an element of $F[M]$
is represented by an $(m,n)$-graph $\Gamma$
with a subset of vertices associated to an $M$-labeling
and a numbering of remaining vertices.
An example is represented in Figure~\ref{Fig:PartiallyLabelledGraph}.

The universal composition product $\mu: \Free(\Free(M))\rightarrow\Free(M)$
has a natural extension to shifted objects:
\begin{equation*}
\mu: F[\Free(M)]\rightarrow F[M].
\end{equation*}
This shifted universal composition product can be defined by the explicit process of~\S\ref{Graphs:PropMonad}.
To give an abstract interpretation of this morphism,
we use the identities $F[\Free(M)](N) = \Free(\Free(M)\oplus N)$
and $F[M](N) = \Free(M\oplus N)$
given by observation~\ref{PropFunctor:ShiftedObjects}.
Let $i: M\hookrightarrow M\oplus N$ and $j: N\hookrightarrow M\oplus N$
be the canonical morphisms associated to the sum $M\oplus N$.
The functor morphism $\mu(N): \Free(\Free(M)\oplus N)\rightarrow\Free(M\oplus N)$
induced by the shifted universal composition product
is identified with the morphism of props
which fits in the diagram:
\begin{equation*}
\xymatrix{ \Free(M)\oplus N\ar[rr]^{(\Free(i),\Free(j)\eta)}\ar[dr]_{\eta} && \Free(M\oplus N) \\
& \Free(\Free(M)\oplus N)\ar@{.>}[ur]_{\mu(N)} & }.
\end{equation*}

\medskip
Our next goal is to prove that the coproduct $\POp\vee\Free(M)$
is an instance of a functor associated to a symmetric $\Sigma_*$-multiobject.
To define this symmetric $\Sigma_*$-multiobject
we use the realization of colimits
by reflexive coequalizers in the category of $\Sigma_*$-biobjects.
The general construction of Proposition~\ref{PropSemiModel:Colimits}
can be simplified for certain particular colimits.
In the case of a coproduct with a free prop~$\POp\vee\Free(M)$,
we have:

\begin{prop}\label{PropFunctor:Coproduct}
The coproduct with a free prop
is realized by a reflexive coequalizer of the form:
\renewcommand{\theequation}{*}
\begin{equation}\label{FreeCoproductPresentation}
\xymatrix{ \Free(\Free(\POp)\oplus M)\ar@<+2pt>[r]^{d_0}\ar@<-2pt>[r]_{d_1} &
\Free(\POp\oplus M)\ar@{.>}[r]\ar@/_8mm/[l]_{s_0} & \POp\vee\Free(M) }.
\end{equation}
\renewcommand{\theequation}{\arabic{equation}}
\end{prop}

\begin{proof}
We simply recall the definition of the morphisms $(d_0,d_1,s_0)$.
The proof that $\coker(d_0,d_1)$ realizes the coproduct $\POp\vee\Free(M)$
reduces to formal verifications
and will be ommitted.

In the definition of~(\ref{FreeCoproductPresentation}),
we identify a $\Sigma_*$-biobject $N$ with a subobject of the free prop~$\Free(N)$
and morphisms towards $N$ with morphisms towards $\Free(N)$.

The face $d_0$ of our coequalizer is the morphism of free props induced
\begin{itemize}
\item
by the identity of~$M$
on the summand $M\subset\Free(\POp)\oplus M$,
\item
and by the composition product
\begin{equation*}
\Free(\POp)\xrightarrow{\lambda}\POp\hookrightarrow\POp\oplus M
\end{equation*}
on $\Free(\POp)\subset\Free(\POp)\oplus M$.
\end{itemize}
The face $d_1$ is induced
\begin{itemize}
\item
by the identity of~$M$ on $M\subset\Free(\POp)\oplus M$,
\item
and by the canonical morphism $\Free(\POp)\rightarrow\Free(\POp\oplus M)$
induced by $\POp\hookrightarrow\POp\oplus M$
on $\Free(\POp)\subset\Free(\POp)\oplus M$.
\end{itemize}
The degeneracy $s_0$ is induced
\begin{itemize}
\item
by the identity of~$M$
on $M\subset\POp\oplus M$
\item
and by the universal morphism $\eta: \POp\rightarrow\Free(\POp)$
on $\POp\subset\POp\oplus M$.
\end{itemize}
The relation $d_0 s_0 = d_1 s_0 = \id$
is immediate.
\end{proof}

The objects of the coequalizer~(\ref{FreeCoproductPresentation})
are functors of~$M$ associated to symmetric $\Sigma_*$-multiobjects
since we have
\begin{equation*}
F(\POp\oplus M)\simeq F[\POp](M)\quad\text{and}\quad F(\Free(\POp)\oplus M)\simeq F[\Free(\POp)](M)
\end{equation*}
by observation~\ref{PropFunctor:SumFunctor}.
We have moreover:

\begin{lemm}\label{PropFunctor:PropCoproductCoequalizerRepresentation}
The natural transformations of Proposition~\ref{PropFunctor:Coproduct}
are induced by morphisms of symmetric $\Sigma_*$-multiobjects
$d_0,d_1: \Free[\Free(\POp)]\rightrightarrows\Free[\POp]$
and $s_0: \Free[\POp]\rightarrow\Free[\Free(\POp)]$
such that $d_0 s_0 = d_1 s_0 = \id$.
\end{lemm}

\begin{proof}
Take:
\begin{itemize}
\item
the morphism $\Free[\lambda]: \Free[\Free(\POp)]\rightrightarrows\Free[\POp]$
induced by the composition product of~$\POp$
for $d_0$,
\item
the shifted universal composition product $\mu: \Free[\Free(-)]\rightarrow\Free[-]$
for $d_1$,
\item
the morphism $\Free[\eta]: \Free[\POp]\rightrightarrows\Free[\Free(\POp)]$
induced by the universal morphism of the free prop
for $s_0$.
\end{itemize}
The equality of the natural transformations of Proposition~\ref{PropFunctor:Coproduct}
with the natural transformations induced by these morphisms
is proved by a straightforward inspection.
\end{proof}

Let $\Env_{\Prop}(\POp)$
be the symmetric $\Sigma_*$-multiobject defined by the reflexive coequalizer of Lemma~\ref{PropFunctor:PropCoproductCoequalizerRepresentation}:
\renewcommand{\theequation}{**}
\begin{equation}\label{FreeCoproductRealization}
\xymatrix{ \Free[\Free(\POp)]\ar@<+2pt>[r]^{d_0}\ar@<-2pt>[r]_{d_1} & \Free[\POp]\ar@{.>}[r]\ar@/_8mm/[l]_{s_0} & \Env_{\Prop}(\POp) }.
\end{equation}
\renewcommand{\theequation}{\arabic{equation}}
In~\S\ref{PropFunctor:Definitions},
we observe that the functor associated to a colimit of symmetric $\Sigma_*$-multiobjects~$S_i$
satisfies the pointwise identity $(\colim_i S_i)(M)\simeq\colim_i S_i(M)$,
for every $M\in\C^{\BCat}$.
This observation
applied to $\Env_{\Prop}(\POp)$
gives immediately:

\begin{prop}\label{PropFunctor:PropCoproductRepresentation}
We have a natural isomorphism $\POp\vee\Free(M)\simeq\Env_{\Prop}(\POp)(M)$,
for every $M\in\C^{\BCat}$.\qed
\end{prop}

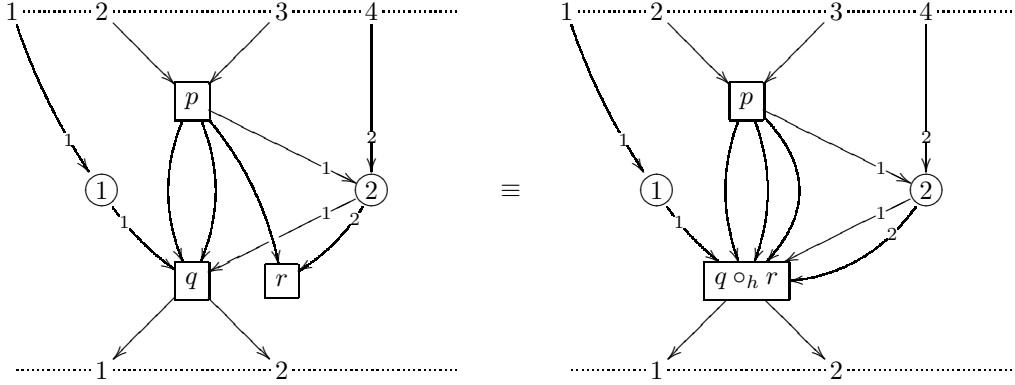
\begin{figure}[h]
\begin{equation*}
\vcenter{\xymatrix@!R=1em@!C=1em@M=2pt{ 1\ar@/_3pt/[ddr]|(0.69){\hole}|(0.7){1}\ar@{.}[r] & 2\ar[dr]\ar@{.}[rr] && 3\ar[dl]\ar@{.}[r] & 4\ar[dd]|(0.7){2}\ar@{.}[r] & \\
&& *+<8pt>[F]{p}\ar[drr]|(0.75){1}\ar@/^6pt/[ddr]\ar@/_9pt/[dd]\ar@/^9pt/[dd] &&& \\
& *+<8pt>[o][F]{1}\ar@/_3pt/[dr]|(0.3){1} &&& *+<8pt>[o][F]{2}\ar[dll]|(0.25){1}|(0.54){\hole}\ar@/^4pt/[dl]|(0.25){2} & \\
&& *+<8pt>[F]{q}\ar[dl]\ar[dr] & *+<8pt>[F]{r} & \\
\ar@{.}[r] & 1\ar@{.}[rr] && 2\ar@{.}[rr] && }}
\quad\equiv\quad\vcenter{\xymatrix@!R=1em@!C=1em@M=2pt{ 1\ar@/_3pt/[ddr]|(0.69){\hole}|(0.7){1}\ar@{.}[r] & 2\ar[dr]\ar@{.}[rr] && 3\ar[dl]\ar@{.}[r] & 4\ar[dd]|(0.7){2}\ar@{.}[r] & \\
&& *+<8pt>[F]{p}\ar[drr]|(0.75){1}\ar@/^20pt/[dd]\ar@/_8pt/[dd]\ar@/^8pt/[dd] &&& \\
& *+<8pt>[o][F]{1}\ar@/_3pt/[dr]|(0.3){1} &&& *+<8pt>[o][F]{2}\ar[dll]|(0.25){1}\ar@/^12pt/[dll]|(0.25){2} & \\
&& *+<8pt>[F]{q\circ_h r}\ar[dl]\ar[dr] & & \\
\ar@{.}[r] & 1\ar@{.}[rr] && 2\ar@{.}[rr] && }}
\end{equation*}
\caption{A composition relation for graphs}\label{Fig:CompositionRelation}
\end{figure}

Intuitively,
the object $\Env_{\Prop}(\POp)$ is generated by graphs, whose vertices are either numbered or $\POp$-labeled,
divided out by relations which are yielded by the evaluation of composable $\POp$-labels
within graphs.
Figure~\ref{Fig:CompositionRelation} gives an example of such a relation for the labeled tree of Figure~\ref{Fig:PartiallyLabelledGraph},
for labels $p,q,r\in\POp$.

To obtain $\POp\vee\Free(M)\simeq\Env_{\Prop}(\POp)(M)$
from $\Env_{\Prop}(\POp)$,
we simply replace numbered vertices by $M$-labeled vertices
in the graphical description of~$\Env_{\Prop}(\POp)(M)$.
If we forget relations,
then we simply obtain the representative in $\Free(\POp\oplus M)$
of the elements of $\POp\vee\Free(M)$.

\medskip
The next propositions
reflect natural identities
\begin{multline*}
\Free(M)\vee\Free(N)\simeq\Free(M\oplus N),
\quad(\colim_i\POp_i)\vee\Free(M)\simeq\colim_i(\POp_i\vee\Free(M))\\
\text{and}\quad(\POp\vee\Free(M))\vee\Free(N)\simeq\POp\vee\Free(M\oplus N)
\end{multline*}
at the functor level.

\begin{prop}\label{PropFunctor:PropFreeCoproduct}
For a free prop $\Free(M)$, we have a natural isomorpism $\Env_{\Prop}(\Free(M))\simeq F[M]$.
\end{prop}

\begin{proof}
In the case of a free prop $\POp = \Free(M)$,
the morphism
\begin{equation*}
d_0: F[\Free(\Free(M))]\rightarrow F[\Free(M)]
\end{equation*}
of Lemma~\ref{PropFunctor:PropCoproductCoequalizerRepresentation}
is given by the image of the universal composition product $\mu: \Free(\Free(M))\rightarrow\Free(M)$
under the functor $F[-]$;
the morphism $d_1: F[\Free(\Free(M))]\rightarrow F[\Free(M)]$
is given by the shifted universal composition product $\mu: F[\Free(N)]\rightarrow F[N]$
on $N = \Free(M)$.
The universal composition product satisfies a natural associativity relation
which implies that the shifted composition product $\mu: F[\Free(M)]\rightarrow F[M]$
coequalizes $(d_0,d_1)$ and induces a morphism $\psi: \Env_{\Prop}(\Free(M))\rightarrow F[M]$.

In the other direction,
the universal morphism of the free prop $\eta: M\rightarrow\Free(M)$
induces a morphism
\begin{equation*}
\xymatrix{ F[M]\ar[r]^(0.45){F[\eta]} & F[\Free(M)]\ar[r] & \Env_{\Prop}(\Free(M)) }
\end{equation*}
and an easy inspection of definitions shows that this morphism
defines an inverse isomorphism of~$\psi$.
\end{proof}

\begin{prop}\label{PropFunctor:ColimitPropCoproduct}
The functor $\POp\mapsto\Env_{\Prop}(\POp)$ preserves reflexive coequalizers and filtered colimits.
\end{prop}

\begin{proof}
In~\S\ref{PropFunctor:ShiftedObjects}
we observe that a functor of the form $F[-]: M\mapsto F[M]$ preserves reflexive coequalizers and filtered colimits in~$M$
and so does the free prop $\Free: M\mapsto\Free(M)$
and the composite $F[\Free(-)]: M\mapsto F[\Free(M)]$.
The proposition follows immediately by interchange of colimits
since $\Env_{\Prop}(\POp)$ is defined by a reflexive coequalizer of functors of this form.
\end{proof}

\begin{prop}\label{PropFunctor:CoproductPropCoproduct}
For a coproduct $\POp\vee\Free(M)$,
we have a natural isomorphism $\Env_{\Prop}(\POp\vee\Free(M))\simeq\Env_{\Prop}(\POp)[M]$.
\end{prop}

\begin{proof}
By Proposition~\ref{PropFunctor:PropFreeCoproduct} and Proposition~\ref{PropFunctor:ColimitPropCoproduct},
the image of the reflexive coequalizer~(\ref{FreeCoproductPresentation})
under the functor $\Env_{\Prop}(-)$
defines a reflexive coequalizer of the form:
\begin{equation}\label{FreeCoproductShiftedPresentation}
\xymatrix{ \Free[\Free(\POp)\oplus M]\ar@<+2pt>[r]^{d_0}\ar@<-2pt>[r]_{d_1} & \Free[\POp\oplus M]\ar@{.>}[r]\ar@/_8mm/[l]_{s_0} &
\Env_{\Prop}(\POp\vee\Free(M)) }.
\end{equation}
The equivalences of Proposition~\ref{PropFunctor:SumFunctor}
gives an isomorphism of reflexive coequalizers between~(\ref{FreeCoproductShiftedPresentation})
and
\begin{equation}
\xymatrix{ \Free[\Free(\POp)][M]\ar@<+2pt>[r]^{d_0}\ar@<-2pt>[r]_{d_1} & \Free[\POp][M]\ar@{.>}[r]\ar@/_8mm/[l]_{s_0} &
\Env_{\Prop}(\POp)[M] },
\end{equation}
from which the conclusion follows.
\end{proof}

In the next section,
we essentially need the explicit definition of the isomorphism
\begin{equation*}
\xymatrix{ \Env_{\Prop}(\POp)[M]\ar@{.>}[r]_(0.4){\simeq}^(0.4){\rho} & \Env_{\Prop}(\POp\vee\Free(M)) }
\end{equation*}
as a quotient of the natural morphism
\begin{equation*}
\xymatrix{ \Free[\POp][M]\ar[r]^{\simeq} & \Free[\POp\oplus M]\ar@{.}[r] & \Env_{\Prop}(\POp\vee\Free(M)) }
\end{equation*}
In the intuitive representation of~$\Env_{\Prop}(\POp\vee\Free(M))$,
this morphism simply identifies $M$-labels in graphs
with indecomposable $F(M)$-labels.
In a sense,
the isomorphism~$\Env_{\Prop}(\POp\vee\Free(M))\simeq\Env_{\Prop}(\POp)[M]$
maps each element of~$\Env_{\Prop}(\POp\vee\Free(M))$
to a reduced form in which every $F(M)$-label is decomposed
into indecomposables.

\section{Symmetric $\Sigma_*$-multiobjects and decomposition of pushouts}\label{PropPushoutDecomposition}
In this section,
we study the symmetric $\Sigma_*$-multiobject $\Env_{\Prop}(\Free(L)\bigvee_{\Free(K)}\POp)$
associated to a pushout of the form (*)
\begin{equation*}
\xymatrix{ \Free(K)\ar[r]^{u}\ar[d]_{\Free(i)} & \POp\ar@{.>}[d]^{f} \\
\Free(L)\ar@{.>}[r]_(0.3){v} & \Free(L)\bigvee_{\Free(K)}\POp }.
\end{equation*}
Essentially,
we adapt the analysis of~\cite[\S\S 19-20]{FresseModuleBook}
in order to prove that the object $\Env_{\Prop}(\Free(L)\bigvee_{\Free(K)}\POp)$
decomposes into a sequence of manageable pushouts.

Since we only have a good construction
of reflexive coequalizers in the category of props,
we apply:

\begin{fact}\label{PropPushoutDecomposition:PushoutPresentation}
In any category,
a pushout
\begin{equation*}
\xymatrix{ S\ar[r]^{u}\ar[d]_{s} & A\ar@{.>}[d] \\ T\ar@{.>}[r] & B }
\end{equation*}
is equivalent to a reflexive coequalizer
\begin{equation*}
\xymatrix{ T\vee S\vee A\ar@<+1mm>[r]^{d_0}\ar@<-1mm>[r]_{d_1} & T\vee A\ar@/_6mm/[l]_{s_0}\ar[r] & B }
\end{equation*}
such that $d_0 = (\id_T,u,\id_A)$, $d_1 = (\id_T,s,\id_A)$ and $s_0 = (\id_T,\id_A)$.
\end{fact}

From which we obtain:

\begin{prop}\label{PropPushoutDecomposition:PushoutPropCoproduct}
The symmetric $\Sigma_*$-multiobject $\Env_{\Prop}(\Free(L)\bigvee_{\Free(K)}\POp)$
fits in a reflexive coequalizer of the form
\begin{equation*}
\xymatrix{ \Env_{\Prop}(\POp)[K\oplus L]\ar@<+2pt>[r]^{d_0}\ar@<-2pt>[r]_{d_1} &
\Env_{\Prop}(\POp)[L]\ar@{.>}[r]\ar@/_8mm/[l]_{s_0} & \Env_{\Prop}(\Free(L)\bigvee_{\Free(K)}\POp) }.
\end{equation*}
\end{prop}

The morphism
\begin{equation*}
\Env_{\Prop}(\POp)[K\oplus L]\xrightarrow{\Env_{\Prop}(\POp)[(i,\id)]}\Env_{\Prop}(\POp)[L]
\end{equation*}
induced by $i: K\rightarrow L$ gives $d_0$.
The morphism
\begin{equation*}
\Env_{\Prop}(\POp)[K\oplus L]\simeq\Env_{\Prop}(\POp\vee\Free(K))[L]
\xrightarrow{\Env_{\Prop}(\id,u)[L]}\Env_{\Prop}(\POp)[L]
\end{equation*}
determined by $u: K\rightarrow\POp$ gives $d_1$.
The morphism
\begin{equation*}
\Env_{\Prop}(\POp)[L]\xrightarrow{\Env_{\Prop}(\POp)[(0,\id)]}\Env_{\Prop}(\POp)[K\oplus L]
\end{equation*}
induced by the canonical embedding $L\hookrightarrow K\oplus L$
gives $s_0$.

The proposition is a straightforward corollary of Fact~\ref{PropPushoutDecomposition:PushoutPresentation}
and Proposition~\ref{PropFunctor:ColimitPropCoproduct}.

\subsubsection{Sequential decomposition}\label{PropPushoutDecomposition:SequentialDecomposition}
The shifted object $S[M]$
associated to any symmetric $\Sigma_*$-multiobject $S$ has a natural filtration
by subobjects $S[M]_e$ such that:
\begin{multline*}
S[M]_e((a_1,b_1),\dots,(a_r,b_r);(m,n))\\
= \bigoplus_{s=0}^{e}\Bigl\{\bigoplus_{(c_*,d_*)} S((a_*,b_*),(c_*,d_*);(m,n))
\otimes_{\Sigma_{(c_*,d_*)}} M(c_*,d_*)^{\otimes s}\Bigr\}/\Sigma_s.
\end{multline*}
The morphisms of Proposition~\ref{PropPushoutDecomposition:PushoutPropCoproduct}
preserve this filtration.
Hence, we have a sequence of reflexive coequalizers
\begin{equation*}
\xymatrix{ \cdots\ar[r] & \Env_{\Prop}(\POp)[K\oplus L]_{n-1}\ar@<-1mm>[d]\ar@<+1mm>[d]\ar[r] &
\Env_{\Prop}(\POp)[K\oplus L]_{n}\ar@<-1mm>[d]\ar@<+1mm>[d]\ar[r] &
\cdots \\
\cdots\ar[r] & \Env_{\Prop}(\POp)[L]_{n-1}\ar@<-2mm>@/_/[u]\ar@{.>}[d]^{(n-1)}\ar[r] &
\Env_{\Prop}(\POp)[L]_{n}\ar@<-2mm>@/_/[u]\ar@{.>}[d]^{(n)}\ar[r] &
\cdots \\
\cdots\ar@{.>}[r]  & \Env_{\Prop}(\Free(L)\bigvee_{\Free(K)}\POp)_{n-1}\ar@{.>}[r]_{j_r} &
\Env_{\Prop}(\Free(L)\bigvee_{\Free(K)}\POp)_{n}\ar@{.>}[r] & \cdots }
\end{equation*}
such that $\colim_n\Env_{\Prop}(\Free(L)\bigvee_{\Free(K)}\POp)_{n} = \Env_{\Prop}(\Free(L)\bigvee_{\Free(K)}\POp)$.

We rearrange coequalizers $(n-1)$ and $(n)$
to obtain:

\begin{lemm}\label{PropPushoutDecomposition:PushoutSequentialDecomposition}
The morphism $j_n: \Env_{\Prop}(\Free(L)\bigvee_{\Free(K)}\POp)_{n-1}\rightarrow\Env_{\Prop}(\Free(L)\bigvee_{\Free(K)}\POp)_{n}$
fits in a pushout of the form
\begin{equation*}
\xymatrix{  \bigoplus_{\substack{p+q=n\\q<n}}
\{\Env_{\Prop}(\POp)(\cdot)\otimes_{\Sigma_{(\cdot)}} (K(\cdot)^{\otimes p}\otimes L(\cdot)^{\otimes q})\}/{\Sigma_p\times\Sigma_q}
\ar[d]_{\lambda}\ar[r] & \Env_{\Prop}(\Free(L)\bigvee_{\Free(K)}\POp)_{n-1}\ar@{.>}[d]^{j_n} \\
\{\Env_{\Prop}(\POp)(\cdot)\otimes_{\Sigma_{(\cdot)}} L(\cdot)^{\otimes n}\}/\Sigma_n
\ar@{.>}[r] & \Env_{\Prop}(\Free(L)\bigvee_{\Free(K)}\POp)_{n} }.
\end{equation*}
\end{lemm}

\begin{proof}
We reproduce the arguments of~\cite[Lemma 18.2.5]{FresseModuleBook}.
Set
\begin{multline*}
S = \Env_{\Prop}(\POp)[K\oplus L]_{n-1},\quad T = \Env_{\Prop}(\POp)[K\oplus L]_{n-1},\\
\begin{aligned}
\text{and}\quad U & = \bigoplus_{\substack{p+q = n\\q<n}}
\{\Env_{\Prop}(\POp)(\cdot)\otimes_{\Sigma_{(\cdot)}} (K(\cdot)^{\otimes p}\otimes L(\cdot)^{\otimes q})\}/{\Sigma_p\times\Sigma_q} \\
V & = \{\Env_{\Prop}(\POp)(\cdot)\otimes_{\Sigma_{(\cdot)}} (L(\cdot)^{\otimes n})\}/{\Sigma_n}.
\end{aligned}
\end{multline*}
The definition of filtrations
imply
\begin{align*}
\Env_{\Prop}(\POp)[K\oplus L]_{n} & = \Env_{\Prop}(\POp)[K\oplus L]_{n-1}\oplus U\oplus V = S\oplus U\oplus V, \\
\Env_{\Prop}(\POp)[L]_{n} & = \Env_{\Prop}(\POp)[L]_{n-1}\oplus V = T\oplus V
\end{align*}
and the coequalizers $C = \Env_{\Prop}(\Free(L)\bigvee_{\Free(K)}\POp)_{n-1}$ and $D = \Env_{\Prop}(\Free(L)\bigvee_{\Free(K)}\POp)_{n}$
fit in a diagram of the form:
\begin{equation*}
\xymatrix{ S\ar@<-1mm>[d]\ar@<+1mm>[d]\ar[r] & S\oplus U\oplus V\ar@<-1mm>[d]\ar@<+1mm>[d] \\
T\ar@<-2mm>@/_/[u]\ar@{.>}[d]\ar[r] & T\oplus V\ar@<-2mm>@/_/[u]\ar@{.>}[d] \\
C\ar@{.>}[r] & D }.
\end{equation*}
The morphisms $d_0,d_1: S\oplus U\oplus V\rightrightarrows T\oplus V$
are the identity on $V$.
The summand $U$ is mapped into $V$ by the morphism $d_0$,
into $T$ by the morphism $d_1$.
Accordingly,
we have a commutative square
\begin{equation*}
\xymatrix{ U\ar[r]^{d_1}\ar[d]_{d_0} & C\ar@{.>}[d] \\ V\ar@{.>}[r] & D }.
\end{equation*}
An easy inspection
shows that this square forms a pushout.
\end{proof}

\subsubsection{Multifold pushout-products}\label{PropPushoutDecomposition:MultifoldPushoutProducts}
We review the definition of multifold pushout-products. We follow the point of view of~\cite[\S 18.2.6-7]{FresseModuleBook}.
Usually,
we have a morphism $i: K\rightarrow L$ in a symmetric monoidal category,
we set $T_0 = K$ and $T_1 = L$,
and we form the $n$-fold cubical diagram
with tensor products $T_{\epsilon_1}\otimes\dots\otimes T_{\epsilon_n}$
on vertices and morphisms
\begin{equation*}
T_{\epsilon_1}\otimes\dots\otimes T_0\otimes\dots\otimes T_{\epsilon_n}
\xrightarrow{T_{\epsilon_1}\otimes\dots\otimes i\otimes\dots\otimes T_{\epsilon_n}}
T_{\epsilon_1}\otimes\dots\otimes T_1\otimes\dots\otimes T_{\epsilon_n}
\end{equation*}
on edges.
In the context of $\Sigma_*$-biobjects,
we take the tensor product in the base category~$\C$
to obtain a multiple $\Sigma_*$-biobject with a component
\begin{equation*}
T_{\epsilon_1}\otimes\dots\otimes T_{\epsilon_n} = T_{\epsilon_1}(a_1,b_1)\otimes\dots\otimes T_{\epsilon_n}(a_n,b_n)
\end{equation*}
for each $n$-tuple $(a_*,b_*) = (a_1,b_1),\dots,(a_n,b_n)$
of pairs $(a_i,b_i)\in\NN^2$.
For the moment,
we can forget this multiple $\Sigma_*$-bigrading.

The tensor product $T_n(L/K) = T_1\otimes\dots\otimes T_1 = L^{\otimes n}$
is associated to the terminal vertex of the cube.
Set
\begin{equation*}
L_n(L/K) = \colim_{(\epsilon_1,\dots,\epsilon_n)<(1,\dots,1)} T_{\epsilon_1}\otimes\dots\otimes T_{\epsilon_n}.
\end{equation*}
The $n$-fold pushout-product of~$i$
is the natural morphism $\lambda: L_n(L/K)\rightarrow T_n(L/K)$.
This terminology is justified by the following easy observation:

\begin{obsv}[{see~\cite[Observation 18.2.7]{FresseModuleBook}}]\label{PropPushoutDecomposition:IteratedPushoutProduct}
The $n$-fold pushout-product
\begin{equation*}
L_n(L/K)\xrightarrow{\lambda} T_n(L/K) = L^{\otimes n}
\end{equation*}
is identified with the pushout-product
\begin{equation*}
L_{n-1}(L/K)\otimes L\bigoplus_{L_{n-1}(L/K)\otimes K} L^{\otimes n-1}\otimes K
\xrightarrow{(\lambda_*,i_*)} L^{\otimes n-1}\otimes L
\end{equation*}
of the $(n-1)$-fold pushout-product $\lambda: L_{n-1}(L/K)\rightarrow L^{\otimes n-1}$
with $i: K\rightarrow L$.
\end{obsv}

\subsubsection{Multifold pushout-products and functors}\label{PropPushoutDecomposition:FunctorMultifoldPushoutProduct}
The symmetric group $\Sigma_n$ acts naturally (on the right) on $T_n(L/K)$ and $L_n(L/K)$.
Moreover the $n$-fold pushout-product $\lambda: L_n(L/K)\rightarrow T_n(L/K)$
is clearly equivariant.
Recall that the objects $L_n(L/K)$ and $T_n(L/K)$
inherit a natural $n$-fold $\Sigma_*$-bistructure
from the tensor products
\begin{equation*}
T_{\epsilon_1}\otimes\dots\otimes T_{\epsilon_n} = T_{\epsilon_1}(c_1,d_1)\otimes\dots\otimes T_{\epsilon_n}(c_n,d_n).
\end{equation*}
The action of $\Sigma_n$ permutes the summands associated to collections $(c_*,d_*) = ((c_1,d_1),\dots,(c_n,d_n))$.

For a symmetric $\Sigma_*$-multiobject $S$,
we set
\begin{align*}
& L_n S[L/K]((a_*,b_*);(m,n))\\
& \qquad\qquad = \Bigl\{\bigoplus_{(c_*,d_*)} T((a_*,b_*),(c_*,d_*);(m,n))
\otimes_{\Sigma_{(c_*,d_*)}} L_n(L/K)(c_*,d_*)\Bigr\}/\Sigma_n, \\
& T_n S[L/K]((a_*,b_*);(m,n)) \\
& \qquad\qquad = \Bigl\{\bigoplus_{(c_*,d_*)} S((a_*,b_*),(c_*,d_*);(m,n))
\otimes_{\Sigma_{(c_*,d_*)}} T_n(L/K)(c_*,d_*)\Bigr\}/\Sigma_n,
\end{align*}
and we form the morphism of symmetric $\Sigma_*$-multiobjects $\lambda_*: L_n S[L/K]\rightarrow T_n S[L/K]$
induced by $\lambda: L_n(L/K)\rightarrow T_n(L/K)$.

We have by definition
\begin{equation*}
T_n\Env_{\Prop}(\POp)[L/K] =
\{\Env_{\Prop}(\POp)(\cdot)\otimes_{\Sigma_{(\cdot)}} L(\cdot)^{\otimes n}\}/\Sigma_n.
\end{equation*}
We prove by a straightforward inspection of definitions:

\begin{obsv}\label{PropPushoutDecomposition:Basis}
The basis morphisms of the pushout of Lemma~\ref{PropPushoutDecomposition:PushoutSequentialDecomposition}
admit factorizations
\begin{equation*}
\xymatrix@!R=2em@!C=2em@M=2pt{ &&& \bigoplus_{\substack{p+q=n\\q<n}}
\{\Env_{\Prop}(\POp)(\cdot)\otimes_{\Sigma_{(\cdot)}} (K(\cdot)^{\otimes p}\otimes L(\cdot)^{\otimes q})\}/{\Sigma_p\times\Sigma_q}
\ar@/_6pt/[]!DL;[dlll]\ar@/^6pt/[]!DR;[drrr]\ar[d] &&& \\
T_n\Env_{\Prop}(\POp)[L/K] &&& L_n\Env_{\Prop}(\POp)[L/K]\ar[lll]^{\lambda_*}\ar@{.>}[rrr] &&&
\Env_{\Prop}(\Free(L)\bigvee_{\Free(K)}\POp)_{n-1} },
\end{equation*}
where the vertical arrow is induced by the canonical morphisms
\begin{multline*}
L(\cdot)\otimes\dots\otimes K(\cdot)\otimes\dots\otimes L(\cdot)
\rightarrow\colim\{L(\cdot)\otimes\dots\otimes K(\cdot)\otimes\dots\otimes L(\cdot)\}
\end{multline*}
towards the object $L_n[L/K] = \colim\{L(\cdot)\otimes\dots\otimes K(\cdot)\otimes\dots\otimes L(\cdot)\}$.
\end{obsv}

Note that the vertical arrow of this diagram is epi.
Hence the factorizations are unique.
This observation implies moreover:

\begin{prop}\label{PropPushoutDecomposition:Pushouts}
The morphism $j_n: \Env_{\Prop}(\Free(L)\bigvee_{\Free(K)}\POp)_{n-1}\rightarrow\Env_{\Prop}(\Free(L)\bigvee_{\Free(K)}\POp)_{n}$
fits in a pushout:
\begin{equation*}
\xymatrix{ L_n\Env_{\Prop}(\POp)[L/K]\ar[d]_{\lambda_*}\ar[r] & \Env_{\Prop}(\Free(L)\bigvee_{\Free(K)}\POp)_{n-1}\ar@{.>}[d]^{j_n} \\
T_n\Env_{\Prop}(\POp)[L/K]\ar@{.>}[r] & \Env_{\Prop}(\Free(L)\bigvee_{\Free(K)}\POp)_{n} },
\end{equation*}
for every $n\in\NN$.\qed
\end{prop}

\subsubsection{The case of topological spaces}\label{PropPushoutDecomposition:Topology}
In the context of topological spaces,
one checks easily that the morphism $\lambda_*: L_n S(\POp)[L/K]\rightarrow T_n S(\POp)[L/K]$
is a topological inclusion whenever $i: K\rightarrow L$
is a cofibration.
The class of inclusions in topological spaces is closed
under pushouts and composites (see~\cite[Proof of lemma 2.4.5]{Hovey}).
Therefore the decomposition of~\S\ref{PropPushoutDecomposition:PushoutSequentialDecomposition}
and Proposition~\ref{PropPushoutDecomposition:Pushouts}
imply that $j: \Env_{\Prop}(\POp)\rightarrow\Env_{\Prop}(\POp\bigvee_{\Free(K)}\POp)$
is a topological inclusion.
Take the leading term of this morphism to conclude that the prop morphism $j: \POp\rightarrow\POp\bigvee_{\Free(K)}\POp$
is too.
The assertion of Lemma~\ref{PropPathObjectArgument:TopologicalPushouts}
follows from this conclusion.

\section{Homotopy of pushouts along morphisms of free props}\label{PropPushoutHomotopy}
The goal of this section is to achieve the proof of Lemma~\ref{MainResult:PropPushouts}.
Thus we consider again a pushout of the form (*)
\begin{equation*}
\xymatrix{ \Free(K)\ar[r]^{u}\ar[d]_{\Free(i)} & \POp\ar@{.>}[d]^{f} \\
\Free(L)\ar@{.>}[r]_(0.35){v} & \Free(L)\bigvee_{\Free(K)}\POp }
\end{equation*}
in the category of props, but from now on we restrict ourselves to $\Sigma_*$-biobjects and props
with non-empty inputs.
We aim at proving that $f$ defines a cofibration (respectively, an acyclic cofibration) in $\C^{\ACat}_0$
if $i$ is a cofibration (respectively, an acyclic cofibration) in $\C^{\BCat}_0$,
as long as $\POp$ is an $\Free(\C^{\BCat})_c$-cell complex.

As explained in the introduction of this appendix,
we do not study the morphism $f: \POp\rightarrow\Free(L)\bigvee_{\Free(K)}\POp$ directly,
but rather the morphism of symmetric $\Sigma_*$-multiobjects
\begin{equation*}
\Env_{\Prop}(f): \Env_{\Prop}(\POp)\rightarrow\Env_{\Prop}(\Free(L)\bigvee_{\Free(K)}\POp)
\end{equation*}
induced by $f$.

The category of symmetric $\Sigma_*$-multiobjects has a natural model structure (like any category of diagrams in a cofibrantly generated model category).
But the morphism $\Env_{\Prop}(F(0)): \Env_{\Prop}(F(0))\rightarrow\Env_{\Prop}(F(M))$
induced by the initial morphism of a free prop $F(M)$ with $M$ cofibrant
does not belong to the class of natural cofibrations of symmetric $\Sigma_*$-multiobjects.
Therefore our first task is to define a better class of cofibrations.
We call it the class of place-cofibrations.
We prove by induction that the object $\Env_{\Prop}(\POp)$
is place-cofibrant if $\POp$ is an $\Free(\C^{\BCat})_c$-cell complex
and we prove that the pushout (*) returns a morphism $f$
such that $\Env_{\Prop}(f)$ is a place-cofibration (respectively, an acyclic place-cofibration).

To reach the conclusion of Lemma~\ref{MainResult:PropPushouts},
we simply have to take the constant term of~$\Env_{\Prop}(f)$ .

\subsubsection{Place-cofibrations of symmetric $\Sigma_*$-multiobjects}\label{PropPushoutHomotopy:PlaceCofibrations}
To define a good class of cofibrations in the category of symmetric $\Sigma_*$-multiobjects,
we forget internal permutations of inputs
and we focus on $\Sigma_r$-actions
\begin{equation*}
S((a_1,b_1),\dots,(a_r,b_r);(m,n))
\xrightarrow{w^*} S((a_{w(1)},b_{w(1)}),\dots,(a_{w(r)},b_{w(r)});(m,n)).
\end{equation*}

The idea is to replace the category $\BCat$
with its underlying discrete category $\ACat$
in the definition of the category of symmetric $\Sigma_*$-multiobjects $\C^{\Sigma_*\wr\BCat\times\BCat}$
and to apply the restriction functor $\phi^*: \C^{\Sigma_*\wr\BCat\times\BCat}\rightarrow\C^{\Sigma_*\wr\ACat\times\ACat}$.
The category $\Sigma_*\wr\ACat\times\ACat$ is formed by collections of non-negative integers
\begin{equation*}
\alpha = ((a_1,b_1),\dots,(a_r,b_r);(m,n))
\end{equation*}
as objects
and permutations
\begin{multline*}
((a_1,b_1),\dots,(a_r,b_r);(m,n))\xrightarrow{w_*}((a_{w(1)},b_{w(1)}),\dots,(a_{w(r)},b_{w(r)});(m,n)),\\
w\in\Sigma_r,
\end{multline*}
as morphisms.
The restriction functor $\phi^*: \C^{\Sigma_*\wr\BCat\times\BCat}\rightarrow\C^{\Sigma_*\wr\ACat\times\ACat}$
is associated to the obvious embedding $\phi: \Sigma_*\wr\ACat\times\ACat\hookrightarrow\Sigma_*\wr\BCat\times\BCat$.

The categories $\C^{\Sigma_*\wr\BCat\times\BCat}$ and $\C^{\Sigma_*\wr\ACat\times\ACat}$
inherit a natural cofibrantly generated model structure
like every category of diagrams.
Define the class of place-cofibrations (respectively, acyclic place-cofibrations)
as the class of morphisms of symmetric $\Sigma_*$-multiobjects $f: S\rightarrow T$
such that $\phi^*(f)$ is a cofibration (respectively, acyclic cofibration) in $\C^{\Sigma_*\wr\ACat\times\ACat}$.
Note that an acyclic place-cofibration is nothing but a place-cofibration
which defines a pointwise weak equivalence
of symmetric $\Sigma_*$-multiobjects.

The restriction functor $\phi^*: \C^{\Sigma_*\wr\BCat\times\BCat}\rightarrow\C^{\Sigma_*\wr\ACat\times\ACat}$
preserves both limits and colimits,
like the restriction functor $\phi^*: \C^{\BCat}\rightarrow\C^{\ACat}$ on the usual category of $\Sigma_*$-biobjects,
because colimits of diagrams are created componentwise.
From this observation,
we deduce that standard properties of (acyclic) cofibrations in model categories
remain valid for (acyclic) place-cofibrations:

\begin{obsv}\label{PropPushoutHomotopy:PlaceCofibrationsStability}
The class of (acyclic) place-cofibrations is stable under retracts, pushouts and (transfinite) composites.
\end{obsv}

\subsubsection{Restriction of functors to objects with non-empty inputs}\label{PropPushoutHomotopy:NonEmptyInputCase}
Let $\C^{\Sigma_*\wr\BCat\times\BCat}_0$ be the full subcategory of~$\C^{\Sigma_*\wr\BCat\times\BCat}$
formed by symmetric $\Sigma_*$-multiobjects $S$
such that $S((a_1,b_1),\dots,(a_r,b_r);(m,n)) = 0$
if $a_i = 0$ for some $i$
or $m = 0$.
The embedding $\C^{\Sigma_*\wr\BCat\times\BCat}_0\hookrightarrow\C^{\Sigma_*\wr\BCat\times\BCat}$,
like $\C^{\BCat}_0\hookrightarrow\C^{\BCat}$,
has an obvious left adjoint which maps any $S\in\C^{\Sigma_*\wr\BCat\times\BCat}$
to the symmetric $\Sigma_*$-multiobject $\overline{S}\in\C^{\Sigma_*\wr\BCat\times\BCat}_0$
such that:
\begin{equation*}
\overline{S}((a_1,b_1),\dots,(a_r,b_r);(m,n)) = \begin{cases} 0,\quad\text{if $a_i = 0$ for some $i$ or $m = 0$}, & \\
S((a_1,b_1),\dots,(a_r,b_r);(m,n)),\quad\text{otherwise}. & \end{cases}
\end{equation*}
We have an analogous subcategory $\C^{\Sigma_*\wr\ACat\times\ACat}_0\hookrightarrow\C^{\Sigma_*\wr\ACat\times\ACat}$
if we forget internal $\Sigma_*$-actions.

The categories $\C^{\Sigma_*\wr\BCat\times\BCat}_0$ and $\C^{\Sigma_*\wr\ACat\times\ACat}_0$
form sub-model categories of $\C^{\Sigma_*\wr\BCat\times\BCat}$ and $\C^{\Sigma_*\wr\ACat\times\ACat}$ respectively.
Note that the functor $S\mapsto\overline{S}$
preserves colimits, weak equivalences, cofibrations and fibrations.
Besides, for a $\Sigma_*$-biobject with non-empty inputs $M$,
we have relations $S(M) = \overline{S}(M)$ and $\overline{S[M]} = \overline{S}[M]$.
These observations allow us to replace the symmetric $\Sigma_*$-multiobjects of~\S\ref{PropPushoutDecomposition}
by their $\overline{S}$-analogues
when we restrict our study to $\Sigma_*$-biobjects with non-empty inputs.
To simplify the writing,
we adopt the notation $\overline{\Env}_{\Prop}(\POp)$
for the object $\overline{\Env}_{\Prop}(\POp)\in\C^{\Sigma_*\wr\BCat\times\BCat}_0$
associated to $\Env_{\Prop}(\POp)\in\C^{\Sigma_*\wr\BCat\times\BCat}$

\medskip
The next proposition
allows us to perform an induction process
leading to the proof of Lemma~\ref{MainResult:PropPushouts}:

\begin{prop}\label{PropPushoutHomotopy:GeneratorPushout}
Suppose we have a prop with non-empty inputs $\POp$ such that $\overline{\Env}_{\Prop}(\POp)$
forms a place-cofibrant symmetric $\Sigma_*$-multiobject.
Let $i: K\rightarrow L$ be a generating cofibration (respectively, an acyclic generating cofibration)
of the category of $\Sigma_*$-biobjects with non-empty inputs.
The morphism $\overline{\Env}_{\Prop}(f): \overline{\Env}_{\Prop}(\POp)\rightarrow\overline{\Env}_{\Prop}(\Free(L)\bigvee_{\Free(K)}\POp)$
yielded by the pushout (*)
forms a place-cofibration (respectively, an acyclic place-cofibration).
\end{prop}

The proof of this proposition is deferred to a series of observations.
Our plan is to study each piece of the decomposition of~\S\ref{PropPushoutDecomposition}.
To begin with:

\begin{lemm}\label{PropPushoutHomotopy:MultifoldPushoutProduct}
Let $S$ be a place-cofibrant symmetric $\Sigma_*$-multiobject such that $S = \overline{S}$.
If $i: K\rightarrow L$ is a generating cofibration (respectively, an acyclic generating cofibration) of $\Sigma_*$-biobjects (with non-empty inputs),
then the morphism $\lambda_*: L_n S[L/K]\rightarrow T_n S[L/K]$
induced by the $n$-fold pushout-product of~$i$
forms a place-cofibration (respectively, an acyclic place-cofibration) of symmetric $\Sigma_*$-multiobjects.
\end{lemm}

\begin{proof}
By definition of generating (acyclic) cofibrations of $\Sigma_*$-biobjects,
we assume that $i$ is a morphism of the form
\begin{equation*}
u\otimes\phi_! G_{p q}: C\otimes\phi_! G_{p q}\rightarrow D\otimes\phi_! G_{p q}
\end{equation*}
where $u$ is a generating (acyclic) cofibration of~$\C$
and $G_{p q}$ is the enriched Yoneda $\A$-diagram
associated to $(p,q)\in\NN\times\NN$.
To get a generating (acyclic) cofibration of $\Sigma_*$-biobjects with non-empty inputs,
we simply have to restrict ourselves to pairs $(p,q)$
such that $p>0$.

By an easy inspection,
we see that the $n$-fold pushout-product of~$i$
is identified with the tensor product
\begin{equation*}
\lambda\otimes\phi_! G_{p q}^{\otimes n}: L_n(D/C)\otimes\phi_!  G_{p q}^{\otimes n}\rightarrow T_n(D/C)\otimes\phi_! G_{p q}^{\otimes n}
\end{equation*}
where $\lambda: L_n(D/C)\rightarrow T_n(D/C)$
refers to the $n$-fold pushout-product of~$u: C\rightarrow D$
in the symmetric monoidal category~$\C$.
In this expression,
we use the multiple $\Sigma_*$-biobject $\phi_! G_{p q}^{\otimes n}$
which has a component
\begin{equation*}
\phi_! G_{p q}(c_*,d_*)^{\otimes n} = \phi_! G_{p q}(c_1,d_1)\otimes\dots\otimes\phi_! G_{p q}(c_n,d_n)
\end{equation*}
for each $n$-tuple $(c_*,d_*) = (c_1,d_1),\dots,(c_n,d_n)$
of pairs $(c_j,d_j)\in\NN^2$.

We have
by definition of enriched Yoneda diagrams:
\begin{equation*}
\phi_! G_{p q}(c,d) = \begin{cases} \unit[\Sigma^{op}_p\times\Sigma_q], & \text{if $c = p$ and $d = q$}, \\
0, & \text{otherwise}. \end{cases}
\end{equation*}
From this definition,
we deduce the identity
\begin{multline*}
\bigl\{\bigoplus_{(c_*,d_*)} S((a_*,b_*),(c_*,d_*);(-,-))
\otimes_{\Sigma_{(c_*,d_*)}}
T_{\epsilon_*}^{\otimes n}\otimes G_{p q}(c_*,d_*)^{\otimes n}\bigr\}/{\Sigma_n}\\
\simeq\bigl\{S((a_*,b_*),(p_*,q_*);(-,-))\otimes T_{\epsilon_*}^{\otimes n}\bigr\}/{\Sigma_n}
\end{multline*}
for each tensor product $T_{\epsilon_*}^{\otimes n} = T_{\epsilon_1}\otimes\dots\otimes T_{\epsilon_r}\in\C$,
where we set $(p_*,q_*) = (p,q),\dots,(p,q)$.
Consequently,
the morphism $\lambda_*: L_n S[L/K]\rightarrow T_n S[L/K]$
can be identified with the morphism
\begin{multline}\label{NonsymmetricPushoutProduct}
\bigl\{S((a_*,b_*),(p_*,q_*);(-,-))\otimes L_n(D/C)\bigr\}/{\Sigma_n}\\
\xrightarrow{\lambda_*}\bigl\{S((a_*,b_*),(p_*,q_*);(-,-))\otimes T_n(D/C)\bigr\}/{\Sigma_n}
\end{multline}
induced by $\lambda: L_n(D/C)\rightarrow T_n(D/C)$,
the $n$-fold pushout-product of~$u$ in $\C$.

Internal $\Sigma_{(c_*,d_*)}$-actions
do not occur anymore in formulas,
and can now be forgotten.
Hence we can see $S$ as an object of the category $\C^{\Sigma_*\wr\ACat\times\ACat}_0$.

To simplify notation,
we apply the conventions of~\S\ref{PropPushoutDecomposition:FunctorMultifoldPushoutProduct}
in the context of objects without internal $\Sigma_{(\cdot)}$-actions
and we adopt the notation $\lambda_*: L_n S[D/C]\rightarrow T_n S[D/C]$
to refer to~(\ref{NonsymmetricPushoutProduct})
when $S$ is supposed to belong to~$\C^{\Sigma_*\wr\ACat\times\ACat}$
or is viewed as an object of~$\C^{\Sigma_*\wr\ACat\times\ACat}$.
Note that the morphism $\lambda_*: L_n S[D/C]\rightarrow T_n S[D/C]$
is natural in $S$.
Consequently,
for any morphism $\rho: R\rightarrow S$ in $\C^{\Sigma_*\wr\ACat\times\ACat}_0$,
we have a pushout-product
\begin{equation*}
T_n R[D/C]\bigoplus_{L_n R[D/C]} L_n S[D/C]\xrightarrow{(\rho_*,\lambda_*)} T_n S[D/C]
\end{equation*}
naturally associated to $\rho$.
We prove that this pushout-product forms a cofibration in~$\C^{\Sigma_*\wr\ACat\times\ACat}$
if $\rho$ does so and $i$ is a cofibration,
an acyclic cofibration if $\rho$ or $i$ is also acyclic.
We apply this result to the initial morphism $0: 0\rightarrow S$
to conclude.

By the argument of~\cite[Lemma 4.2.4]{Hovey},
it is sufficient to address the case where $\rho$ is a generating (acyclic) cofibration of~$\C^{\Sigma_*\wr\ACat\times\ACat}_0$.
In this case,
our claim is proved by a straightforward inspection
using the form of generating (acyclic) cofibrations
given in Proposition~\ref{PropSemiModel:DiagramModelCategory}.
\end{proof}

The next assertions
are immediate consequences of the stability claim of observation~\ref{PropPushoutHomotopy:PlaceCofibrationsStability}:

\begin{fact}\label{PropPushoutHomotopy:PushoutStability}
If $\lambda_*: L_n\overline{\Env}_{\Prop}(\POp)[L/K]\rightarrow T_n\overline{\Env}_{\Prop}(\POp)[L/K]$
is a place-cofibration (respectively, an acyclic place-cofibration),
then so is the morphism $j_n$ yielded by the pushout of Proposition~\ref{PropPushoutDecomposition:Pushouts}:
\begin{equation*}
\xymatrix{ L_n\overline{\Env}_{\Prop}(\POp)[L/K]\ar[d]_{\lambda}\ar[r] & \overline{\Env}_{\Prop}(\Free(L)\bigvee_{\Free(K)}\POp)_{n-1}\ar@{.>}[d]^{j_n} \\
T_n\overline{\Env}_{\Prop}(\POp)[L/K]\ar@{.>}[r] & \overline{\Env}_{\Prop}(\Free(L)\bigvee_{\Free(K)}\POp)_n }.
\end{equation*}
\end{fact}

\begin{fact}\label{PropPushoutHomotopy:CompositeStability}
If any $j_n: \overline{\Env}_{\Prop}(\Free(L)\bigvee_{\Free(K)}\POp)_{n-1}\rightarrow\overline{\Env}_{\Prop}(\Free(L)\bigvee_{\Free(K)}\POp)_{n}$
is a place-cofibration (respectively, an acyclic place-cofibration),
then so is the composite of Proposition~\ref{PropPushoutDecomposition:PushoutSequentialDecomposition}:
\begin{multline*}
\overline{\Env}_{\Prop}(\POp) = \overline{\Env}_{\Prop}(\POp)_0\rightarrow\dots\\
\dots\rightarrow\overline{\Env}_{\Prop}(\Free(L)\bigvee_{\Free(K)}\POp)_{n-1}\rightarrow\overline{\Env}_{\Prop}(\Free(L)\bigvee_{\Free(K)}\POp)_{n}\rightarrow\dots\\
\dots\rightarrow\colim_n\Bigl\{\overline{\Env}_{\Prop}(\Free(L)\bigvee_{\Free(K)}\POp)_n\Bigr\} = \overline{\Env}_{\Prop}(\Free(L)\bigvee_{\Free(K)}\POp).
\end{multline*}
\end{fact}

This assertion achieves the proof of Proposition~\ref{PropPushoutHomotopy:GeneratorPushout}
since the composite of Fact~\ref{PropPushoutHomotopy:CompositeStability}
arises from a decomposition of~$\overline{\Env}_{\Prop}(f)$.\qed

\begin{prop}\label{PropPushoutHomotopy:GeneralCofibrationPushout}
The conclusion of Proposition~\ref{PropPushoutHomotopy:GeneratorPushout}
holds when $i: K\rightarrow L$ is any (acyclic) cofibration of~$\Sigma_*$-biobjects
and not only a generating (acyclic) cofibration.
\end{prop}

\begin{proof}
Let $\I$ (respectively, $\J$) be the set of generating (respectively, acyclic) cofibrations of $\C^{\BCat}$.
Use that:
\begin{itemize}
\item
any (respectively, acyclic) cofibration forms a retract of a relative $\I$-cell (respectively, $\J$-cell) complex,
\item
pushouts along coproducts decompose naturally into composites of pushouts,
\item
pushouts preserve retracts, pushouts and composites,
\item
the class of (acyclic) cofibrations is stable under retracts, pushouts and composites
\end{itemize}
to conclude.
\end{proof}

The next lemma allows us to initiate an inductive application of Proposition~\ref{PropPushoutHomotopy:GeneralCofibrationPushout}
to cell complexes in the category of props:

\begin{lemm}\label{PropPushoutHomotopy:InitialProp}
For the initial prop $I$,
we have an identity $\Env_{\Prop}(I) = F$,
where $F$ is the symmetric $\Sigma_*$-multiobject of~\S\ref{PropFunctor:FreePropFunctor}
which represents the free prop $\Free: \C^{\BCat}\rightarrow\Prop$.

The associated reduced symmetric $\Sigma_*$-multiobject $\overline{\Env}_{\Prop}(I) = \overline{F}$
is place-cofibrant.
\end{lemm}

\begin{proof}
The identity $\Env_{\Prop}(I) = F$ is a particular case of the result of Proposition~\ref{PropFunctor:PropFreeCoproduct},
because we have $I = \Free(0)$
and $F[0] = F$.

The last assertion is a consequence of the analysis of Proposition~\ref{PropFunctor:FreePropSymmetricStructure}.
\end{proof}

By an immediate induction from the result of Lemma~\ref{PropPushoutHomotopy:InitialProp},
we obtain:

\begin{lemm}\label{PropPushoutHomotopy:Object}
The symmetric $\Sigma_*$-multiobject $\overline{\Env}_{\Prop}(\POp)$
associated to an $\Free(\C^{\BCat})_c$-cell complex in the category of props with non-empty inputs
is place-cofibrant.\qed
\end{lemm}

Then Proposition~\ref{PropPushoutHomotopy:GeneralCofibrationPushout}
gives as an immediate corollary:

\begin{prop}\label{PropPushoutHomotopy:CoefficientConclusion}
Let $\POp$ be an $\Free(\C^{\BCat})_c$-cell complex in the category of props with non-empty inputs.
Let $i: K\rightarrow L$ be a cofibration (respectively, an acyclic cofibration)
of the category of $\Sigma_*$-biobjects with non-empty inputs.
The morphism $\overline{\Env}_{\Prop}(f): \overline{\Env}_{\Prop}(\POp)\rightarrow\overline{\Env}_{\Prop}(\Free(L)\bigvee_{\Free(K)}\POp)$,
yielded by the pushout (*)
forms a place-cofibration (respectively, an acyclic place-cofibration).\qed
\end{prop}

From which we conclude:

\begin{prop}\label{PropPushoutHomotopy:Conclusion}
Let $\POp$ be an $\Free(\C^{\BCat})_c$-cell complex in the category of props with non-empty inputs.
If $i: K\rightarrow L$ forms a cofibration (respectively, an acyclic cofibration) in the category $\Sigma_*$-biobjects with non-empty inputs,
then the prop morphism $f: \POp\rightarrow\Free(L)\bigvee_{\Free(K)}\POp$
yielded by the pushout (*)
forms a cofibration in~$\C^{\ACat}$.
\end{prop}

\begin{proof}
Corollary of Proposition~\ref{PropPushoutHomotopy:CoefficientConclusion}.
Note simply that the leading term $f: R(\emptyset;(m,n))\rightarrow S(\emptyset;(m,n))$
of a place-cofibration (respectively, an acyclic place-cofibration)
of symmetric $\Sigma_*$-multiobjects forms a cofibrations in~$\C^{\ACat}$.
\end{proof}

This assertion achieves the proof of Lemma~\ref{MainResult:PropPushouts}
and Theorem~\ref{PropSemiModel:Result}.\qed

\subsubsection{Remark}\label{PropPushoutHomotopy:CanonicalCoproductForm}
In the case of operads and their variants (cyclic operads, modular operads, properads, \dots)
a pushout of the form $\QOp = \POp\vee\Free(M)$
can be described as a colimit of tensor products over graphs
\begin{equation*}
\Bigl\{\bigotimes_{v\in V_{\POp}(\Gamma)} \POp(\In_h,\Out_h)\Bigr\}
\otimes\Bigl\{\bigotimes_{v\in V_{M}(\Gamma)} M(\In_h,\Out_h)\Bigr\},
\end{equation*}
where $V_{\POp}(\Gamma)\amalg V_{M}(\Gamma) = V(\Gamma)$
is a partition of the vertices of $\Gamma$,
such that no pair of vertices $v,w\in V_{\POp}(\Gamma)$
is adjacent in $\Gamma$.
The quotient relation of the colimit does not involve the multiplicative structure of $\POp$,
because only adjacent vertices can be composed in operad-like structures.
The improvement of Theorem~\ref{PropSemiModel:Result} alluded to in~\S\ref{PropSemiModel:Difficulties}
can be obtained from this observation,
as we can forget composition structures and study the homotopy of pushouts
from axioms of symmetric monoidal model categories.

The horizontal composition operation $\circ_h: \POp(k,m)\otimes\POp(l,n)\rightarrow\POp(k+l,m+n)$
makes the case of props more complicated,
because it gives interactions between non-adjacent vertices
in graphical descriptions of coproducts.
For instance,
consider the graph
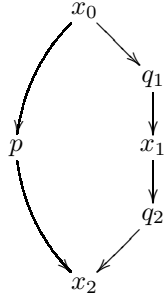
\begin{equation*}
\xymatrix@!C=2pt@!R=2pt@M=2pt{ & x_0\ar[dr]\ar@/_6pt/[ddl] & \\
&& q_1\ar[d] \\
p\ar@/_6pt/[ddr] && x_1\ar[d] \\
&& q_2\ar[dl] \\
& x_2 & }
\end{equation*}
which represents a formal composite of elements $x_0,x_1,x_2\in M$, $p,q_1,q_2\in\POp$.
In $\POp\vee\Free(M)$,
we have the relations
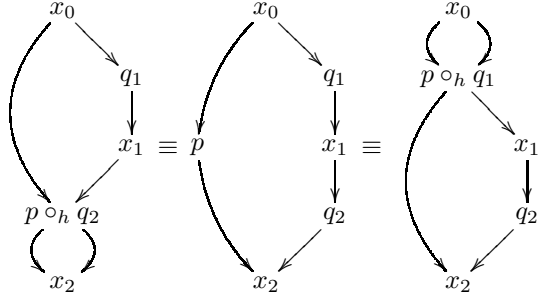
\begin{equation*}
\vcenter{\xymatrix@!C=2pt@!R=2pt@M=2pt{ x_0\ar[dr]\ar@/_2em/[ddd] & \\
& q_1\ar[d] \\
& x_1\ar[dl] \\
p\circ_h q_2\ar@/^12pt/[d]\ar@/_12pt/[d] & \\
x_2 & }}
\equiv\vcenter{\xymatrix@!C=2pt@!R=2pt@M=2pt{ & x_0\ar[dr]\ar@/_6pt/[ddl] & \\
&& q_1\ar[d] \\
p\ar@/_6pt/[ddr] && x_1\ar[d] \\
&& q_2\ar[dl] \\
& x_2 & }}
\equiv\quad\vcenter{\xymatrix@!C=2pt@!R=2pt@M=2pt{ x_0\ar@/^12pt/[d]\ar@/_12pt/[d] & \\
p\circ_h q_1\ar@/_2em/[ddd]\ar[dr] & \\
& x_1\ar[d] \\
& q_2\ar[dl] \\
x_2 & }}
\end{equation*}
and the graphs on the left and right hand side cannot be reduced further (this would create a directed loop in the graph).
From these identities,
we see that the elements of $\POp\vee\Free(M)$
have no canonical form and any realization of $\POp\vee\Free(M)$
necessarily involves a quotient by relations involving the horizontal product of~$\POp$.

\end{document}